\documentclass[12pt,a4paper]{article}

\usepackage[english]{babel}
\usepackage{newtxtext}
\usepackage[varvw]{newtxmath}
\usepackage[scaled=0.90]{helvet} % ss
\usepackage{courier} % tt
\normalfont
\usepackage[T1]{fontenc}
\usepackage{latexsym,amsmath}
\usepackage{makeidx}
\usepackage{a4}
\usepackage{dsfont}
\newtheorem{Thm}{Theorem}
\newtheorem{Defi}[Thm]{Definition}
\newtheorem{Cor}[Thm]{Corollary}
\newtheorem{Lemma}[Thm]{Lemma}
\newtheorem{Prop}[Thm]{Proposition}
\newtheorem{Rem}[Thm]{Remark}
\newtheorem{Conj}[Thm]{Conjecture}
\newtheorem{Prelim}[Thm]{Preliminary}

\newenvironment{thm}[0]{\begin{Thm}\noindent}%
{\end{Thm}}
\newenvironment{defi}[0]{\begin{Defi}\noindent\rm}%
{\end{Defi}}
\newenvironment{cor}[0]{\begin{Cor}\noindent}%
{\end{Cor}}
\newenvironment{lemma}[0]{\begin{Lemma}\noindent}%
{\end{Lemma}}
\newenvironment{prop}[0]{\begin{Prop}\noindent}%
{\end{Prop}}
\newenvironment{rem}[0]{\begin{Rem}\noindent\rm}%
{\end{Rem}}
{\end{Conj}}
{\end{Prelim}}

\def\proof{\par\noindent{\it Proof.}{\ }{\ }}
\def\qed{~\hfill\hbox{$\Box$}\medbreak}

\def\naam#1{\label{#1}}

\def\medno{\medbreak\noindent}
\def\text#1{\;\;\;\;{\rm \hbox{#1}}\;\;\;\;}
\def\qquad{\quad\quad}
\def\qqquad{\quad\quad\quad}

\def\itema{\vspace{-1mm}\item[{\rm (a)}]}
\def\itemb{\item[{\rm (b)}]}
\def\itemc{\item[{\rm (c)}]}
\def\itemd{\item[{\rm (d)}]}

\def\msy#1{{\mathbb #1}}
\def\C{{\msy C}}
\def\N{{\msy N}}
\def\Z{{\msy Z}}
\def\R{{\msy R}}

\def\ga{\alpha}
\def\gb{\beta}
\def\gd{\delta}
\def\ge{\varepsilon}
\def\gf{\varphi}
\def\gg{\gamma}

\def\gl{\lambda}

\def\gs{\sigma}

\def\gD{\Delta}

\def\gS{\Sigma}

\def\got#1{\mathfrak #1}
\def\fa{{\got a}}
\def\fb{{\got b}}
\def\fg{{\got g}}
\def\fh{{\got h}}

\def\fk{{\got k}}

\def\fm{{\got m}}
\def\fn{{\got n}}

\def\fp{{\got p}}
\def\fq{{\got q}}
\def\fs{{\got s}}
\def\ft{{\got t}}

\def\implies{\Rightarrow}
\def\to{\rightarrow}
\def\Re{{\rm Re}\,}
\def\Im{{\rm Im}\,}
\def\inp#1#2{\langle#1\,,\,#2\rangle}

\def\Ad{{\rm Ad}}
\def\End{{\rm End}}
\def\Hom{{\rm Hom}}

\def\ad{{\rm ad}}

\def\after{\,{\scriptstyle\circ}\,}

\def\implies{\Leftarrow}

\def\iC{{\scriptscriptstyle \C}}

\def\cA{{\mathcal A}}

\def\cC{{\mathcal C}}
\def\cD{{\mathcal D}}
\def\cE{{\mathcal E}}
\def\cF{{\mathcal F}}

\def\cK{{\mathcal K}}

\def\cO{{\mathcal O}}
\def\cP{{\mathcal P}}
\def\cR{{\mathcal R}}
\def\cS{{\mathcal S}}

\def\bs{{\backslash}}

\numberwithin{Thm}{section}
\numberwithin{equation}{section}

\def\rdiff{{\rm\hspace{1pt};\hspace{1pt}}}

\def\Ind{{\rm Ind}}
\def\supp{{\rm supp}}
\def\dotvar{\, \cdot\,}
\def\Vtau{V_\tau}

\def\fZ{{\mathfrak Z}}

\def\implies{\Rightarrow}
\def\Wh{{\rm Wh}}
\def\ev{{\rm ev}}

\def\faPdc{{\fa_{P\iC}^*}}
\def\faPd{{\fa_P^*}}

\def\fhRdc{{\fh_{R\iC}^*}}

\def\ds{{\rm ds}}
\def\gL{\Lambda}

\def\Cartan{\theta}

\numberwithin{equation}{section}

\def\fad{{\fa^*}}
\def\cl{{\rm cl}}
\def\fZ{{\mathfrak Z}}

\def\gL{\Lambda}
\def\ft{\mathfrak t}

\def\bp{{\,}^\backprime}
\def\spec{{\rm spec}}

\def\faPdc{\fa_{P\iC}^*}
\def\fh{{\mathfrak h}}
\def\Cartan{\theta}

\def\faP{\fa_P}
\def\faQ{\fa_Q}

\def\twohol{{\rm II}_{\rm hol}}
\def\twoholpr{{\rm II}'_{\rm hol}}

\def\faQd{{\fa_Q^*}}
\def\faQdc{\fa_{Q\iC}^*}

\def\cPst{{\cP_{\rm st}}}

\def\faQp{\fa_Q^+} 
\def\faQ{{\fa_Q}}

\def\cPst{\cP_{\rm st}}
\def\faR{{\fa_R}}
\def\faP{{\fa_P}}
\def\faRd{\fa_R^*}
\def\faRdc{\fa_{R \iC}^*}
\def\st{{}^*\!}
\def\ist{{\rm st}}

\def\dMPds{{\widehat M}_{P,\ds}}
\def\dMQds{{\widehat M}_{Q,\ds}}
\def\dMRds{{\widehat M}_{R,\ds}}

\def\faQp{\fa_Q^+}

\def\cJ{\mathcal{J}}

\def\iQ{{\scriptscriptstyle Q}}

\def\faPp{\fa_P^+}

\def\cAtwoP{\cA_{2, P}}
\def\Wave{\cJ}
\def\fz{\mathfrak{z}}

 \def\Tup{T^\uparrow}
 \def\Tdown{T^\downarrow}
 \def\ccH{H}

 \def\Kcirc{K_\circ}
 \def\oppcPst{\bar{\cP}_{\rm st}}
 \def\od{\circ}
  \def\rmm{{\rm m}}
 \def\co{{\rm co}}
 \def\nFou{\cF^\circ}
 \def\Eta{\eta_*}

\begin{document}
%%% ----------------------------------------------------------------------
\title{Maass--Selberg relations for Whittaker functions\\on a real reductive group}
\author{Erik P. van den Ban}
%\date{17 November 2025}
\maketitle
\tableofcontents
\setcounter{section}{-1}
\section{Introduction}
In this paper we give a complete proof of the Maass--Selberg relations for Whittaker integrals on a 
real reductive Lie group. These relations were announced in 1982 by Harish-Chandra as an important
part of the Plancherel formula for Whittaker functions. Because of his untimely death in 1983, no further details appeared
until 2018, when V.S. Varadarajan and R. Gangolli published an edited version of Harish-Chandra's manuscripts,  \cite[pp. 141-307]{HCwhit}.
The results of the present paper are based on a result of  \cite{HCwhit} for a basic case. Harish-Chandra's proof for this basic case involves an asymptotic analysis of boundary terms reflecting the non-symmetry of the Casimir operator over an expanding  $K$-invariant domain in $G/N_0$ whose radial part  is a  simplex. In particular, this involves the application of Gauss' divergence theorem on a simplex. 
 A detailed account of Harish-Chandra's arguments is given in Sections \ref{s: radial part Casimir} and \ref{s: result HC} of the present paper. 

Starting with the proof for the basic case, there appears to be a complete proof of the Maass-Selberg relations in  \cite{HCwhit}; however, we have not been able to understand the details.
In the present paper we follow a different approach by combining 
the result for the basic case with ideas from 
the theory of reductive symmetric spaces, in which the action of the so-called  standard  intertwining operators plays a central role. We believe the obtained information is 
of separate interest. As an application of the Maas-Selberg relations we prove that the normalized Fourier and Wave packet transforms are continuous linear maps between appropriate Schwartz spaces. 

N.~Wallach  \cite{Wwhit} independently developed another approach to the Whittaker--Plancherel formula, in which the Whittaker Maass-Selberg relations do not seem to play a role. 

The results of the present paper will be of key importance in a proof of the Plancherel
theorem that I have outlined in several lectures in recent years. The accompagnying slides
are available on my website. Details will appear in a follow up paper. 

We will now describe the results of our paper in more detail. It is assumed that $G$ is a real reductive Lie group
of the Harish-Chandra class, that $K$ is a maximal compact subgroup, and that $G = K A N_0$ is an
Iwasawa decomposition. Furthermore, $\chi$ is a fixed unitary character of $N_0$ which is regular in the sense
that for any simple root $\ga$ of  $\fa$ in $\fn_0$ the restriction of $\chi_*:= d\chi(e)$ 
to the root space $\fg_\ga$ is non-zero. Here $\fa$ and $\fn_0$ denote the Lie algebras of $A$ and $N_0,$ respectively, in accordance with the convention to denote
Lie groups by roman capitals and their 
Lie algebras by the corresponding gothic lower cases.

The root system of $\fa$ in $\fg$ is denoted by $\gS.$ Furthermore, $\gS^+$ denotes the positive system 
for which $\fn_0$ is the sum of the associated root spaces $\fg_\ga, $ for $\ga \in \gS^+.$ 

We denote by $\cP$ the finite set of parabolic subgroups of $G$
containing $A$ and by $\cPst$ the subset of the  standard  ones among them, i.e., the $P \in \cP$ such that $P$ contains 
the minimal parabolic subgroup $P_0:= Z_G(\fa) N_0.$ 
Every $Q \in \cP$ is conjugate to a unique $P \in \cPst$, under the action of $N_K(\fa).$ The action of the latter group on $\fa$ induces an isomorphism 
from $N_K(\fa)/Z_K(\fa)$ onto $W(A),$ the Weyl group of the root sytem $\gS.$ 

Given $Q \in \cP$ we denote its Langlands decomposition by 
$Q = M_Q A_Q N_Q.$ By $\dMQds$ we denote the set of equivalence classes
of representations in the discrete series of $M_Q.$ 

For $Q \in \cP$, $\gs\in \dMQds$ and $\nu \in \faQdc$ 
we define $C^\infty(G/ Q: \gs :\nu)$ to be the Fr\'echet space of smooth functions $\gf: G \to H_\gs$ transforming 
according to the rule 
$$
\gf(x m an ) = a^{-\gl - \rho_Q} \gs(m)^{-1} \gf(x),\qquad (x \in G, (m,a ,n) \in M_Q \times A_Q \times N_Q).
$$
Equipped with the left regular representation this space 
realizes the space of smooth vectors for the normalized induced representation 
\begin{equation}
\label{e: normalized induced rep}
\Ind_{Q}^G(\gs \otimes \nu \otimes 1).
\end{equation}
Let $C^{-\infty}(G/ Q:\gs :\nu)$ denote the 
continuous conjugate linear dual of the Fr\'echet space $C^{\infty}(G/Q:\gs :-\bar \nu).$  Via the  standard  $G$-equivariant
sesquilinear pairing by integration over $K,$ 
$C^\infty(G/Q:\gs:\nu)$ is injectively and  $G$-equivariantly mapped into $C^{-\infty}(G/ Q:\gs:\nu).$ Accordingly, 
the latter is viewed as the space of generalized vectors for (\ref{e: normalized induced rep}).
We write 
\begin{equation}
\label{e: type chi}
C^{-\infty}(G/ Q:\gs : \nu)_\chi 
\end{equation}
for the subspace of $C^{-\infty}(G/ Q:\gs :\nu)$ consisting of $\gf$ such that
$\gf(n x) = \chi(n) \gf(x),$ for $ x\in G, n \in N_0$.
Its elements are called the generalized Whittaker vectors of type $\chi.$ 

If $P$ is an opposite  standard  parabolic subgroup then $N_0 P$ is open in $G.$  In  \cite[Thm. 8.6]{Bunitemp} it is shown that every function $\gf  \in C^{-\infty}(G/ P:\gs :\nu)_\chi$ 
restricts to a continuous function 
$N_0 P \to H_\gs^{-\infty},$ satisfying $\gf (n m) = \chi(n) \gs(m)^{-1}  \gf(e)$ 
for $n\in N_0$ and $m \in M_P.$ In particular, $\gf(e) \in H_\gs^{-\infty}$ satisfies 
$$ 
\gs(n) \gf(e) = \chi(n) \gf(e), \qquad (n \in  M_P \cap N_0 ).
$$ 
We define $\chi_P:= \chi|_{(M_P \cap N_0)}$ and put
$$ 
\ccH^{-\infty}_{\gs, \chi_P} = \{ \eta \in H_\gs^{-\infty} \mid \forall n \in  M_P \cap N_0: 
\gs(n) \eta = \chi_P(n) \eta\;\}.
$$ 
In  \cite[Prop.~8.15]{Bunitemp} it was proven that the evaluation map $\ev_e: C^{-\infty}(G/P : \gs : \nu)_\chi \to \ccH_{\gs, \chi_P}^{-\infty}$ is a bijective linear map of finite dimensional linear spaces, for every $\nu \in \faPdc.$ 
The inverse of $\ev_e$ is denoted by 
\begin{equation}
\label{e: j for opp standard}
j(P : \gs: \nu):  \ccH_{\gs,\chi_P}^{-\infty}  \to C^{-\infty}(G/P: \gs: \nu)_\chi.
\end{equation}
Furthermore, according to  \cite[Prop.~8.14]{Bunitemp}, 
for every $\eta \in\ccH_{\gs,\chi_P}^{-\infty}$ the function $\nu \to j(P, \gs, \nu)\eta$ 
is holomorphic as a function on $\faPdc$ with values in $C^{-\infty}(K/K_P: \gs_P)$ (in the compact picture).

In  \cite[Prop. 8.10]{Bunitemp}  the Whittaker integrals for $P \in \cPst$ are essentially defined as finite
sums of matrix coefficients of $K$-spherical matrix coefficients with the generalized Whittaker vectors for 
$\Ind_{\bar P}^G(\gs \otimes \bar \nu \otimes 1).$
The Maass--Selberg relations give information about their asymptotic behavior towards infinity,
on the closed positive Weyl chamber $\cl A^+$. In the  theory of reductive symmetric spaces
these Maass--Selberg relations can be reformulated in terms of the action of  standard  intertwinining
operators on the analogues of the generalized vectors $j(\bar P:\gs:\bar \nu)\eta.$ The known product 
decomposition of these operators then reduce the Maass-Selberg relations to a basic case, where
they can be established more directly. For this approach to work in the Whittaker setting, one needs to 
define the generalized Whittaker vectors for $\Ind_{Q}^G(\gs \otimes \nu \otimes 1),$ 
with $Q \in \cP$ not necessarily opposite  standard. This is worked out in Section \ref{s: defi j}, making use of the existence
of an element  $v \in  N_K(\fa)$ such that $vQv^{-1}$ is opposite  standard.
The orbit $N_0 v \bar Q$ is the unique open $N_0$-orbit on $G/Q.$ In this setting, evaluation at $v$ defines
a bijective linear map $\ev_v : C^{-\infty}(G/Q: \gs : \nu)_\chi \to \ccH^{-\infty}_{\gs, \chi_Q},$ whose inverse
$j(Q: \gs: \nu)$ gives the appropriate generalization of $j(P ,\gs , \nu),$ $P\in \oppcPst.$ 
At the end of Section \ref{s: defi j} we give the definition of the corresponding Whittaker integral as a finite sum of matrix coefficients of
$K$-spherical vectors with the Whittaker vectors of the induced representations $\Ind_{\bar Q}^G(\gs \otimes -\nu \otimes 1).$

In Section \ref{s: interaction with W} it is shown that for parabolic subgroups $P,Q\in \cP$ with equal split components (i.e., $\faP = \faQ$) there is a unique meromorphic map
$\nu \mapsto B(Q , P, \gs, \nu),$ $\faPdc \to \Hom(\ccH^{-\infty}_{\gs, \chi_P}, \ccH^{-\infty}_{\gs, \chi_Q})$
such that
$$
A( Q ,  P, \gs, \nu) j(P, \gs, \nu)  = j(Q, \gs, \nu) B(Q, P, \gs, \nu), \qquad (\nu \in \faPdc).
$$  
Each space $\ccH^{-\infty}_{\gs,\chi_Q},$ for $Q\in \cP,$ carries a natural structure of Hilbert space.  
In terms of these structures, the Maass-Selberg relations can be formulated as
\begin{equation} 
\label{e: MS for B intro}
B(Q, P, \gs, -\bar\nu)^* B(Q, P, \gs, \nu) = \eta(Q, P, \gs, \nu),
\end{equation}
where $\eta(Q,P, \gs , \nu)$ is the scalar meromorphic function on $\faPdc = \faQdc$ determined by
$$
A(Q , P, \gs, -\bar \nu)^* A(Q , P, \gs, \nu) = \eta( Q, P,\gs, \nu).
$$
The interaction of these structures with the Weyl group $W(\fa)$ is discussed.

In Section \ref{s: Bmatrix, reduction} we study the operators  $B(Q, P, \gs, \nu)$ in detail. The Maass--Selberg relations (\ref{e: MS for B intro}) for $B$ are formulated in Theorem \ref{t: MS for B}. By using 
the well-known product decomposition of the  standard  intertwining operators
in terms of those with adjacent $P$ and $Q$ we reduce the proof
of the relations (\ref{e: MS for B intro}) to the setting in which $P,Q$ are adjacent. We discuss the well-known technique 
of chosing a subgroup $G^{(\ga)}$ of $G$ in which $P$ and $Q$  determine opposite maximal parabolic subgroups
$P({\ga})$ and $Q({\ga}).$ It is then shown that $B(Q,P,\gs, \nu)$ is essentially equal to the $B$-matrix for
$G^{(\ga)},  Q^{(\ga)}, P^{(\ga)}, \gs, \nu,$ see Lemma \ref{l: comparison B and B ga}. 
The proof of that lemma  requires comparison of distributions on $G$ with distributions on $G^{(\ga)}.$ This makes
it long and technical, see  Sections \ref{s: parabolics} - \ref{s: comparison B}.

In the end, the proof of the Maass-Selberg relations for $B$ is reduced to 
the basic setting  in which $G$ has compact center and  $P$ and $Q$ 
are opposite maximal parabolic subgroups, see Lemma  \ref{l: basic case MS}.

In Section \ref{s: c functions MS} we introduce the $C$-functions as coefficients in certain constant terms
of the Whittaker integrals, see Theorem \ref{t: asymp char C functions}. These $C$-functions can be expressed
in terms of the  standard  intertwining operators and the $B$-matrices, see Lemma \ref{l: C func and A B}.
We show that the Maass--Selberg relations for $B$ imply similar Maass-Selberg relations 
for $C$-functions $C_{Q|P}(s,\nu),$ 
with $Q,P$ associated and $s \in W(\faQ |\faP).$ 
They eventually  take the form $$C_{Q\mid P}(s,-\bar \nu)^*C_{Q\mid P}(s, \nu)= \eta_*(P,\bar P,-\nu)$$
see Theorem \ref{t: MS for C}.

Conversely, it is not a priori clear that the Maass-Selberg
relations for the $C$-functions imply those for the $B$-matrices. However, in the basic setting they
do, as is explained in Section \ref{s: MS basic setting}. 
The completion of the proof of the Maass-Selberg relations thus depends on their validity
for the $C$-functions in the mentioned basic setting. 
The latter case is addressed in the next two sections, \ref{s: radial part Casimir}  and \ref{s: result HC},
which are based on Harish-Chandra's work in \cite{HCwhit}. 

Section \ref{s: radial part Casimir} is preparatory, determining a useful formula for the radial part of the Casimir operator, which leads to a formula given without proof in  \cite[p.~208]{HCwhit}. 
That formula allows the application of Gauss' divergence
theorem for a simplex, which in turn leads to asymptotic information in Section 
\ref{s: result HC}, see Thm.  \ref{t: integral asympt HC}.  At the end of Section \ref{s: result HC}, the obtained  asymptotic 
information turns out to imply the Maass-Selberg relations for the $C$-functions in the basic setting.

In Section \ref{s: normalized Whittaker} we introduce the normalized Whittaker integrals $\Wh^{\circ}(P,\psi)$
and the associated normalized $C$-functions $C^{\circ}_{Q|P},$ for $P,Q \in \cPst,$ following the definitions of Harish-Chandra, \cite{HCwhit}. The Maass-Selberg relations imply the following relations for
the normalized $C$-functions, for $s \in W(\faQ|\faP),$
$$
C^\circ_{Q|P}(s, -\bar\nu)^*C^\circ_{Q|P}(s, \nu) = I\qquad (\nu \in \faPdc).
$$ 
These, combined with the uniformly tempered estimates obtained in \cite{Bunitemp}, allow us to show that the normalized Whittaker integrals are (finite sums of) 
functions of type $\twoholpr.$  In view of results in our paper \cite{Bcterm} (to apppear in the near future) this allows us, in Section \ref{s: nFou and Wave}, 
to define for each standard parabolic $P \in \cPst$ a normalized Fourier transform $\nFou_P$
which is continuous linear from the Harish-Chandra type Schwartz space $\cC(\tau: G/N_0: \chi)$ 
to the Euclidean Schwartz space $\cS(i\faPd, \cA_{2,P}).$ 
The conjugate Wave packet transform $\Wave_P$ is continuous linear between these Schwartz spaces in the converse direction.

In Section \ref{s: functional eqn} we establish the functional equations for 
the normalized Whittaker integrals as given by Harish-Chandra \cite[\S 17.1]{HCwhit}
In turn these imply transformation
formulas for the normalized $C$-functions, the normalized Fourier transform,
and the Wave packet transform. 

\section{Definition of the map $j(Q, \gs, \nu)$}
\label{s: defi j}
A parabolic subgroup $Q \in \cP$ is said to be opposite  standard  if $\bar Q$ is  standard.
The set of $Q \in \cP$ with $\bar Q \in \cPst$ is denoted by $\oppcPst.$
For this paper it will be necessary to describe the action of  standard  interwining operators 
on Whittaker vectors of parabolically induced representations. 
To make this possible, we need to extend the definition of the map $j(Q, \gs, \nu)$ to the setting 
of all parabolic subgroups $Q$ from $\cP,$ beyond those from $\oppcPst.$  
To prepare for this we start with the description of the open orbits $N_0 v Q$ 
in $G,$ for $v \in N_K(\fa).$ We assume that $Q \in \cP.$ 

\begin{lemma}
\label{l: v and open coset}
$G$ is a finite union of double cosets of the form $N_0 v Q,$ 
for $v \in N_K(\fa).$ The coset $N_0 v Q$ is open in $G$ if and only if $vN_Qv^{-1} \subset \bar N_0,$ which in turn is equivalent to the condition that $v Qv^{-1}\in \oppcPst$.
\end{lemma}

\proof
There exists an $s \in N_K(\fa)$ such that $P: = s Q s^{-1}$ is opposite  standard. By the Bruhat decomposition, $G$ is the disjoint union 
of the sets $P_0 v \bar P_0$ for $v \in W(\fa).$ Since $ P_0 v \bar P_0 = N_0 v \bar P_0 \subset N_0 v s  Q s^{-1},$ 
it follows that $G$ is a finite union of sets of the form $N_0 v Q s^{-1}$, with $v \in N_K(\fa).$  Hence, $G = G s$ is 
a finite union of orbits $N_0 v Q,$ for $v\in N_K(\fa).$

Put $\st\fn_0 := \fn_0 \cap \fm_Q.$ Then $N_0 v  Q$ is open in $G$ if and only if
the map $N_0 \times  Q \to G, $ $(n, q)\mapsto n v q$ is submersive at $(e,e).$  This is equivalent to
to $\fn_0 + \Ad(v) (\fq)  = \fg,$ which in turn is equivalent to
$\Ad(v)^{-1}\fn_0 + \fq = \fg.$ Since $\Ad(v)$ maps $\fa$-root spaces to $\fa$-root spaces
the latter assertion is equivalent to
$$ 
\Ad(v)^{-1} \fn_0 + \st \fn_0 + \st \bar\fn_0 +  \fn_Q =  \st \fn_0 + \st \bar\fn_0 +\bar \fn_Q + \fn_Q.
$$
This in turn is equivalent to 
$
\Ad(v)^{-1} \fn_0 \supset \bar \fn_Q,
$ 
hence to $v \bar N_Q  v^{-1} \subset N_0$ and to $v Q v^{-1} \supset \bar P_0.$ 
\qed

We denote by $W_Q(\fa)$ the centralizer of $\faQ$ in $W(\fa).$ 
The following lemma is well known through its formulation in terms of
root systems. 

\begin{lemma}
\label{l: W on parabs}
Let $Q \in \cP,$ $s, t \in W(\fa).$ 
\begin{enumerate}
\itema
If $s \in W(\fa)$ is such that $sQs^{-1} = Q$ then $s\in W_Q(\fa).$ 
\itemb
The group $Q$ is $W(\fa)$-conjugate to a unique $P\in \oppcPst.$ 
\end{enumerate}
\end{lemma}

\proof
We  start by proving (b). The existence part of (b) is well-known (and also follows from the previous  lemma).
So, there exists an $s \in W(\fa)$ and a $P \in \oppcPst$ such that $sQ s^{-1} = P.$
It follows that $s(\faQp) = \faPp.$ If $Q$ is opposite  standard , then $\faQp \subset - \cl(\fa^+).$ 
Fix $X \in \faQp;$ then both $X$ and $ sX$ belong to $  - \cl( \fa^+).$ Since the latter set is a fundamental domain for the  action of $W(\fa)$ on $\fa$, we conclude that   $ sX = X.$ Now
$s$ can be written as a product of simple reflections in roots vanishing on $- X$ hence on 
$\faPp.$ 
Therefore, $s \in W_P(\fa).$ This in turn implies that $\faQp = \faPp$ hence $P = Q$ and uniqueness follows.

Assume now that $Q$ is general, and let $t \in W(\fa)$ 
be such that $P':= t Q t^{-1}$ is  standard. Then $P'$ is $W(\fa)$-conjugate 
to $P$ and from the argument above 
it follows that $P = P'.$  This establishes uniqueness of $P.$ 

For (a) we fix $t\in W(\fa)$ and $P \in \oppcPst$ 
such that $t Q t^{-1} = P.$ Then $sQs^{-1} = Q$ 
implies that conjugation by $t s t^{-1}$ fixes $P.$ By (b) this implies that $t^{-1} s t $ 
belongs to $W_P(\fa).$ This in turn implies that $s \in t W_P(\fa)   t^{-1} = W_Q(\fa).$ 
\qed 

Let $Q \in \cP.$ Then by the lemma above there exists a unique $P \in \cPst$ that is Weyl conjugate to $Q.$ We fix $v\in N_K(\fa)$ such that $v Qv^{-1} = P.$ Then $vN_Qv^{-1} = N_P\subset N_0$ and by Lemma \ref{l: v and open coset} it follows that $N_0 v\bar Q$ is open 
in $G.$ The image of $v$ in $W_P(\fa)\bs W(\fa)$ is independent of the possible choices 
of $v.$  Likewise, the image of $v$ in $W(\fa)/W_Q(\fa)$ is also independent of such 
choices.

\begin{cor}
\label{c: open N v Q}
Let $Q \in \cP(A).$ Then precisely one of the $N_0$-orbits on $G/Q$ is open. 
This orbit equals $N_0 v Q$ for any $v\in N_K(\fa)$ such that $v Q v^{-1}$ is opposite  standard.
\end{cor} 

Corollary \ref{c: open N v Q} allows us to make the following choice once and for all.

\begin{defi}
\label{d: choice of vQ}
For the remainder of this paper we fix a map $Q \mapsto v_Q, $ $\cP(A) \to N_K(\fa)$ such that
\begin{enumerate}
\itema  for every $Q\in \cP(A)$ 
the double coset $N_0 v_\iQ  Q$ is open in $G;$
\itemb  if $Q \in \oppcPst $ then $v_Q = e.$ 
\end{enumerate}
\end{defi}

We define 
\begin{equation}
\label{e: type chi new}
C^{-\infty}(G/ Q:\gs : \nu)_\chi 
\end{equation}
to be the subspace of $C^{-\infty}(G/Q:\gs :\nu)$ consisting of $\gf$ such that
$\gf(n x) = \chi(n) \gf(x),$ for $ x\in G, n \in N_0$.
The elements of (\ref{e: type chi new}) are called the generalized functions of type $\chi.$ 
 
For $w \in N_K(\fa)$ we define the representation $w\gs$ of $wM_Qw^{-1}$ 
in $H_\gs$ by $w\gs := \gs \after w^{-1}.$ Furthermore, we write $H_{w\gs}$ 
for $H_\gs$ equipped with the representation $w\gs.$ Let 
$$ 
R_w: C^{-\infty}(G/ Q: \gs :\nu) \to C^{-\infty}(G/ w Q w^{-1} : w\gs : w\nu) 
$$ 
be the unique continuous linear $G$-intertwining operator which is given 
by the right regular action by $w$ on the subspace of smooth functions.
 
It maps functions of type $\chi$ for the left regular action by $N_0$ 
bijectively onto functions of the same type in the image space. 

Suppose  $Q \in \cP(A)$ and let $P \in \oppcPst$ be the unique opposite  standard 
parabolic subgroup that is $W(\fa)$-conjugate to $Q.$ Let $v \in N_K(\fa)$ be such that $vQv^{-1} = P;$ this condition
is equivalent to $v \in N_{K_P} (\fa) v_Q$ and to $v \in v_Q N_{K_Q}(\fa) .$

Clearly, $vM_Q v^{-1} = M_P.$ Since $M_P = K_P (M_P \cap A) (M_P \cap N_0) $
is an Iwasawa decomposition for $M_P$ it follows that 
\begin{equation}
\label{e: Iwasawa MQ} 
M_Q = K_Q (M_Q \cap A) (M_Q \cap v^{-1}N_0 v)
\end{equation}
is an Iwasawa decomposition for $M_Q.$  Its $N$-component is given by  
$
M_Q \cap v^{-1}N_0 v = v^{-1}(M_P \cap N_0) v.
$
The associated character 
$M_Q\cap  v^{-1}N_0 v \to \C$, given by 
\begin{equation}
\label{e: char vmin chi P}
 v^{-1}\chi_P:   \;\;n \mapsto \chi(v n v^{-1}),
\end{equation}
will in general depend on the particular choice of $v \in N_K(\fa).$ To avoid any ambiguity 
we agree to exclusively use the notation $\chi_Q$ for the character (\ref{e: char vmin chi P}) defined with $v = v_Q,$ see Definition \ref{d: choice of vQ}. 
\begin{lemma}
Put $v = v_Q.$ 
The character $\chi_Q$ of $M_Q\cap v(M_P \cap N_0)v^{-1}$defined by (\ref{e: char vmin chi P}) is regular with respect to the Iwasawa decomposition 
(\ref{e: Iwasawa MQ}).
\end{lemma}

\proof 
This is immediate from the regularity of $\chi|_{M_P\cap N_0}$
relative to the Iwasawa decomposition $M_P = K_P(M_P\cap A )(M_P\cap N_0),$ see 
 \cite[(8.12)]{Bunitemp}. 
\qed

In analogy with  \cite[(8.12)]{Bunitemp}, we define 
\begin{equation}
\label{e: Q eval set}
\ccH_{\gs, \chi_Q}^{-\infty} :=
\{
\eta \in H_\gs^{-\infty}\mid \forall m \in M_Q \cap v_Q^{-1}N_0v_Q :  \gs(m)\eta  = \chi_Q(m)\eta
\}.
\end{equation} 
In general, if $L$ is a closed subgroup of $G,$ $\pi$ a unitary representation of $L$ and $\xi$ a unitary 
character of a closed subgroup $N \subset L,$ we agree to write 
 $$
 \ccH_{\pi, \xi}^{-\infty}: = \{ \eta \in H_{\pi}^{-\infty} \mid\;\; \forall m \in N:\;\;\; \pi^{-\infty}(m)\eta  = \xi(m)\eta\}.
 $$%
Given $v \in G,$ we denote by $v\pi$ the unitary representation of $vLv^{-1}$ in $H_\pi$ given by $v\pi(y) = \pi(v^{-1}y v).$
 Furthermore, $v \xi$ denotes the character of $vNv^{-1}$ given by $v \xi(z) = \xi(v^{-1} z v).$ 
 It is readily verified that
 \begin{equation}
 \label{e: w on H pi xi}
 \ccH^{-\infty}_{\pi, \xi} = \ccH^{-\infty}_{v\pi, v\xi}.
\end{equation}%
Indeed, the space on the left consists of all $\eta \in H_\pi^{-\infty}$ such that
$\pi^{-\infty}(m)\eta = \xi(m)\eta$. Substituting $m' = v m v^{-1}$ we see that
the condition on $m' \in vNv^{-1}$ is $(v\pi)^{-\infty} (m') = v[ \pi^{-\infty}](m') = v\xi(m')$ 
for all $m'\in vNv^{-1}.$ This in turn is equivalent to $\eta \in  H^{-\infty}_{v\pi, v\xi}.$
  
\begin{cor}
\label{c: space chi Q and P}
 Let $Q \in \cP$  and let $P \in \oppcPst$ be $W(\fa)$-conjugate to $Q.$ 
 Let $v = v_Q.$  Then for all $\gs \in \dMQds,$
 $$ 
\ccH_{\gs, \chi_Q}^{-\infty} = \ccH_{v \gs, \chi_P}^{-\infty}.
$$
\end{cor}

\vspace{-12pt}
\proof
From (\ref{e: char vmin chi P}) we see that $\chi_Q = v^{-1} \chi_P.$ Now apply (\ref{e: w on H pi xi}).
\qed
 
\begin{lemma}
\label{l: chi invariance implies continuous}
Let $ Q \in \cP(A)$ and put $v = v_Q$ and $P = v Q v^{-1}.$ Then $P \in \oppcPst.$ 
For every $\gf \in C^{-\infty}(G/Q: \gs: \nu)_\chi$ the following assertions are valid:
\begin{enumerate}
\itema
the restriction of $\gf$ to the open subset $N_0 v Q$ is a continuous function $N_0 v  Q \to H_{\gs}^{-\infty};$
\itemb
$\gf(v) \in \ccH^{-\infty}_{\gs, \chi_Q}.$
\end{enumerate}
\end{lemma}

\proof
This can be proven in the same fashion as  \cite[Thm.\ 8.6]{Bunitemp} where $Q$ is assumed to be opposite  standard  and where $v = e.$ Alternatively, one may apply the mentioned
result as follows. The generalized function $R_v \gf \in C^{-\infty}(G/P: v\gs: v \nu)_\chi$, 
when restricted to $N_0 P,$ yields a continuous function 
$N_0  P \to H_{v\gs}^{-\infty}.$ Consequently, $\gf$ is continuous on
$N_0  P v$ with values in $H_{v\gs}^{-\infty} = H_\gs^{-\infty}.$
Now $N_0 P v =N_0 v Q$ so (a) follows. 

For (b) we note that $\gf(v)= R_\gf(e) \in \ccH_{v\gs, \chi_P}^{-\infty},$ by the mentioned result.
We now apply Corollary \ref{c: space chi Q and P}.
\qed
\begin{cor}
\label{c: diagram ev}
Let $Q\in \cP(A),$ $v = v_Q,$ and $P: = v Q v^{-1}.$ 
Then for every $\nu \in \faQdc$ the following  is a commutative diagram of linear maps.
All appearing maps are linear isomorphisms between finite dimensional spaces.
$$ 
\begin{array}{ccc}
C^{-\infty}(G/Q: \gs :\nu)_\chi & \;\;\;\;{\buildrel R_v \over \longrightarrow} \;\;\;\;& C^{-\infty}(G/P : v\gs : v\nu)_\chi \\
 \downarrow {\scriptstyle \ev_v } &  &\downarrow {\scriptstyle \ev_e } \\
 \ccH_{\gs, \chi_Q}^{-\infty}\;\;\;\; &{\buildrel = \over \longrightarrow} &\!\!\!\!\!\!\!\ccH_{ v\gs,\chi_P}^{-\infty}
 \end{array}
 $$
\end{cor}

\proof
Fix $\nu\in \faQdc.$  It follows from the arguments of the 
above lemma that for $\gf\in C^{-\infty}(G/Q: \gs :\nu)_\chi $ we have 
$\ev_e R_v  \gf = \gf(v) = \ev_v  \gf.$ Hence the diagram commutes.
It follows from  \cite[Cor.~14.5]{Bunitemp} that the vertical map on the right is a linear isomorphism 
of finite dimensional linear spaces. Clearly the horizontal maps are linear isomorphisms.
Therefore, the vertical map on the left is a linear isomorphism.  Since the linear spaces
on the right are finite dimensional, all appearing spaces are. 
\qed

We agree to write $\ev_Q$ for the evaluation map 
$$ 
\ev_{v_Q}: 
C^{-\infty}(G/Q:\gs :\nu)_\chi \to \ccH^{-\infty}_{\gs, \chi_Q}
$$%
appearing in the left column of the diagram of Corollary \ref{c: diagram ev}.
The above definition depends on our choice of $v_Q \in N_K(\fa).$ Any alternative choice
$v\in N_K(\fa)$ must satisfy $v Q v^{-1} = P$ or, equivalently, $v = v_Q u$ with $u \in N_{K_Q}(\fa).$ 
The following result expresses the dependence of our definition of $\ev_Q$ on the choice of $v_Q.$

\begin{lemma}
\label{l: change of choice vQ}
Let $Q\in \cP$ and let $v = v_Q u$ with $u\in N_{K_{Q}}(\fa).$ 
Then for every $\nu \in \faQdc$ the following diagram commutes
$$ 
\begin{array}{ccc}
&C^{-\infty}(G/Q:\gs:\nu)_\chi & \\
&&\\
\;\;\;\;\;\;\;\;\;\;\;\;{\scriptstyle\ev_Q}\!\!\!\!\!\!\!\!\!\!\!\!\!\!&\swarrow \qquad \searrow&\!\!\!\!\!\!\!\!\!\!\!\!\!\!{\scriptstyle \ev_v}\;\;\;\;\;\;\;\;\;\;\;\;\\
&&\\
\;\;\;\;\;\;\ccH_{\gs, \chi_Q}^{-\infty}\!\!\!\!\!\!\!\!\!\!\!\!\!\!&{\buildrel\gs(u)^{-1}\over \longrightarrow}& \!\!\!\!\!\!\!\!\!\!\!\!\!\!\ccH_{u\gs, \chi_Q}^{-\infty}.\;\;\;\;\;\;
\end{array}
$$
The  horizontal map at the bottom is given by $\eta \mapsto \gs(u)^{-1} \eta.$
All maps in the diagram are linear isomorphisms.
\end{lemma}

\proof
Let $\gf \in C^{-\infty}(G/Q:\gs:\nu)_\chi.$ Then $\ev_v(\gf) = \gf(v) = \gf(v_Q u) = \gs(u)^{-1} \gf(v_Q) = \gs(u)^{-1}\ev_Q(\gf).$ Now $\ev_Q$ is a linear isomorphism and the map 
$\tau: \eta \mapsto \gs(u)^{-1}\eta$ is a linear automorphism of $H_\gs^{-\infty}.$ 
Hence, it suffices to show that $\tau$ maps 
 $\ccH_{\gs, \chi_Q}^{-\infty}$ onto $\ccH_{u\gs, \chi_Q}^{-\infty}.$ Let $\eta \in 
 \ccH_{\gs, \chi_Q}^{-\infty},$ and suppose that $n \in M_Q\cap v_Q^{-1} N_0 v_Q.$ 
 Then 
 $$
 u \gs(n)[\tau(\eta)] = 
 u\gs(n) [\gs(u)^{-1} \eta ]= \gs(u^{-1} n)\eta = \gs(u^{-1}) \chi_Q(n)\eta = 
 \chi_Q(n)\tau(\eta),
 $$
so that $\tau(\eta) \in \ccH_{u\gs, \chi_Q}^{-\infty}.$ In a similar fashion
it is shown that $\tau^{-1} (\ccH_{u \gs, \chi_Q}^{-\infty}) \subset \ccH_{\gs, \chi_Q}^{-\infty}.$
 \qed

\begin{defi}
Let $Q \in \cP.$ 
For each $\nu \in \faQdc$ the map 
$$
j(Q, \gs ,\nu): \ccH_{\gs, \chi_Q}^{-\infty} \to C^{-\infty}(G/ Q: \gs :\nu)_\chi
$$ 
is defined to be the inverse of 
$\ev_Q := \ev_{v_Q} :  C^{-\infty}(G/ Q: \gs :\nu)_\chi \to \ccH_{\gs,\chi_Q}^{-\infty} .
$
\end{defi}

\begin{cor}
\label{c: diagram j}
Let $Q\in \cP,$ and put $v = v_Q$ and $P:= vQv^{-1}.$ Then $P \in \oppcPst.$ 
For every $\nu \in \faQdc,$ the following  is a commutative diagram of linear maps.
$$
\begin{array}{ccc}
C^{-\infty}(G/Q: \gs :\nu)_\chi & \;\;\;\;{\buildrel R_v \over \longrightarrow} \;\;\;\;& C^{-\infty}(G/ P : v\gs : v\nu)_\chi \\
 \uparrow {\scriptstyle j(Q, \gs,\nu)} &  &\uparrow {\scriptstyle \scriptstyle j(P, v\gs, v\nu) } \\
 \ccH^{-\infty}_{\gs, \chi_Q}\;\;\;\; &{\buildrel = \over \longrightarrow} & 
 \!\!\!\!\!\!\!\!    \ccH^{-\infty}_{ v\gs, \chi_P}
 \end{array}
 $$
 All maps are  linear isomorphisms of finite dimensional spaces.
\end{cor}
\proof
This is immediate from Cor. \ref{c: diagram ev}.
\qed

We chose a non-degenerate $\Ad(G)$-invariant symmetric bilinear form 
\begin{equation}
\label{e: defi B} 
B: \fg\times \fg \to \R
\end{equation}
as in \cite[(2.1)]{Bunitemp} and define an $\Ad(K)$ positive definite inner product
on $\fg$ by $\inp{X}{Y}: = -B(X,\theta Y).$ 
For $Q \in \cP$ and $R > 0$ we put
\begin{equation}
\label{e: defi fa Q R}
\fad(Q, R) = \{\nu \in \faQdc\mid \inp{\Re \nu}{\ga} > R \; (\forall \ga \in \gS(\fn_Q, \fa_Q))\}.
\end{equation}
Note that for $w \in N_K(\fa)$ we have $w \fad(Q,R) = \fad( w Q w^{-1}, R).$
\begin{cor}
\label{c: holomorphy j Q}
Let $Q\in \cP,$ $\gs \in \dMQds$ and $\eta \in H^{-\infty}_{\gs, \chi_Q}.$ For every $R \in \R$ there exists a 
positive integer $s$ such that the assignment
$\nu \mapsto j( Q ,\gs, \nu)$ is holomorphic as a map 
$\fad(\bar Q, R) \to C^{-s}(K/K_Q:\gs_Q).$  
\end{cor}

\proof 
If $Q$ is opposite  standard , this follows from  \cite[Lemma 14.3]{Bunitemp}. For general $Q \in \cP$ we observe that
it follows from  Cor.\ \ref{c: diagram j} that for a fixed $\eta \in \ccH^{-\infty}_{\gs, \chi_Q}$ 
we have $j(Q, \gs, \nu)\eta = R_{v^{-1}} \after j(P, v\gs, v\nu)$ for all $\nu \in \faQdc.$ 
Now $R_{v^{-1}}$ restricts to a continuous linear map $C^{-\infty}(K:v \gs) \to C^{-\infty}(K:\gs),$ 
independent of $\nu.$ Therefore, the required result follows from the established case.
\qed

 \begin{rem}
 \label{r: order of j}
  If $\Omega \subset \faQdc$ is a bounded open subset, then it follows from Cor. \ref{c: holomorphy j Q}
that there exists a positive integer $s$ such that for all $\eta$ the assignment $\nu \mapsto j( Q ,\gs, \nu)\eta$ is holomorphic as a map 
$\Omega \to C^{-s}(K/K_Q:\gs_Q).$ The smallest $s$ with this property will be called the order of $j(Q,\gs,\dotvar)$
over  $\Omega.$ 
\end{rem}

\begin{rem}
The definition of the space $C^{-s}(K/K_Q:\gs_Q)$ appearing in the preceding statements
 is explained in \cite[\S 7]{Bunitemp}. 
 \end{rem}

If $P\in \oppcPst$ we fix an arbitrary positive invariant density $d\bar m_P$ on $M_P/M_P\cap N_0.$
If $Q\in \cP$ is conjugate to $P,$ then $v_Q Q v_Q^{-1} = P$ and conjugation by $v_Q$ induces a
diffeomorphism from $M_Q/M_Q \cap v_Q^{-1} N_0 v_Q$ onto $M_P/M_P\cap N_0.$ The pull-back 
of $d\bar m_P$ under $C_{v_Q}$ is a positive invariant density on $M_Q/M_Q \cap v_Q^{-1} N_0 v_Q$
which we denote by $d\bar m_Q.$ 

\begin{lemma}
Let $Q, Q'\in \cP$ be $N_K(\fa)$-conjugate to the same 
opposite  standard  parabolic subgroup $P.$ 
Put $w = {v_Q'}^{-1}v_Q.$ Then $Q' = wQ w^{-1},$  
$w(M_Q \cap v_Q^{-1} N_0 v_Q)w^{-1} = (M_{Q'} \cap v_{Q'}^{-1} N_0 v_{Q'})$ and 
$$
\cC_w^*(d \bar m_{Q'}) = d \bar m_Q.
$$ 
\end{lemma}

\proof
The first assertions are evident.
By definition of $w$ we have $\cC_w^* = \cC_{v_Q}^* \cC_{v_{Q'}}^{*-1}.$
Hence 
$$
\cC_w^*(d\bar m_{Q'}) = \cC_{v_Q}^*(d \bar m_P) = d \bar m_Q.$$
\qed

Let $Q\in \cP,$ $\gs \in \dMQds.$ According to  \cite[Lemma 9.2]{Bunitemp}
with $M_Q$ in place of $G,$ we have the matrix coefficient map 
\begin{equation}
\label{e: matrix coefficient map mu for Q}
\mu_\gs= \mu_{Q,\gs}: H_{\gs} \otimes \overline{\ccH_{\gs, \chi_\iQ}^{-\infty}} \to
L^2(M_Q/M_Q \cap v_Q^{-1}  N_0 v_Q:\chi_Q)
\end{equation} given by
$
\mu_{\gs}(z \otimes \eta) (m) = \inp{\gs(m)^{-1} z}{\eta},
$
for $z \in H_\gs^\infty, \eta \in \ccH_{\gs, \chi_Q}^{-\infty}$ and $m \in M_Q.$ 
Here the $L^2$-norm is defined with respect to the invariant measure $d\bar m_Q$ 
on $M_Q/M_Q\cap v_Q^{-1} N_0 v_Q.$ 

From  \cite[Cor. 9.5]{Bunitemp} it follows that the finite dimensional space 
$\ccH_{\gs, \chi_\iQ}^{-\infty}$ carries a unique inner product such that the
map (\ref{e: matrix coefficient map mu for Q}) is an isometric linear map
onto a closed subspace, which we denote by 
$$
L^2(M_Q/M_Q \cap v_Q^{-1}  N_0 v_Q:\chi_Q)_\gs.
$$ 

\begin{defi}
\naam{d: unique ip on H spaces}
Let $Q\in \cP$ and $\gs \in \dMQds.$ From now on we assume that $H_{\gs, \chi_Q}^{-\infty}$ 
is equipped with the unique Hermitean inner product that makes $\mu_{Q, \gs}$ isometric.
\end{defi}

Let $Q'$ be conjugate to $Q$ and put $w= v_{Q'}^{-1}v_Q.$ Then it is readily verified that
conjugation by $\cC_w^{-1}$ defomes a diffeomorphism from $M_Q$ onto 
from $M_{Q'},$ which maps $M_Q\cap v_Q^{-1} N_0 v_Q$ onto 
the similar intersection with everywhere $Q$ replaced by $Q'.$ 
In turn this implies that pull-back under conjugation $C_{w^{-1}}$
induces an isometric isomorphism
$$
A_w : L^2(M_Q/M_Q\cap v_Q^{-1}N_0v_Q: \chi_Q)_\gs
\to 
 L^2(M_{Q'}/M_{Q'}\cap v_{Q'}^{-1}N_0v_{Q'}: \chi_{Q'})_{w\gs}.
 $$ 

\begin{lemma} 
\label{l: identities unitary for cHgs} 
Let $Q,Q' \in \cP$ be such that $wQw^{-1} = Q'$ for $w = v_{Q'}^{-1} v_Q.$ 
Let $\gs \in \dMQds.$ 

The linear spaces  $\ccH_{\gs, \chi_Q}^{-\infty}$ and $\ccH_{w\gs, \chi_{Q'}}^{-\infty}$ 
are equal as subspaces of $\ccH_{\gs}^{-\infty}.$ The Hermitean inner products, as specified in Definition
\ref{d: unique ip on H spaces} are the same. 
\end{lemma}

\proof
Put $v_Q = v$ and $v' = v_{Q'}. $ Then $w = {v'}^{-1}v.$ 
It follows from Corollary 2.6 \ref{c: space chi Q and P} that 
$$ 
\ccH_{\gs, \chi_Q}^{-\infty} = \ccH_{v \gs, \chi_P}^{-\infty}  
= \ccH_{{v'}^{-1} v \gs , \chi_{Q'}}^{-\infty}=  \ccH_{w \gs, \chi_Q'}^{-\infty}.
$$  
In view of our choices of measure, the map $A_w$ given above is an isometry. For $z \in H_{\gs}$ and 
$\eta \in \ccH_{\gs, \chi_Q}^{-\infty}$ 
we have that
$$
A_w \after [\mu_{Q, \gs}(z \otimes \eta)](m) = \mu_{Q, \gs}(z \otimes \eta)(w^{-1} m w) = \mu_{Q', w\gs} (w\gs \otimes \eta)(m).
$$ 
The maps $\mu_{Q, \gs}$ and $\mu_{Q', w\gs}$ are unitary by definition, and we see that the identity 
induces a unitary map 
$$
H_\gs \otimes \overline{\ccH_{\gs, \chi_Q}^{-\infty}} \to H_{w\gs} \otimes \overline{\ccH_{w\gs, \chi_{Q'}}^{-\infty}}.
$$
Since the identity map $H_\gs \to H_{w\gs}$ is unitary, we conclude that
the identity map $\ccH_{\gs, \chi_Q}^{-\infty}\to \ccH_{w\gs, \chi_{Q'}}^{-\infty}$ is unitary as well.
\qed

\begin{lemma} 
\label{l: gs u inv isom}
Let $Q \in \cP$ and $u \in N_{K_Q}(\fa).$ Then the map $\eta \mapsto \gs(u)^{-1}\eta$
is an isometry from $\ccH_{\gs, \chi_Q}^{-\infty}$ onto $\ccH_{u \gs, \chi_Q}^{-\infty}.$ 
\end{lemma}

\proof
From Lemma \ref{l: change of choice vQ} it follows that  the linear map 
$\tau: H_\gs^{-\infty}\to H_\gs^{-\infty},$ $\xi \mapsto 
\gs(u)^{-1} \xi$ maps $H^{\infty}_{\gs, \chi_Q}$ onto $H^{\infty}_{u \gs, \chi_Q}.$ 
We will finish the proof by showing that $\tau$ is isometric on 
$\ccH^{-\infty}_{\gs, \chi_Q}.$ 

For $v \in \ccH_\gs^{\infty}$ and $\eta \in \ccH_{\gs, \chi_Q}^{-\infty}$ the matrix coefficient
attached to $v\otimes \eta$ is the function in 
$C^\infty(M_Q/M_Q  \cap v_Q^{-1} N_0 v_Q: \chi_Q)$ defined 
by 
$$
\mu_{v\otimes \eta}(m)  = \inp{\gs(m)^{-1}v}{\eta}, \qquad (m \in M_Q). 
$$ 
The sesquilinear map $(v, \eta) \mapsto m_{v \otimes \eta}$ induces a linear isometry
from the pre-Hilbert space $H_\gs^\infty \otimes \ccH_{\gs, \chi_Q}^{-\infty}$ to $L^2(M_Q/M_Q\cap v_Q^{-1} N_0 v_Q: \chi_Q). $ This implies that for all $v \in H_\gs^\infty$ 
and $\eta \in \ccH_{\gs, \chi_Q}^{-\infty},$
\begin{equation}
\label{e: L2 norm of matrix coeff}
\|v\|^2_\gs \|\eta\|_{\gs, \chi_Q}^2= \int_{M_Q/M_Q\cap v_Q^{-1}N_0v_Q}|\mu _{v \otimes \eta}(m)|^2 \; d\bar m_Q.
\end{equation}

The representation $u\gs $ of $M_Q$ defined by $u\gs(m) = \gs(u^{-1}mu)$ 
is irreducible unitary and belongs to the discrete series of $M_Q$ again. 
We write $\ccH_{u\gs}$ for the Hilbert space $\ccH_\gs$ equipped with the representation
$u\gs$ and will discuss the induced inner product on $\ccH_{u\gs, \chi_Q}^{-\infty}.$ 
The identity map is unitary from $H_\gs$ to $H_{u\gs}.$ 
If $v \in H_\gs^\infty,$ then $v \in H_{u\gs}^\infty.$  If $\eta' \in \ccH_{u\gs, \chi_Q}^{-\infty}$ 
then the associated matrix coefficient $\mu_{v \otimes \eta'}'$ is given by 
$$ 
\mu_{v\otimes \eta'}' (m) = \inp{[u\gs](m)^{-1}v}{\eta'} = \inp{\gs(u^{-1})\gs(m^{-1}) \gs(u)v }{\eta'} =
\mu_{v \otimes \gs(u)\eta'}(u^{-1} m).
$$ 
Substituting $\eta' = \tau (\eta)$ we find that
$$ 
\mu_{v\otimes\tau(\eta)}'(m)   = \mu_{v\otimes \eta}(u^{-1} m).
$$
Using the analogue of (\ref{e: L2 norm of matrix coeff}) for $u\gs$ and $\mu',$
combined with the left invariance of the measure $d\bar m_Q,$ we infer that
$$
\|\tau(\eta)\|_{u\gs, \chi_Q}^2 = \|\eta\|_{\gs, \chi_Q}^2, \qquad (\eta \in\ccH^{-\infty}_{\gs, \chi_Q}).
$$ 
Hence, $\tau$ is isometric as stated.
\qed
In  \cite{Bunitemp} the Whittaker integral $\Wh(P,\dotvar)$ for  $P \in \cPst$ is expressed in terms of matrix coefficients
involving $j(\bar P,\dotvar).$ Guided by this definition we will now use matrix coefficients involving $j(\bar Q, \dotvar)$ to define the notion of Whittaker integral $\Wh(Q, \dotvar),$ 
for $Q \in \cP$  arbitrary.

As in  \cite[\S 9]{Bunitemp}, which in turn relies on Harish-Chandra 
 \cite[Lemmas 7.1,9.1]{HCha3} we aim at defining a linear isomorphism 
$$
T \mapsto \psi_T, \;\;C^\infty(\tau: K/K_Q: \gs_Q) \otimes \overline{\ccH^{-\infty}_{\gs, {\chi_\iQ}}}
\to L^2(\tau_Q : M_Q / M_Q \cap v_Q^{-1}N_0 v_Q: \chi_Q)_\gs.
$$ 
Here $\gs \in \dMQds$ and $\gs_Q$ denotes the restriction of $\gs$ to $K_Q.$ 
Furthermore, $(\tau, \Vtau)$ is a finite dimensional unitary representation of
$K$ and  $\tau_Q$ denotes the restriction of $\tau $ to $K_Q = K \cap M_Q.$

We define the space of spherical functions
\begin{equation}
\label{e: spherical L^2 space MQ}
L^2(\tau_Q: M_Q/M_Q \cap v_Q^{-1} N_0 v_Q: \chi_Q)
\end{equation}
to be the subspace of $K_Q$-fixed elements in $L^2(M_Q/M_Q \cap v_Q^{-1} N_0 v_Q: \chi_Q) \otimes \Vtau.$
Viewing (\ref{e: spherical L^2 space MQ}) naturally as a space of functions $M_Q\to \Vtau$ we shall express the spherical behavior of
its functions by $f(k m) = \tau(k) f(m) $, for $m \in M_Q$ and $k \in K_Q.$ The space is equipped with
the restriction of the tensor product Hilbert structure, and thus is a Hilbert space of its own right.
The space of functions in (\ref{e: spherical L^2 space MQ}) which belong to $L^2(M_Q/M_Q \cap v_Q^{-1} N_0 v_Q: \chi_Q)_\gs \otimes \Vtau$
is indicated by the subscript $\gs$ on the right. Since only finitely many representations of the discrete
series of $M_Q$ have a $K_Q$-type in common with $\tau_Q,$ it follows that
$$
L^2(\tau_Q: M_Q/M_Q\cap v_Q^{-1} N_0 v_Q: \chi_Q) =
\bigoplus_{\gs \in \dMQds}  L^2(\tau_Q: M_Q/M_Q\cap v_Q^{-1} N_0 v_Q: \chi_Q)_\gs
$$ 
is finite dimensonal. In particular the (orthogonal) sum over the $\gs$ is finite. 
From this it also follows that
$$
L^2(\tau_Q: M_Q/M_Q\cap v_Q^{-1} N_0 v_Q: \chi_Q) 
=
\cC(\tau_Q: M_Q/M_Q\cap v_Q^{-1} N_0 v_Q: \chi_Q),
$$ 
where the definition of the space on the right is obvious.
By finite dimensionality, the center $\fZ_{Q}$ of $U(\fm_Q)$ acts finitely on the space on the
right. For this  reason, that space is also denoted by $\cA_{2,\bar Q}.$

The subspace 
$\cC(M_Q/ M_Q \cap v_Q^{-1} N_0 v_Q:\chi_Q)_\gs \otimes \Vtau \cap \cA_{2, \bar Q}$ is denoted by 
$\cA_{2, \bar Q, \gs}.$ We have the finite orthogonal direct sum 
$$
\cA_{2,\bar Q} = \oplus_{\gs \in \dMQds} \cA_{2, \bar Q,\gs}.
$$

After these preparations 
we define, for $\gf \in C^\infty(\tau_Q: K/K_Q : \gs_Q)$ and $\eta \in  \ccH_{\gs, \chi_Q}^{-\infty}$ 
the function $\psi_{\gf\otimes \eta}: M_Q \to V_\tau$ by 
$$ 
\psi_{\gf \otimes \eta}(m) =\inp{ (\gs(m^{-1}) \otimes I)\gf(e)}{\eta}_{\gs,1},\qqquad (m \in M_Q),
$$ 
where the sesquilinear pairing $\inp{\dotvar}{\dotvar}_{\gs, 1}: (H_\gs \otimes \Vtau) \times \ccH^{-\infty}_{\gs, \chi_Q}
 \to V_\tau$ 
is given by $\inp{z\otimes v}{\eta}_{\gs, 1} = \inp{z}{\eta}_\gs v.$ 

\begin{lemma}
\label{l: identity on cHgs is isometry}
The map $(\gf, \eta)\mapsto \psi_{\gf \otimes \eta}$ induces an isometric linear isomorphism
$T \mapsto \psi_T$ from $C^\infty(\tau : K /K_Q : \gs_Q)  \otimes \overline{\ccH_{\gs, \chi_Q}^{-\infty}}$
onto $\cA_{2, \bar Q, \gs}.$ 
\end{lemma}
\proof
The proof, which relies on an application of Frobenius reciprocity, is identical to the proof 
for the case that $Q=P \in\cPst,$  in  \cite[Lemma 9.8]{Bunitemp}.
\qed

Finally, we are prepared to define the Whittaker integral associated with $Q \in \cP.$
\begin{defi}
\label{d: defi j Q}
Let $Q \in \cP.$ The Whittaker integral $\Wh(Q, \psi, \nu),$ for $\psi \in \cA_{2, Q}$ and for generic $\nu \in \faQdc$ 
is defined to be the function in $C^\infty(\tau: G/N_0:\chi)$ determined by the following requirements.
\begin{enumerate}
\itema $ \Wh(Q, \psi, \nu)$ depends linearly on $\psi \in \cA_{2,  Q}$;
\itemb for $\gs \in \dMQds$ and $T =\gf \otimes  \eta  \in C^\infty(\tau: K/K_Q: \gs_Q)\otimes
 \ccH_{\gs, \chi_{\bar Q}}^{-\infty},$ we  have
 $$ 
\Wh(Q, \psi_T, \nu)(x) = \inp{ \pi_{\bar Q, \gs, -\nu}(x)^{-1} \gf}{j(\bar Q, \gs, \bar \nu)\eta}, \qquad (x \in M_Q).
$$ 
\end{enumerate} 
\end{defi}

We retain our assumption that $Q \in \cP,$ $\gs \in \dMQds$ and put
$v = v_{\overline Q}.$ Then $P:= v Q v^{-1}$ is the unique  standard  parabolic subgroup 
which is conjugate to $Q.$ The Whittaker integral $\Wh(Q)$ can now be expressed in terms of $\Wh(P).$

\begin{lemma}
\label{l: existence cR Q}
There exists a unique isometric linear isomorphism  $\cR_Q: \cA_{2, Q} \to \cA_{2, P}$ such that for all $\gs \in \dMQds$ the following is valid
\begin{enumerate} 
\itema $\cR_Q$ maps $\cA_{2, Q, \gs}$ onto $\cA_{2, P, v \gs}.$
\itemb for all $T   \in C^\infty(\tau: K/K_Q: \gs_Q)\otimes
 \ccH_{\gs, \chi_{\bar Q}}^{-\infty}$ we have
\begin{equation}
\label{e: formula cR Q}
 \cR_Q \psi_T = \psi_{(\cR_v  \otimes I) T}.
 \end{equation}
Here $R_v: C^\infty(\tau: K/K_Q: \gs_Q) \to C^\infty(\tau: K/K_P: (v\gs)_P)$ is
the map induced by right translation by $v.$  \end{enumerate}
\end{lemma}

\proof
Since $\cA_{2, Q}$ decomposes
as the orthogonal finite direct sum of the subspaces $\cA_{2, Q, \gs}$ 
and since for each $\gs$ the map $\Psi_\gs: T \mapsto \psi_T$ is an isometry
from $V_\gs:= C^\infty(\tau: K/K_Q: \gs_Q)\otimes
 \ccH_{\gs, \chi_{\bar Q}}^{-\infty}$ onto $\cA_{2,Q,\gs},$ uniqueness is obvious.
  
 For each $\gs \in \dMQds$ the map $R_v: C^\infty(\tau: K/K_Q: \gs_Q) \to C^\infty(\tau: K/K_P: (v\gs)_P)$ is induced by right translation by $v,$ which is clearly an isometry.
Furthermore, $I$ denotes the identity map from $\ccH^{-\infty}_{\gs, \chi_{\bar Q}}$ to $\ccH^{-\infty}_{v\gs,  \chi_{\bar P}} ,$ which is an isometry by Lemma \ref{l: identity on cHgs is isometry}.  We thus see that 
${}_\gs R_v: T \mapsto  (R_v \otimes I) T$ 
 is an isometry from $V_\gs$ onto $V_{v\gs}.$ 
The condition on $\cR_Q$ may now be reformulated as 
$\cR_Q \after \Psi_\gs  = \Psi_{v\gs} \after {}_\gs R_v$  on $V_\gs.$
Thus, $\cR_Q  = \Psi_{v\gs} \after {}_\gs R_v \after \Psi_\gs ^{-1}$ 
on $\cA_{2, Q, \gs}.$ This establishes the existence.
\qed

\begin{prop}
Let $\cR_Q$ be as in Lemma \ref{l: existence cR Q}.
For all $\gs \in \dMQds$ and all $\psi \in \cA_{2, Q, \gs}$ we have
$$
\Wh(Q, \psi, \nu, x) = \Wh(P, \cR_Q \psi, v \nu, x),
$$
for all $\nu \in \faQdc$ and $x \in G.$ 
\end{prop}
\proof
Let $\gs \in \dMQds$ and $T = \gf \otimes \eta \in C^\infty(\tau:K/K_Q:\gs_{\bar Q}) \otimes \ccH_{\gs, \chi_Q}^{-\infty}.$
Then by Definition \ref{d: defi j Q} and Corollary \ref{c: diagram j} we have 
\begin{eqnarray*}
\Wh(Q, \psi_T, \nu, x) 
& = & 
\inp{\pi_{\bar Q, \gs, -\nu}(x)^{-1} \gf}{j(\bar Q, \gs, \bar \nu, \eta)}\\
&=& 
\inp{R_v \pi_{\bar Q, \gs, -\nu}(x)^{-1} \gf}{j(\bar P, v\gs, v\bar \nu,  \eta)}\\
&=&
\inp{\pi_{\bar P, v\gs, -v \nu}(x)^{-1} R_v \gf}{j(\bar P, v\gs, v\bar \nu,  \eta)}\\
&=&
\Wh(P,  \psi_{(R_v \otimes I)T}, v \nu, x).
\end{eqnarray*}
The proof is completed by using (\ref{e: formula cR Q}).
\qed

\begin{lemma} 
Let $Q\in \cP$ and let $v = v_{\bar Q}.$ Then $P = v Q v^{-1}$ is  standard.
Furthermore, for $\psi \in \cA_{2, Q}$ and $\nu \in \faQdc$ 
such that $\Re \inp{v\nu }{\ga}> 0$ for all $\ga \in \gS^+$ we have
$$
\Wh(Q,\psi, \nu)(x) = \int_{N_Q} [ v^{-1}\chi_{\bar Q}](n )^{-1} \psi_{\bar Q, \gs, -\nu}(x v n) \; dn
\qquad(x \in G).
$$ 
\end{lemma}

\begin{rem} For $Q$  standard  and $v = e,$ we retrieve Harish-Chandra's
formula.
\end{rem}
 
 \proof 
 First of all, by Definition \ref{d: choice of vQ}, $v\bar Qv^{-1} \in \oppcPst$ and it follows that
 $P$ is  standard. By linearity, it suffices to consider the case that
 $\psi = \psi_T,$ with $T = \gf \otimes \eta \in C^{\infty}(\tau: K/K_Q: \gs_Q) \otimes 
\ccH_{\gs, \chi_{\bar Q}}^{-\infty},$ where $\gs \in \dMQds.$ 
Write $\psi_{-\nu}$ for the function in $C^{\infty}(\tau: G/\bar Q:\gs :-\nu)$ given by 
$\psi_{-\nu}|_{K} =\psi.$ 
Then 
\begin{eqnarray*}
\Wh(Q, \psi_T, \nu)(x) &=& \inp{\pi_{\bar Q, \gs, -\nu}(x)^{-1}\gf}{j(\bar Q, \gs, -\nu)\eta}
\\
&=& \int_{K/K_Q}  \inp{\gf_{-\nu} (x k)}{j(\bar Q, \gs, -\nu,\eta)(k)} \; dk\\ 
&= & \int_{N_Q}  \inp{\gf_{-\nu} (x v n )}{j(\bar Q, \gs, -\nu,\eta)(v n)}\; dn\\
&=& \int_{N_Q} \chi(vn v^{-1})\;  \inp{\gf_{-\nu} (x v n )}{\eta}\; dn.
\end{eqnarray*}
The required representation now follows from the equality 
$$
\inp{\gf_{-\nu} (y)}{\eta} = \psi_{\bar Q, \gs, - \nu}(y),\qquad (y \in G).
$$ 
By left $\tau$-sphericality and right $\bar N_Q$-invariance,
it suffices to prove the latter equality for $y = m a \in M_Q A_Q.$ This in turn
follows from
$$ 
\inp{\gf_{-\nu} (ma)}{\eta} =  a^{\nu + \rho_Q} \inp{\gs(m)^{-1} \gf(e)}{\eta} =
a^{\nu + \rho_Q} \psi_T(m) = \psi_{\bar Q, \gs, - \nu}(ma).
$$ 
 \qed
 
\section{Interaction with the Weyl group}  
\label{s: interaction with W}

We assume that $Q\in \cP$ and that $\gs \in \dMQds.$ 

\begin{lemma}
\label{l: diagram composition Rw and j}
Let $w \in N_K(\fa),$ and put $Q' = wQ w^{-1}.$ 
There exists a unique linear map 
$$
\cR_{w, Q}: \ccH_{\gs, \chi_\iQ}^{-\infty}\to \ccH^{-\infty}_{w\gs, \chi_{\iQ'}}
$$ 
such that  for every $\nu \in \faQdc$ the 
following diagram commutes:
\begin{equation}
\label{e: diagram with Rw and ev}
\begin{array}{ccc}
C^{-\infty}(G/Q : \gs : \nu)_\chi  & \;\;\;\;{\buildrel R_{w}  \over \longrightarrow} \;\;\;\;& 
C^{-\infty}(G/Q': w\gs :w\nu)_\chi \\
\downarrow {\scriptstyle \scriptstyle \ev_Q}
 &  &
 \downarrow {\scriptstyle \ev_{Q'}}\\
 \!\!\!\!\!\!\!\!   \ccH_{\gs, \chi_Q }^{-\infty}&{\buildrel \cR_{w,Q}  \over \longrightarrow} & 
 \!\!\!\!\!\!\!\!  \ccH_{w\gs, \chi_{Q'}  }^{-\infty}\;\;\;\; 
\end{array}
\end{equation}
The map $\cR_{w,Q}$ is a unitary linear isomorphism.  
\end{lemma}

\proof
Uniqueness of the map follows from the fact that the remaining
maps in diagram (\ref{e: diagram with Rw and ev}) are linear isomorphisms.
Let $P$ be the unique parabolic subgroup in $\oppcPst$ such that $Q$ is $W(\fa)$-conjugate 
to $P.$ Let $v = v_Q$ and $v' = v_{Q'}.$ Then conjugation by $v' w$ maps
$Q$ to $P;$ hence $ v' w = v u$ for a suitable element $u \in N_{K_Q}(\fa).$  

We put $w' = w u^{-1}$ so that $w' = (v')^{-1} v$ and 
observe that $R_w$ equals the composition
 $R_{w'} R_u,$ with $R_u: C^{-\infty}(G/Q: \gs: \nu )_\chi \to C^{-\infty}(G/Q: u\gs: u\nu)_\chi$
and with $R_{w'}: C^{-\infty}(G/Q: u\gs: u\nu)_\chi \to C^{-\infty}(G/Q: w\gs: w\nu)_\chi.$
If $f \in  C^{-\infty}(G/Q: \gs: \nu )_\chi,$ then 
$$
\ev_Q  R_u f = f(v_Q u) =\gs(u)^{-1} \ev_Q(f).
$$ 
 From
Lemma \ref{l: gs u inv isom}
we know that $\tau: \eta \mapsto \gs(u)^{-1} \eta$ defines an isometric linear isomorphism 
$$
\tau:  \ccH^{-\infty}_{\gs, \chi_Q} \to \ccH^{-\infty}_{u\gs, \chi_Q}.
$$%
This gives us the first square of maps in the diagram  below, for every $\nu \in \faQdc.$
\begin{equation}
\label{e: split up diagram of Rv and ev}
\begin{array}{ccccc}\!\!\!\!\!\!
C^{-\infty}(G/Q : \gs : \nu)_\chi \!\!\!\!\! & \;\; {\buildrel R_{u} \over \longrightarrow} \;\;\;\; & \!\!\!\!\! 
 C^{-\infty}(G/Q : u\gs : \nu)_\chi    & {\buildrel R_{w'}  \over \longrightarrow} &  C^{-\infty}(G/Q': w\gs :w\nu)_\chi \!\! \!\! \!\! \!\! \!\! \\
\downarrow {\scriptstyle \scriptstyle \ev_Q}
 &  &
 \downarrow {\scriptstyle \ev_{Q}}
 & & \downarrow  {\scriptstyle \ev_{Q'}} \\
 \!\!\!\!\!\ccH_{\gs, \chi_Q }^{-\infty}& {\buildrel \tau  \over \longrightarrow} & 
  \ccH_{u\gs, \chi_Q}^{-\infty}\;\;\;\; &{\buildrel\cR_{w', Q}\over \longrightarrow} & \ccH_{w\gs,\chi_{Q'}}^{-\infty}.
\end{array}
\end{equation}
To understand the second square of maps, let $g \in C^{-\infty}(G/Q: u\gs: u\nu)_\chi.$ 
Then $\ev_{Q'} (R_{w'} g) = g (v' w') = g(v) = \ev_Q(g).$  
Since $w' = {v' }^{-1}v,$ 
it follows from  Lemma \ref{l: identities unitary for cHgs} that there is a unique linear isomorphism  $\cR_{w',Q}$ from  $\ccH^{-\infty}_{u \gs, \chi_Q}$ onto $\ccH^{-\infty}_{w\gs,\chi_{Q'}},$ which makes the second square of maps commutative. Furthermore, 
$\cR_{w', Q}$ is isometric. 

The map $R_{w, Q} : = R_{w', Q} \after \tau$ makes the diagram (\ref{e: diagram with Rw and ev})  commutative, and is the composition of two isometries, hence an isometry of
its own right.
\qed

\begin{cor}
\label{c: diagram composition Rw and j}
Let $Q\in \cP,$ $w \in N_K(\fa)$ and let $\cR_{w,Q}: \ccH^{-\infty}_{\gs, \chi_{Q}} \to 
\ccH^{-\infty}_{w\gs, \chi_{w Q w^{-1}}}$ be the isometry of Lemma \ref{e: diagram with Rw and ev}.
The following diagram commutes
for every $\nu \in \faQdc$:
\begin{equation}
\begin{array}{ccc}
C^{-\infty}(G/ Q : \gs : \nu)_\chi  & \;\;\;\;{\buildrel R_{w}  \over \longrightarrow} \;\;\;\;& 
C^{-\infty}(G/w Qw^{-1}: w\gs :w\nu)_\chi \\
\uparrow {\scriptstyle \scriptstyle j( Q , \gs, \nu)}
 &  &
 \uparrow {\scriptstyle j(w Q w^{-1}, w\gs, w\nu)}\\
 \!\!\!\!\!\!\!\!   (\ccH_{\gs }^{-\infty})_{\chi_{Q}}&{\buildrel \cR_{w,Q}  \over \longrightarrow} & 
 \!\!\!\!\!\!\!\!  (\ccH_{w\gs}^{-\infty})_{      \chi_{w  Qw^{-1}}            }\;\;\;\; 
\end{array}
\end{equation}
The map  $\cR_{w, Q}$ is an isometric linear isomorphism. 
\end{cor}

\medbreak
Before proceeding we list some properties of the  standard  intertwining operators
for two parabolic subgroups $Q_j \in \cP(A)$ $(j=1,2)$ with equal split components.
Let $Q_j = M_j A_j N_j$, be their Langlands decompositions, then $A_1 = A_2$ and $M_1 =M_2.$ We denote by $\gS(Q_j)$ the set of $\fa_j$-roots in $\fn_{Q_j}.$
For $R \in \R$ we define
\begin{equation}
\label{e: defi fa Q21}
\fad(Q_2|Q_1, R) := \{\nu \in \fa_{Q_1\iC}^* \mid \inp{\Re\nu}{\ga} > R\;\;(\forall \ga \in \gS(\bar Q_2)\cap \gS(Q_1))\}.
\end{equation}
Fix $\gs\in \widehat M_{1\ds}$ and $\nu\in \fad(\bar Q_2|Q_1, 0).$
The  standard  intertwining operator  $A(Q_2, Q_1 , \gs, \nu)$ from $C^\infty(G/Q_1 : \gs: \nu)$ 
to $ C^\infty(G/Q_2: \gs: \nu)$  is given by the usual integral formula
\begin{equation}
\label{e: integral formula A}
A(Q_2, Q_1 , \gs, \nu) f(x) = \int_{N_2 \cap \bar N_1} f(x \bar n)\; d \bar n
\end{equation}
for $f \in C^{\infty}(G/Q_1: \gs : \nu).$ 

We agree to equip any nilpotent subalgebra  $\fn_*$ of $\fg$  with the Riemannian 
inner product obtained by restriction of the positive definite inner product
$ - B(\dotvar, \theta( \dotvar)).$  The associated 
analytic subgroup $N_*$ of $G$ is equipped with the associated  bi-invariant unit Haar measure. 
As a result, for every $k \in K$ the map $n \mapsto knk^{-1}$ is measure preserving 
from $N_*$ to $k N_* k^{-1}.$ 
In particular, all  standard  intertwining operators will be normalized in this way.
Then for all $w \in N_K(\fa),$ all $\nu \in \fad(Q_2|Q_1,0)$ and all
$f \in C^\infty(G/Q_1: \gs : \nu) $ we have
\begin{equation}
\label{e: Weyl twisted A}
A(Q_2, Q_1, \gs, \nu) f =   R_w^{-1} A(w Q_2 w^{-1}, w Q_1 w^{-1}, w\gs, w \nu) R_w  f.
\end{equation}
Indeed this follows from the integral formula (\ref{e: integral formula A}) in view of the 
normalization of measures specified in the above. 

We shall now describe the well known meromorphic continuation of 
the intertwining operator in terms of the compact pictures of the induced representations,
where we exploit the topological linear isomorphism 
$$
C^\infty(G/Q_j : \gs :\nu) \simeq C^\infty (K/K_{Q_j}: \gs_{Q_j}),
$$
induced by restriction to $K.$ 
By transfer under this isomorphism, the left regular representation $L$ in the
first space becomes a $\nu$-dependent representation $\pi_{Q_j, \gs, \nu}$ 
of $G$ in the second space. The operator $A(Q_2, Q_1, \gs, \nu)$ can now
be viewed as a continuous linear operator of  Fr\'echet spaces 
$C^\infty(K/K_{Q_1}: \gs)\to C^\infty(K/ K_{Q_2}: \gs).$ 
The dependence of the intertwining operator of  $\nu \in \fa_{1\iC}^*$ is known to be meromorphic, by  \cite{VW}, 
with singular locus $\cS(Q_2, Q_1, \gs)$ a  locally finite union of affine root hyperplanes of the form
$\inp{\nu}{\ga} = c,$ with $\ga \in \gS(\fg, \fa_j)$ and $c \in \R.$  As a result, the equality
(\ref{e: Weyl twisted A}) is valid as equality of meromorphic functions 
on $\fa_{1\iC}^*$ with values in  $C^{\infty}(K/K_{Q_1}:\gs_{Q_1}).$

Via the equivariant sesquilinear pairing $\inp{\dotvar}{\dotvar}$ of the space $C^\infty(Q:\gs:\nu)$ with 
$C^\infty (Q: \gs -\bar \nu)$ by integration over $K/K_{Q}$ we may 
embed the first space in the conjugate continuous linear dual of the second, denoted
$C^{-\infty}(Q: \gs: \nu).$ In view of the formula 
$A(Q_1, Q_2,  \gs: -\overline \nu)^* = A(Q_2, Q_1, \gs, \nu)$ one sees that the intertwining
operator $A(Q_2, Q_1, \gs ,\nu)$ has a continuous linear extension to a continuous
linear intertwining operator $C^{-\infty}(G/Q_1: \gs : \nu) \to C^{-\infty}(G/Q_2: \gs : \nu)$ 
for non-singular values of $\nu.$ In  \cite{BSestc} this extended operator is shown to 
depend meromorphically on $\nu \in \fa_{1\iC}^*$ in the following way, in terms of the compact
picture.

\begin{lemma}
\label{l: mero prop A}
For every $R \in \R$ there exists a polynomial function $q\in P(\fa_{1\iC}^*),$
with zero set contained in $\cS(Q_2, Q_1, \gs),$ and a constant $r \in \N$
such that for every positive integer $s$ the assignment
$$
\nu \mapsto p(\nu) A(Q_2 , Q_1 , \gs , \nu)
$$ 
defines a holomorphic function on $\fad(Q_2|Q_1, R)$
with values in the Banach space $B(C^{-s}, C^{-s-r})$ of bounded linear maps
$C^{-s}(K/P_1: \gs_{K_{P_1}}) \to C^{-s-r}(K/P_2: \gs_{K_{P_2}}).$ 
\end{lemma}

By continuity and density it now readily follows that (\ref{e: Weyl twisted A}) is valid for all
$f\in C^{-\infty}(G/Q_1: \gs :\nu), $ provided $\nu \in  \fa_{1\iC}^*\setminus \cS(Q_2, Q_1, \gs).$ 

\begin{cor} 
Let $\Omega$ be a bounded open subset of $\fa_{1\iC}^*.$ Then there exists an $r \in \N$ such that
for every $s \in \N$ the assignment $\nu \mapsto A(Q_2, Q_1,\gs, \nu)$ defines a meromorphic function on $\Omega$
with values in  in the Banach space $B(C^{-s}, C^{-s-r})$ of bounded linear maps
$C^{-s}(K/P_1: \gs_{K_{P_1}}) \to C^{-s-r}(K/P_2: \gs_{K_{P_2}}).$ 
\end{cor}

\begin{defi}
\label{d: order of sm loss}
The smallest $r \in \N$  for which the above is valid will be called the order (of smoothness loss) of 
the family $A(Q_2, Q_1, \gs, \dotvar)$ over $\Omega.$ 
\end{defi}
\begin{lemma}
Let $\Omega \subset \faQdc$ be open and suppose that $f: \Omega \to C^{- \infty}(K/K_Q:\gs_Q)$ is a  holomorphic function such that 
\begin{equation}
\label{e: transformation property under N 0}
\pi_{Q, \gs, \nu}(n) f_\nu = \chi(n) f_\nu, \qquad (\nu \in \Omega, n \in N_0).
\end{equation}
Then $\nu \mapsto \ev_Q (f_\nu)$ 
is a holomorphic function $\Omega \to \ccH^{-\infty}_{\gs, \chi_Q}.$ 
\end{lemma}

\proof
Fix $\nu_0 \in \Omega.$  Put $I^{\pm \infty} = C^{\pm \infty}(K/K_Q:\gs_Q).$
Furthermore, put $V: =  \ccH^{-\infty}_{\gs, \chi_Q}.$ 
Then $j(\nu): = j(Q,\gs, \nu)$ is a linear map $V \to I^{-\infty}$ which depends holomorphically on $\nu.$
From $\ev _Q \after  j(\nu) ={ \rm id}_V$ it follows that $j(\nu)$ is injective for every $\nu \in \Omega.$ 

Let $E:= j(\nu_0)(V)$ and let $E^\perp$ denote the annihilator of $E \subset I^{-\infty}$ 
 in $I^\infty.$ 
We fix a linear subspace $W_0 \subset  \overline{I^{\infty}}$ which is complementary to $E^\perp.$ 
Then the restriction map $r: I^{-\infty} \to W_0'$, $\xi \mapsto \xi|_{W_0}$ is continuous 
linear and restricts to a linear isomorphism from $E$ onto $W_0'$.

Since $r$ is continuous linear, the function  $J: \nu \mapsto r\after j(\nu)$ 
is holomorphic with values in $\Hom(V, W_0').$ As $J(\nu_0)$ is bijective linear $V \to W_0',$ 
the same is true for $\nu$ in a  sufficiently small open neighborhood $\Omega_0 \ni \nu_0$ in $\Omega.$ 
Furthermore, the function $\nu \mapsto J(\nu)^{-1}$ is holomorphic on $\Omega_0$ with values 
in $\Hom(W_0', V).$ 

Let $f$ satisfy the hypotheses, and put $a(\nu) = \ev_Q f_\nu.$ Then $f_\nu =  j(\nu)(a(\nu))$ for $\nu \in \Omega.$ Furthermore, since $\nu \mapsto f_\nu$ is holomorphic as a map with values in $I^{-\infty}$
it follows that  $\nu \mapsto r (f_\nu) = r \after j(\nu) [a (\nu)]$ is holomorphic in $\nu$ with values in 
$W_0'.$ It follows that $\nu \mapsto a(\nu) = J(\nu)^{-1} [ r(f_\nu) ]$ is holomorphic on $\Omega_0$ with values
in $V.$ Hence, $\ev_Q\after f $ is a holomorphic function $\Omega_0 \to V.$ 
\qed
\medbreak
Let $Q_1, Q_2\in \cP$ have equal split components: $\fa_{Q_1} = \fa_{Q_2}.$ 
In analogy with the theory of symmetric spaces, we define, 
for a regular point $\nu\in \fa_{Q_1\iC}^*$ of $A(Q_2,  Q_1, \gs, \dotvar),$ the linear map 
$B(Q_2, Q_1,\gs , \nu): \ccH^{-\infty}_{\gs, \chi_{Q_1}} \to \ccH^{-\infty}_{\gs, \chi_{Q_2}}$ 
by 
\begin{equation}
\label{e: defi B by A j}
B( Q_2, Q_1,\gs , \nu) \eta := \ev_{Q_2}A( Q_2, Q_1, \gs, \nu) j(Q_1, \gs, \nu) \eta, \qquad (\eta \in 
\ccH^{-\infty}_{\gs, \chi_{Q_1}}).
\end{equation}
Here $\ev_{Q_2} = \ev_{v_2},$ with $v_2 = v_{Q_2},$  so that $v_2 Q_2 v_2^{-1}$ belongs to $\bar \cP_{\rm st}$.

\begin{lemma}
\label{l: A B and j}
\ Let $Q_1, Q_2 \in \cP$ have the same split component, and let 
$\gs \in \widehat M_{Q_1,\ds}.$ 
\begin{enumerate}
\itema
The function $B( Q_2, Q_1, \gs, \dotvar ): \fa_{1\iC}^* \to
\Hom(\ccH^{-\infty}_{\gs, \chi_{Q_1}},   \ccH^{-\infty}_{\gs, \chi_{Q_2}})$ is meromorphic.
\itemb 
For every $\eta \in \ccH^{-\infty}_{\gs , \chi_{Q_1}},$ 
\begin{equation}
\label{e: A and j and B} 
A(Q_2, Q_1, \gs, \nu) j(Q_1, \gs, \nu) \eta = j(Q_2, \gs, \nu) B(Q_2,Q_1, \gs, \nu) \eta
\end{equation}
as an identity of meromorphic 
$C^{-\infty}(K/K_{Q_2}: \gs_{Q_2})$-valued functions of $\nu \in \fa_{Q_1 \iC}^*$.
\end{enumerate}
\end{lemma}

\proof
For (a), assume that $\nu_0 \in  \fa_{Q_1\iC}^*$ and let $\Omega$ be a sufficiently small bounded neighborhood 
of $\nu_0$ in  $\fa_{Q_1\iC}^*.$ Then there exists a holomorphic function $\gf: \fa_{Q_1\iC}^* \to \C$ 
such that the function $$
\tilde A: \nu \mapsto \gf(\nu) A(Q_2, Q_1, \gs, \nu)
$$
 is holomorphic on $\Omega$  with values
in $\End(C^{-\infty}(K/K_{Q_1} : \gs_{Q_1}))$ in the sense that there exists an $r > 0$ such that 
for every $k > 0$ the function $\tilde A$ defines a holomorphic function on $\Omega$ with values in 
the Banach space $B(C^{- k}, C^{-k -r})$ of bounded linear operators from 
$C^{-k}(K/K_{Q_1}:\gs_{Q_1}) \to C^{-k-r}(K/K_{Q_2}:\gs_{Q_2}).$ Furthermore, 
by holomorphy of $j(Q_1, \gs, \dotvar)$  and  boundedness of $\Omega,$ there exists $s > 0$ such that $j( Q_1, \gs, \dotvar )$ defines a holomorphic function 
$\Omega \to B(\ccH_{\gs, \chi_{Q_1}}^{-\infty}, C^{- s}).$ 
Hence, 
$
f: \nu \mapsto  A(\nu) j(Q_1, \gs, \nu) \eta
$ 
defines a holomorphic function on $\Omega$ with values in $C^{-s-r}(K/K_{Q_2}:\gs_{Q_2}).$ 
By the transformation property under $N_0$ of $j( Q_1, \gs, \nu)\eta$ and the equivariance of $A(Q_2, Q_1, \gs, \nu)$ it follows that $f_\nu$ satisfies the transformation property 
(\ref{e: transformation property under N 0}) with $Q = Q_2.$  
By application of Lemma 2.4 it now follows that $\ev_{2} \after A(\nu) ( Q_1, \gs, \nu)$ 
is holomorphic in $\nu \in \fa_{Q_1\iC}^* $ with values in $\ccH_{\gs, \chi_{Q_2}}^{-\infty}$. 
Hence $\nu \mapsto \gf(\nu) B(Q_2, Q_1, \gs, \nu) \eta$ is holomorphic, and assertion 
(a) follows. 

The validity of (b) is checked by applying the evaluation $\ev_{Q_2}$ to both sides
of the equation (\ref{e: A and j and B}) and using that $\ev_{Q_2}$ is a bijection from
$C^{-\infty}(G/ Q_2: \gs: \nu) _\chi$ onto $\ccH^{-\infty}_{\sigma, \chi_{Q_2}},$ for every $\nu \in \fa_{Q_2 \iC}^*.$ 
\qed

%\end{document}

\begin{cor}
\label{c: commutator Rw and B}
Let $w \in N_K(\fa)$ and let $Q_1,Q_2 \in \cP$ have common split component. 
Let $\gs \in \widehat M_{Q_1 \ds}.$ Then 
$$
\cR_{w, Q_2}\after B(Q_2, Q_1, \gs, \nu) 
= B(w Q_2 w^{-1}, w Q_1 w^{-1}, w\gs, w\nu) \after \cR_{w, Q_1}
$$ 
as meromorphic functions of $\nu \in \fa_{1\iC}^*$ with values in $\Hom(\ccH^{-\infty}_{\gs, \chi_{Q_1}},
\ccH^{-\infty}_{\gs, \chi_{Q_2}}).
$
\end{cor}
%\comment{until here, Sept 3, 2025}
\proof
From (\ref{e: Weyl twisted A}) for generalized functions it follows, for generic $\nu \in \fa_{Q_1\iC}^*$ that
$$ 
R_w A( Q_2,  Q_1, \gs, \nu)j( Q_1,\gs, \nu) =
A(w Q_2w^{-1}, w  Q_1 w^{-1}, w\gs, w\nu) R_w j( Q_1,\gs, \nu).
$$%
Applying $\ev_{w Q_2 w^{-1}}$  to  both sides of the equation and using Lemma \ref{l: diagram composition Rw and j} 
and Corollary \ref{c: diagram composition Rw and j}
we
find
the asserted equality.
\qed

\begin{cor}
Let $Q_1, Q_2 \in \cP$ have the same split component.
The following statements are equivalent, for any meromorphic function $\eta: \fa_{Q_1\iC}^* \to \C$ 
and any $w \in N_K(\fa).$ 

\begin{enumerate}
\itema
$B( Q_2, Q_1, \gs, -\bar \nu)^*B( Q_2,  Q_1, \gs, \nu) = \eta(\nu) I$ \;\;
{\rm for generic }\;\;$\nu\in \fa_{Q_1\iC}^*,$
\itemb 
$B(w Q_2w^{-1} , w  Q_1 w^{-1} , w\gs, - w\bar \nu)^*B(w  Q_2w^{-1},w  Q_1w^{-1}, w\gs, w\nu) = \eta(\nu) I\;\;\;$ 
for generic $\nu\in \fa_{Q_1\iC}^*.$
\end{enumerate}
\end{cor}

\proof 
Let $\nu \in \fa_{Q_1\iC}^*$ be a regular value for each of the finitely many functions involved.
It follows from Lemma \ref{l: diagram composition Rw and j}  that  $\cR_{w,Q}: \ccH_{\gs, \chi_Q}^{-\infty} \to \ccH_{w\gs,\chi_{wQw^{-1}}}^{-\infty}$ is unitary. The result now follows by a simple argument, using Corollary 
\ref{c: commutator Rw and B}. 
\qed

\section{The B-matrices, reduction arguments}
\label{s: Bmatrix, reduction}
Let $P,Q \in \cP$  have the same split component. 
Then $M_Q = M_P.$ If $\gs \in \dMPds$ then for generic $\nu \in \faPdc$ the composition 
 $A(P, Q, \gs, \nu) A(Q, P, \gs, \nu)$ of  standard  intertwining operators is a self intertwining
 operator of $\Ind_P^G(\gs \otimes \nu \otimes 1).$ Since the latter representation is irreducible
 for generic $\nu,$ it follows that
 \begin{equation}
 \label{e: definition eta}
 A(P, Q, \gs, \nu) A(Q, P, \gs, \nu) = \eta(Q, P, \gs, \nu) {\rm id} , \qquad (\nu \in \faPdc),
 \end{equation}
 for a unique meromorphic function $\eta(Q,P, \gs, \nu).$ It is easy to see that
 $\eta(Q, P, \gs, \nu) = \eta(P, Q, \gs, \nu).$ 
 Furthermore, suppressing ${\rm id}$, 
 $$ 
 A(Q,P ,\gs, -\bar \nu)^*A(Q,P ,\gs, \nu) = \eta(Q, P, \gs, \nu), \qquad (\nu\in \faPdc).
 $$ 
For this and other properties of $\eta,$ we refer to  \cite{KS2}. 
In the course of this and the next sections we will prove the following
manifestation of the Maass -- Selberg relations.

\begin{thm}
\label{t: MS for B}
Let $P,Q \in \cP$ and suppose that $\faP = \faQ.$  
If $\gs \in \dMPds$, we have the following identity of 
meromorphic functions $\faPdc \to \End(\ccH_{\gs, \chi_P}^{-\infty}):$
\begin{equation}
\label{e: Maass-Selberg for B}
B( Q,  P, \gs, -\bar \nu)^*B( Q, P,\gs, \nu) = \eta( Q, P, \gs, \nu),\qquad (\nu \in \faPdc).
\end{equation}
\end{thm}

This result will be proven in the course of the next sections, through reduction 
to a basic setting where $G$ has compact kernel and $P$ is a maximal standard parabolic subgroup.
\begin{rem} 
It immediately follows from  (\ref{e: definition eta}) that $B( P,  Q, \gs, \nu) B( Q,  P , \gs, \nu) = \eta( Q, P, \gs, \nu).$
Therefore, the identity (\ref{e: Maass-Selberg for B}) is valid if and only if
$$
B( Q,  P, \gs, -\bar \nu)^* = B( P,  Q, \gs, \nu), \qquad (\nu \in \faPdc).
$$
\end{rem}

\begin{rem} 
Note that for a non-cuspidal parabolic subgroup $P \in \cP$ the group $M_P$ has no discrete series,
so that $\dMPds= \emptyset.$ This means that for trivial reasons the assertion of 
Theorem \ref{t: MS for B} is automatically fulfilled for $Q\in \cP$ with $\faQ = \faP.$ 
\end{rem}
 
We will write $G = M_G A_G$ for the Langlands decomposition
of $G,$ viewed as a parabolic subgroup of $G.$ 
Here $A_G:= \exp \fa_G,$ where $\fa_G$ is the intersection
of the center of $\fg$ with $\fp.$ The group $G$ has compact center if and only if $A_G = \{e\}.$
The group $M_G$ is also denoted ${}^\circ G$ and equals the intersection 
of the kernels $\ker \xi$ where $\xi$ runs over the collection $X(G)$of multiplicative characters $G \to \R_+.$ 
Write ${}^\circ \!A = A \cap {}^\circ \!G.$ Then ${}^\circ G$ is of the Harish-Chandra class, with Iwasawa docomposition
$$
{}^\circ G = K\, {}^\circ \!A N_0.
$$
 Let ${}^\circ{\cP}$ denote the finite set 
of parabolic subgroups of ${}^\circ G$ containing ${}^\circ A.$ Then the map $Q \mapsto {}^\circ Q := Q \cap {}^\circ G$ is a bijection $\cP \to {}^\circ \cP,$ with inverse ${}^\circ Q \mapsto 
{}^\circ Q A_G.$ Note that every $v \in N_K(\fa)$  leaves 
${}^\circ A$ invariant. We define ${}^\circ v: {}^\circ\cP \to N_K({}^\circ \!A)$ 
by ${}^\circ v_{{}^\circ Q}  = {}^\circ\! v_Q|_{{}^\circ A},$ for $Q\in \cP.$ 

We note that for $Q\in \cP,$ ${}^\circ Q$ has the Langlands decomposition 
${}^\circ Q = M_Q {}^\circ A_Q N_Q,$ where ${}^\circ A_Q= A_Q \cap {}^\circ G.$
Thus, if $P, Q\in \cP$ then $P$ and $Q$ have the same  split components as
elements of $\cP$ if and only if ${}^\circ P$ and $ {}^ \circ Q$ have the same split components as elements of ${}^\circ \cP.$ 

Furthermore, if $Q \in \cP,$ $M_Q\subset {}^\circ G.$ Therefore, if 
$\gs \in \dMQds,$ we may identify the space $\ccH_{\gs, \chi_Q}^{-\infty}$ 
with the similar space for the pair ${}^\circ G, {}^\circ Q.$ The same is true for 
the inner products on these spaces. 

Accordingly, we may define the $B$ matrices for ${}^\circ G,$ denoted 
${}^\circ B({}^\circ Q, {}^\circ P, \gs, \mu)$ in the
obvious fashion as endomorphisms of $\ccH_{\gs, \chi_P}^{-\infty}, $
for $P, Q\in \cP$ with $\faP = \faQ$ and for $\gs \in \dMPds$ and 
$\mu \in {}^\circ\faPdc.$
Then $ \mu \mapsto {}^\circ B({}^\circ Q, {}^\circ P, \gs, \mu)$
is a meromorphic function on ${}^\circ \faPdc$ with values 
in $\End (\ccH_{\gs, \chi_Q}^{-\infty}).$ 
We agree to write ${}^\circ \nu:= \nu |_{{}^\circ \fa_{P\iC}}$ for $\nu \in \faPdc.$ 

\begin{lemma}
\label{l: comparison with cc}
Let notation be as in the above text. Then for generic $\nu \in \faPdc,$ 
\begin{enumerate}
\itema $\eta(Q,P, \gs, \nu) =  {}^\circ \eta({}^\circ Q, {}^\circ P, \gs, {}^\circ \nu);$ 
\itemb $B( Q, P, \gs, \nu) =  {}^\circ B({}^\circ Q, {}^\circ P, \gs, {}^\circ \nu).$
\end{enumerate}
\end{lemma}

\proof
The proof is straightforward, but a bit tedious. Details are left to the reader.
\qed

\begin{cor}
\label{c: reduction to cc} 
The assertions of Theorem \ref{t: MS for B} are valid for $G$ if and only if they are valid for
${}^\circ G.$
\end{cor}

Let $P \in \cP.$ A  root of $P$ is defined to be  a non-trivial linear functional $\ga \in \faP^*$ such that
the space $\fg_\ga := \cap_{H\in \faP} \ker(\ad H -\ga(H))$ is contained in $\fn_P.$ Equivalently this means
that $\ga$ is the restriction of a root $\gb \in \gS(\fa)$ with $\fg_\gb \subset \fn_P.$ 
The set of $P$-roots is denoted by $\gS(P).$ 
We note that $\faPp$ is the set of points $H \in \faP$ such that $\ga(H) > 0$ for all $\ga \in \gS(P).$ 

A $P$-root $\ga \in \gS(P)$ is called reduced if the multiples of $\ga$ in $\gS(P)$ are all of the form
$c \ga,$ with $c \geq 1.$

For two parabolic subgroups $P,Q \in \cP$  with the same split component, a $P$-root $\ga$ is said 
to separate $\faPp$ and $\faQp$ if the sign of $\ga$ on $\fa_{Q}$ is negative, or equivalently, $\ga \in \gS(\bar Q).$ 

The distance $d(P,Q)$ 
is defined to be the number of reduced $P$-roots that separate $\faPp$ and $\faQp.$ 
$P$ and $Q$ are said to be adjacent if $d(P,Q) = 1.$ Equivalently this means
that all roots in $\gS(P) \cap \gS(\bar Q)$ are proportional. 

If $P,Q$ have the same split component and are different, then there is a parabolic subgroup 
$R \in \cP$ with split component $\faR = \faP$ such that
$P$ and $R$ are adjacent, and $d(R,Q) < d(P,R).$ 
It is well known, see e.g.  \cite[Cor.~7.7]{KS2}, that 
for $\gs \in \dMPds$ one has in this case that 
\begin{equation}
\label{e: prod deco A}
A(Q ,P,  \gs, \nu) =   A(Q, R, \gs, \nu) A(R, P, \gs, \nu),
\end{equation}
for generic $\nu \in \faPdc.$ 

\begin{lemma}
\label{l: product deco B}
In the above setting, $B(Q ,P,  \gs, \nu) =   B(Q, R, \gs, \nu) B(R, P, \gs, \nu).$
\end{lemma} 
\proof
This follows from \ref{e: prod deco A} by application of Lemma \ref{l: A B and j} and 
(\ref{e: defi B by A j}).
\qed

\begin{lemma}
\label{l: reduction B to adjacent}
If the identity of Thm. \ref{t: MS for B} holds for all adjacent $Q,P\in \cP$ then it holds
for all $Q,P \in \cP$ with $\faP = \faQ.$ 
\end{lemma}

\proof
We assume that $G$ is fixed and that the identity of Thm. \ref{t: MS for B} holds for all 
$P,Q\in \cP$ with $d(P,Q) = 1.$ 
Arguing by induction on $d(P,Q)$ we will show that the identity holds for all $P, Q$ 
with the same split component.

Let $k \geq 1$ and suppose that the identity holds for $P,Q \in \cP$ with $d(P,Q) \leq k.$ 
Assume now that $d(P,Q) = k+1.$ Then there exists $R$ as in the text leading to (\ref{e: prod deco A}). 
By induction we know that
$$
B(Q,R, \gs, -\bar \nu)^* B(Q,R , \gs, \nu) = \eta(Q,R, \gs, \nu).
$$ 
Therefore, 
\begin{eqnarray*}
\lefteqn{
B(Q, P, \gs, -\bar \nu)^*B(Q, P, \gs, \nu) = }\\
&=& B(R, P, \gs,- \bar \nu)^* B(Q, R, \gs, -\bar \nu)^* \circ   
B(Q, R, \gs, \nu) B(R, P, \gs, \nu) \\
&=& B(R, P, \gs,- \bar \nu)^* \eta(Q,R, \gs, \nu) B(R, P, \gs, \nu)\\
&=& \eta(Q,R, \gs, \nu)\cdot B(R, P, \gs,- \bar \nu)^* B(R, P, \gs, \nu).
\end{eqnarray*}
Since $d(P,R) = 1,$ the expression on the last entry of the array equals the product
$\eta(Q,R, \gs, \nu) \eta(R,P, \gs, \nu).$ In turn, as a consequence of (\ref{e: prod deco A}),
 this product equals
$ \eta(Q,P, \gs, \nu).$ 
\qed

\begin{lemma}
\label{l: MS and v conjugacy}
Suppose the identify of Theorem \ref{t: MS for B} holds for $G, P, Q, \gs$ and all $\nu \in \faPdc.$ If $v\in N_K(\fa)$,
then the identity of the theorem also holds with
$G, vPv^{-1}, vQv^{-1}, v \gs$ in place of $G,P,Q,\gs$ respectively, and all $\nu \in v \faPdc.$ 
\end{lemma}

\proof
Let the hypothesis be fulfilled, then it suffices to prove the following identity for all $\nu \in \faPdc$,
\begin{equation}
\label{e: conjugate of MS rels for B}
B(v Q v^{-1}, v P v^{-1},  v\gs, -v \bar \nu)^* B(v Q v^{-1}, v P v^{-1},  v\gs, v\nu) = \eta(v Q v^{-1}, v P v^{-1},  v\gs, v  \nu).
\end{equation}
Using Cor. \ref{c: commutator Rw and B} and the unitarity of $ R_{v, Q}$ and $R_{v,P}$ we may rewrite the expression on the left-hand side of \ref{e: conjugate of MS rels for B} as
\begin{eqnarray*}
\lefteqn{
\!\!  \!\!  \!\!  \!\!  \!\!  \!\!  \!\!  \!\!   R_{v,P}^* B(Q,P, \gs, -\bar \nu)^* R_{v, Q}^* R_{v,Q} B(Q , P,  \gs, \nu) R_{v,P} =}\\
&=& R_{v,P}^* B(Q,P, \gs, -\bar \nu)^* B(Q , P,  \gs, \nu) R_{v,P} \\
&=&\eta(Q,P, \gs, \nu).
\end{eqnarray*} 
The result now follows from the observation that 
$\eta(v Q v^{-1}, v P v^{-1},  v\gs, v \nu) =  \eta(Q, P ,  \gs, \nu), $ in view of  (\ref{e: Weyl twisted A})  and 
the definition of $\eta(Q,P, \gs, \nu).$ 
\qed

\section{Reduction to maximal parabolic subgroups}
\label{s: parabolics}
In this section we will discuss a method that will allow reduction
of the proof of the Maass-Selberg relations to those for maximal parabolic subgroups of lower dimensional
groups of the Harish-Chandra class. 

To prepare for this, we will first discuss well known aspects of the structure of parabolic subgroups. We briefly write $\gS(\fa):= \gS(\fg, \fa).$ 
If $P$ is a parabolic subgroup of $G$, we denote its Langlands 
decomposition by $P = M_P A_P N_P.$ The collection of $\fa_P$-roots in $\fn_P$ is defined as in the text below Cor. \ref{c: reduction to cc}  and denoted by $\gS(P).$ 

The positive chamber $\fa_P^+$ consists
of the points $X \in\faP$ such that $\ga(X) > 0$ for all $\ga \in \gS(P).$ 

Given a point $X \in \fa$ we put
$$ 
 \gS^+(\fa, X):= \{\ga \in \gS(\fa)\mid  \ga(X) >  0\}.
 $$ 
We note that for a root $\ga\in \gS(\fa),$ one has $\ga(X) < 0 \iff -\ga \in \gS^+(\fa, X)$ and
$\ga(X) = 0$ iff $\pm \ga \notin \gS^+(\fa, X).$ 

We define the equivalence relation $\sim$ on $\fa$ by 
\begin{equation}
\label{e: equivalence determining parabs}
X\sim Y \iff \gS^+(\fa, X) = \gS^+(\fa, Y),\qquad (X , Y \in \fa).
\end{equation}
For $X\in \fa$ we define the subspace  $\fp_X = \fm_X \oplus \fa_X \oplus \fn_X$ 
of $\fg$ by 
\begin{equation}
\label{e: defi psa fp X}
\fp_X = \fm + \fa + \oplus_{\ga \in \gS^+(\fa, X)} \;\;\fg_\ga.
\end{equation}
Then, clearly, $\fp_X = \fp_Y \iff X \sim Y.$ 

Suppose now that $\gS^+$ is a positive system for 
$\gS(\fa)$ and $\fa^+$ the associated open positive chamber in $\fa.$ 
If $X \in \cl(\fa^+)$ then one readily verifies that $\fp_X = \fm_X \oplus \fa_X \oplus \fn_X$ is the  standard  parabolic subalgebra  $\fp_F$ with $F$ the collection of simple roots in $\gS^+$ vanishing on $X.$ 
The indicated Langlands decomposition is determined by 
$$ 
\fm_X ={}^\circ {\fz}_\fg(X),\;\; \fa_X^+ = [X], \;\; \fn_X = \oplus_{\ga \in 
\gS^+(\fa, X)}  \;\; \fg_\ga,
$$
where $[X]$ denotes the class of $X$ for the equivalence relation $\sim.$
Write $\spec(X)$ for the spectrum of $\ad(X) \in \End(\fg),$ and $\spec(X)_+$ for 
its positive part.
Then it is readily checked that 
$$
\fp_{X} = \ker \ad (X) + \oplus_{\gl \in \spec(X)_+} \ker (\ad(X) - \gl I).
$$ 
By using $\Ad(K)$ conjugacy, we see that this definition gives  a parabolic subalgebra of
$\fg$ for any $X\in \fg$ with $\theta(X) = - X.$ Furthermore, its Langlands components 
are given by 
$$ 
\fm_X + \fa_X = \fz_\fg(X),\quad \fn_X =  \oplus_{\gl >0} \ker \,(\ad(X) - \gl I).
$$
As usual, the parabolic subgroup with algebra $\fp_X$ is defined
by $P_X := N_G(\fp_X).$ Its Langlands components are given by 
$$ 
M_X = {}^\circ Z_G(X), \;\;A_X = \exp (\fa_X) , \;\;N_X = \exp (\fn_X).
$$

The following lemma will be used repeatedly in the sequel.

\begin{lemma}
\label{l: root deco principle}
If $\fb$ is an abelian subspace of $\fs$ and 
$X \in \fb$ then for each $\gl\in\spec(X)$
the associated eigenspace $\ker(\ad(X) - \gl)$ in $\fg$ equals the direct sum 
of the $\ad(\fb)$-weight spaces $\fg_\gb$ for $\gb \in \gS(\fb) \cup \{0\}$ with 
$\gb(X) = \gl.$ 
\end{lemma}
\proof
For each $X\in \fs,$ the endomorphism $\ad X$ of $\fg$ is symmetric with respect
to the inner product $\inp{\dotvar}{\dotvar}$ hence semisimple with real eigenvalues.
The proof is now straightforward.
\qed

Let $P\in \cP$ have split component $\fa_P.$   The complement
of the union of the finitely many hyperplanes $\ker \gb$ with $\gb \in \gS(\fg, \fa_P)$
consists of finitely many convex polyhedral components, called the chambers of $\fa_P.$
These chambers are readily seen to be equivalence classes for the equivalence relation
$\sim$ given by (\ref{e: equivalence determining parabs}). For each chamber $[Y]$  the associated parabolic subgroup
$P_{[Y]}:= P_Y$ has split component $\fa_P$ and positive chamber $\fa_{P_{[Y]}}^+ = [Y].$
Conversely, for each parabolic subgroup $Q$ with split component $\fa_Q= \fa_P$
the positive chamber $\fa_Q^+$ is a chamber in $\fa_P.$ 

Two parabolic subgroups $P,Q \in \cP$ with the same split components
are said to be adjacent if their positive chambers are separated by precisely one hyperplane
from the collection of hyperplanes $\ker \ga \subset \fa_P$ for $\ga \in \gS(\fn_P, \fa_P).$ 
A root $\gb \in \gS(\fg, \fa_P)$ is said to be reduced iff all its real multiples in $ \gS(\fg, \fa_P)$
are of the form $c\gb$ with $|c| \geq 1.$ Thus, $P$ and $Q$ are adjacent if there is
a unique reduced $\fa_P$-root $\ga \in \gS(\fn_P, \fa_P)$ such that
$\ga < 0$ on $\fa_Q^+.$ Note that $-\ga$ is the unique reduced root in $\gS(\fn_Q, \fa_Q)$ 
which is negative on $\fa_P.$ In this situation it is easy to see that
\begin{equation}
\label{e: relation roots adjacent psg} 
\fn_P \cap \bar \fn_Q = \fn_\ga := \oplus_{c \geq 1}\;\; \fg_{c \ga}, 
\end{equation} 
where the summation is over the real $c \geq 1$ such that $c \ga$ is a root of $\fa_P.$ 
We note that the reduced root $\ga$ has the property that
$\ker \ga \cap \cl(\faPp)$ has non-empty interior in $\ker \ga.$ 
Conversely if such a reduced root
$\ga$ is given, then there is a unique parabolic subgroup $Q$ with split component
$\faP$ that is adjacent to $P.$ It is determined by the requirement that 
(\ref{e: relation roots adjacent psg}).

Let $P \in \cP.$ We assume that $P$ is {\em not maximal} {\ } and 
fix a reduced root $\ga \in \gS(\fn_P, \faP)$ such that $\ker \ga \cap \cl(\faPp)$
has non-empty interior as a subset of the hyperplane $\ker\ga$ in $\faP.$ 
A point $X$ in this interior will be called $(P, \ga)$-regular if
it has the property that for all $\gb \in \gS(\fg, \faP),$
$$ 
\gb(X) = 0 \implies \gb|_{\ker \ga} = 0 . 
$$ 
As we explained, the pair $(P,\ga)$ uniquely determines 
an adjacent $Q.$ A point $X$ in the interior of 
$\ker \ga \cap \cl(\faPp) =  \ker \ga\cap \cl(\faQp)$
is $(P, \ga)$-regular if and only if it is $(Q, - \ga)$-regular. 
We may select a $(P, \ga)$-regular point $X$ in the interior of $\ker \ga \cap \cl (\faPp).$ 
We fix such an $X$ and define the parabolic subalgebra $\fp_X$  as in (\ref{e: defi psa fp X}).
The corresponding parabolic subgroup $N_G(\fp_X)$ has the 
Levi decomposition $M_{1X} N_X.$ We put
$G^{(\ga)}:= M_{1X},$ ignoring the precise dependence on the choice of
the $(P,\ga)$-generic element $X.$ 

\begin{lemma}
The Lie algebra of $G^{(\ga)}$ is given by 
\begin{equation}
\label{e: defi fg ga}
\fg^{(\ga)} = \bar \fn_\ga \oplus \fm_{1P} \oplus \fn_\ga.
\end{equation}
\end{lemma}

\proof
By Lemma \ref{l: root deco principle},
$
\fg^{(\ga)} = \fm_{1X} = \oplus_{\gb}\; \fg_\gb
$
where the sum is taken over the $\gb\in \gS(\faP)\cup \{0\}$
for which either (a): $\gb|_\faP = 0$ or (b): $\gb|_{\faP}\neq 0$ and $\gb(X) = 0.$
Condition (a) is equivalent to $\fg_\gb \subset \fm_{1P}.$ Condition (b) is equivalent to
$\gb|_{\ker \ga} = 0$ which in turn is equivalent to $\gb = c \ga$ for a constant $c \in \R.$ 
The validity  of (\ref{e: defi fg ga}) follows.
\qed

\begin{lemma}
\label{l: deco P in ga X}
$
\fn_P =\fn_\ga \oplus \fn_X.
$
\end{lemma}

\proof
$\fn_P$ is the direct sum of the weight spaces $\fg_\gb$ 
for $\gb \in \gS(\faP)\cup \{0\}$ such that $\gb|_{\faPp} > 0.$ 
This collection of $\gb$ splits into (a) those such that $\gb(X) > 0$ 
and (b) those such that $\gb(X) = 0$ and $\gb|_{\faPp} > 0.$ 
The terms satisfying (a) are contained in $\fn_X,$ those satisfying 
(b) satisfy $\gb|_{\ker\ga} =0$ and $\gb \neq 0$ hence $\gb = c \ga$ for $c > 0.$ 
It follows that $\fn_P \subset 
\fn_\ga \oplus \fn_X.$  Conversely, $\fn_X$ is the direct sum of the spaces
$\fg_\gb$ with $\gb(X) > 0.$ The latter condition implies $\gb>0$ on $\faPp$ hence
$\fg_\gb \in \fn_P.$  We see that $\fn_X \subset \fn_P.$ 
Finally, $\fn_\ga$ is the sum of the root spaces $\fg_{c\ga}$ with $c \geq 1.$ 
Since $\ga > 0$ on $\faPp,$ it follows that $\fn_\ga \subset \fn_P.$ This proves
the required identity.
\qed

\begin{lemma}
\label{l: fp in fpX}
$\fp \subset \fp_X.$
\end{lemma}

\proof
$\fp$ is the direct sum of the spaces $\fg_\gb$ were $\gb \in \gS(\fa)\cup\{0\}$ is 
such that  $\gb\geq 0$ on $\fa_P^+.$ For all such $\gb$ one has $\gb(X) \geq 0$ 
so that $\fg_\gb \subset \fp_X.$ 
\qed

From $\fp  \subset \fp_X$ it follows that
$$
P^{(\ga)}  := P \cap M_{1X} =  P \cap G^{(\ga)}
$$ 
is a parabolic subgroup of $G^{(\ga)}$ 
with Langlands decomposition 
$$
P^{(\ga)} =  M_P  A_P N_\ga, 
$$
The centralizer of $G^{(\ga)}$ in $\fa_P$ equals $\ker \ga,$ which has  codimension $1.$
From this we see that $P^{(\ga)}$ is a maximal parabolic subgroup of $G^{(\ga)}.$ 

Let now $Q$ be the adjacent parabolic subgroup determined by 
the  pair $(P, \ga).$ Then $X$ is $(Q, -\ga)$ generic, and it follows
from an easy adaptation of the proof of Lemma \ref{l: deco P in ga X} that
\begin{equation}
\label{e: deco fn Q in ga X}
\fn_Q = \bar \fn_\ga + \fn_X,
\end{equation}
and that $Q^{(\ga)} = 
Q\cap G^{(\ga)}$ is  a parabolic subgroup of $G^{(\ga)},$ with the Langlands decomposition
\begin{equation}
\label{e: Q ga and P ga}
Q^{(\ga)} = M_P A_P \bar N_\ga = \overline{P^{(\ga)}}.
\end{equation}
\vspace{-12pt}

\begin{lemma}
\label{l: properties G ga P ga Q ga}
The group $G^{(\ga)}$ is of the Harish-Chandra class, and 
$P^{(\ga)}$ and $Q^{(\ga)}$ are maximal parabolic subgroups of 
$G^{(\ga)}.$ They have the common split component $A_P$ and are adjacent
and opposite. 
\end{lemma}
\proof
Since the group $G^{(\ga)} = M_{1X}$ is the centralizer of $\fa_X$ in $G,$ 
it belongs to the Harish-Chandra class. From (\ref{e: Q ga and P ga})
it follows that both $P^{(\ga)}$ and $Q^{(\ga)}$ 
have split component $A_P.$ 

The maximality of the parabolic sugroups $P^{(\ga)}$ and $Q^{(\ga)}$ was established
in the above.  These parabolics are opposite. 
In view of their maximality it follows that they are adjacent as well.
\qed

The following results will turn out to be key for the argument reducing the proof of the Maass--Selberg relations
for $B$ to the case of maximal parabolic subgroups. 
It makes that certain data for $G$ and for $G^{(\ga)}$ are suitably compatible.

\begin{lemma}
\label{l: G ga preserves fn X}
The group $G^{(\ga)}$ normalizes $\fn_P\cap \fn_Q.$ 
\end{lemma}

\proof
The  intersection $\fn_P \cap \fn_Q$ is the sum of the root spaces
$\fg_\gb$, ($\gb \in \gS(\fa)$) with $\gb >0$ on both $\faPp$ and $\faQp.$ From the choice of 
$X$ on $\cl \faPp \cap \cl \faQp$ we see that the condition on $\gb$ is equivalent
to $\gb(X) > 0.$ The latter condition in turn is equivalent to $\fg_\gb \subset \fn_X.$ 
It follows that $\fn_P \cap \fn_Q = \fn_X.$ As $G^{(\ga)} = M_{1X},$ the result follows.
\qed 

\begin{lemma}
\label{e: several decos parabs}
We have the following direct sums as linear spaces:
\begin{enumerate} 
\itema
$\fp = \fm_{1P} \oplus \fn_\ga \oplus \fn_X,$
\itemb
$\fq = \fm_{1P} \oplus \bar \fn_\ga \oplus \fn_X.$
\end{enumerate}
If $P$ is standard, then 
\begin{enumerate}
\itemc
$\fn_0 = \fn_0^{(\ga)} \oplus \fn_X.$ 
\end{enumerate}
\end{lemma}

\proof
Since $P$ and $Q$ are adjacent, $\fm_{1P} = \fm_{1Q}.$ Consequently, 
(a) and (b) follow from Lemma \ref{l: deco P in ga X} and (\ref{e: deco fn Q in ga X}). If $P$ is standard, so is 
$P_X,$ in view of Lemma \ref{l: fp in fpX}. It follows that $\fn_0 =
(\fn_0 \cap \fm_{1X}) \oplus \fn_X = \fn_0^{(\ga)} \oplus \fn_X.$ 
\qed

The algebra $\fa$ is maximal abelian in the $-1$ eigenspace of the Cartan involution
$\theta|_{\fg^{(\ga)}}.$ It follows that the normalizer of $A$ in $K^{(\ga)}$ 
maps onto the Weyl group $W^{(\ga)}$ of the root system $\gS^{(\ga)} = \gS(\fg^{(\ga)}, \fa).$ 
It follows from Lemma \ref{l: G ga preserves fn X} that $W^{(\ga)}$ preserves $\fn_X = \fn_P \cap \fn_Q.$

\begin{lemma}
\label{l: compatibility open orbits G and G ga}
Let $v \in N_{K^{(\ga)}}(\fa)$ be such that $N_0^{(\ga)} v \, \bar Q^{(\ga)}$ is open in $G^{(\ga)}\!.$
If $P$ is standard, then $N_0 v \bar Q$ is open in $G.$
\end{lemma}

\proof
The orbit $N_0^{(\ga)} v \bar Q^{(\ga)}$ is open in $G^{(\ga)}$ iff 
$ \fn_0^{(\ga)} + \Ad(v) \bar \fq^{(\ga)} = \fg^{(\ga)}.$
By adding $\fn_X+ \bar \fn_X$ to the left and the right of the latter expression,
we find, using that $\Ad(v)$ normalizes $\bar \fn_X$,
\begin{equation}
\label{e: first orbit deco fg}
(\fn_0^{\ga} + \fn_X)  + \Ad(v)(\bar\fq^{(\ga)} + \bar \fn_X) = \fg.
\end{equation}
From Lemma \ref{e: several decos parabs}  (b), we see that $\fq^{(\ga)} + \fn_X = \fq.$ Using this and
(\ref{e: first orbit deco fg}), we find, by taking Lemma \ref{e: several decos parabs} (c) into account that
$$ 
\fn_0 + \Ad(v)\bar \fq = \fg.
$$ 
This implies that $N_0 v \bar Q$ is open in $G.$ 
\qed

From now on, we will assume that $P$ is standard. 
Since $\Cartan X = - X$, the group $G^{(\ga)} = M_{1X} = Z_G(X)$ 
is invariant under $\theta.$ The restriction $\theta^{(\ga)}: = \theta|_{G^{(\ga)}}$ 
is a Cartan involution. Its group of fixed points is the maximal 
compact subgroup $K^{(\ga)} = K \cap G^{(\ga)}$ of $G^{(\ga)}.$ 

Furthermore,
$\fa$  is a maximal abelian subspace
of $\fs^{(\ga)} = \fs \cap \fg^{(\ga)}.$ 
The algebra  $\fg^{(\ga)} = 
\fm_{1X}$ is the direct sum of the weight spaces $\fg_\gb$ for $\gb \in \gS(\fa)\cup \{0\}$ 
such that $\gb(X) = 0.$ It follows that 
$$
\gS(\fg^{(\ga)}, \fa) = \{ \gb \in \gS(\fa)\mid \gb(X) = 0\}.
$$
Since $P$ is standard, $\cl(\faPp) \subset \cl(\fa^+),$ hence also
$X \in \cl( \fa^+).$ Therefore,
$$
\gS^+(\fg^{(\ga)}, \fa) := \gS(\fg^{(\ga)}, \fa) \cap \gS^+(\fg, \fa)
$$ 
is a positive system for $\gS(\fg^{(\ga)}, \fa).$ It is well-known that the associated set of simple roots is given by
$$
\gD^{(\ga)} := \{ \gb \in \gD \mid \gb(X) = 0\} .
$$ 
The standard minimal parabolic subgroup $P_0$ of $G$ is contained
in $P_X.$ Hence, $P_0^{(\ga)} = P_0 \cap G^{(\ga)}$  is a minimal parabolic
subgroup of $G^{(\ga)} = M_{1X}.$ The nilpotent
radical of $P_0^{(\ga)}$  equals $N_0^{(\ga)} = N_0 \cap G^{(\ga)}.$ Accordingly, the Iwasawa 
decompositions 
$$
G = K A N_0 \quad {\rm and} \quad G^{(\ga)} = K^{(\ga)} A N_0^{(\ga)} 
$$ 
are compatible. The restriction $\chi^{(\ga)} := \chi|_{N_0^{(\ga)}}$ is a unitary 
character of $N_0^{(\ga)}.$

\begin{lemma}
\label{l: char chi ga reg}
The character $\chi^{(\ga)}$ is regular with respect to $G^{(\ga)}, A, N_0^{(\ga)}.$ 
\end{lemma}

\proof 
The $\fa$-roots in $ \fn_0^{(\ga)}$ form the positive system $\gS^+(\fg^{(\ga)}, \fa).$ 
Let $\gb$ be a simple root for this positive system. Then $\gb \in \gD$ and $\gb(X) = 0.$ The simple
root space $\fg^{(\ga)}_\gb$ equals $\fg_\gb.$ The derivative $d\chi^{(\ga)}(e) $ is the restriction
of $d\chi(e)$ to $\fn_0^{(\ga)}.$ Since $\chi$ is regular and $\gb$ simple,
$$
d\chi^{(\ga)}(e)|_{\fg_\gb} = d\chi(e)|_{\fg_\gb} \neq 0.
$$ 
It follows that $\chi^{(\ga)}$ is regular.
\qed

The groups  $P^{(\ga)} = M_P A_P  N_\ga$ and $Q^{(\ga)} = \theta P^{(\ga)}$ 
are opposite maximal parabolic subgroups of $G^{(\ga)}$ with the same split component.
They are adjacent and compatible with $Q$ and $P.$ Let $\gs$ be a representation
of the discrete series of $M_{P}= M_{Q}.$ We consider the characters
$$
\chi_P = \chi|_{M_P\cap N_0} \quad {\rm and}\quad \chi^{(\ga)}_{P^{(\ga)} } = \chi^{(\ga)}|_{M_P \cap N^{(\ga)}_0}.
$$
\newline
\begin{lemma}
\begin{enumerate}
\itema 
$N_0 \cap M_P = N_0^{(\ga)} \cap M_P;$
\itemb
$\chi_P = \chi^{(\ga)}_{P^{(\ga)}};$
\itemc
$\ccH^{-\infty}_{\gs, \chi_P} = \ccH^{-\infty}_{\gs, \chi^{(\ga)}_{P^{(\ga)}}};$
\itemd
the Hermitian inner products on the spaces in (c) are equal.
\end{enumerate}
\end{lemma}

\proof
If $\gb \in \gS^+(\fa)$ is such that $\fg_\gb \subset  \fm_P$ then $\gb|_{\faP} = 0$
so that  $\gb(X)= 0$ which in turn implies that $\fg_\gb \subset \fn_0^{(\ga)}.$ This implies (a).
Since clearly 
 $$
\chi|_{M_P \cap N^{(\ga)}_0} =  \chi^{(\ga)}|_{M_P \cap N^{(\ga)}_0}
$$
it follow from (a) that (b). Assertion (c) is now immediate. For (d) we note that the inner product on 
$\ccH^{-\infty}_{\gs, \chi_P}$ is determined by the requirement that the matrix coefficient map 
$$
\mu_\gs: H_\gs^\infty \otimes \overline{ \ccH^{-\infty}_{\gs, \chi_P}} \to L^2(M_P/M_P\cap N_0: \chi_P)
$$
is an isometry. In view of (a), (b)  and (c), the matrix coeffient map $\mu_\gs$ coincides with 
the matrix coefficient map 
$$
\mu_\gs^{(\ga)}: H_\gs^\infty \otimes \overline{ \ccH^{-\infty}_{\gs, \chi^{(\ga)}_{P^{(\ga)}}}} \to 
L^2(M_P/M_P\cap N_0^{(\ga)}: \chi^{(\ga)}_P).
$$
\begin{rem} 
For (d) it is essential that we agree to equip $M_P/M_P\cap N_0^{(\ga)}$ with the same positive
invariant measure as $M_P/M_P\cap N_0.$
\end{rem}
We have now introduced all ingredients needed for the definition of $B^{(\ga)}(Q^{(\ga)}, P^{(\ga)}, \gs, \nu)$
for the group $G^{(\ga)},$ the adjacent parabolic subgroups $P^{(\ga)}, Q^{(\ga)}$ and any $\gs \in \dMPds$,
as an $\End(
\ccH^{-\infty}_{\gs, \chi^{(\ga)}_{P^{(\ga)}}
}
)$-valued meromorphic function of $\nu \in \faPdc.$ 
The following is a crucial reduction result.

\begin{lemma}
\label{l: comparison B and B ga}
Let $P,Q\in \cP$ be as before, and assume $P$ is standard. Then 
for all $\gs \in \dMPds$ we  have
$$
B(\bar Q,\bar P,\gs, \nu) =  B^{(\ga)}(P^{(\ga)}, \bar P^{(\ga)}, \gs, \nu),
$$ 
as $\End(\ccH^{-\infty}_{\gs, \chi_P})$-valued meromorphic functions of
 $\nu \in \faPdc.$ 
\end{lemma}

\begin{rem}
\label{r: consistent choice of Haar measures}
Here it is important that the Haar measures on $N_Q\cap \bar N_P$ and on $\bar N_\ga$, used in the definition of the standard intertwining operators for $G$  and for $G^{(\ga)},$ are equal.
\end{rem}
  
The proof of Lemma \ref{l: comparison B and B ga} will be given in the next few sections. In the text below we
shall explain its role in completing the proof of  Theorem \ref{t: MS for B}.
                                     
By a Whittaker datum we shall mean a triple $(G , KAN_0, \chi)$  with                                                                                                                                                                                                                                                     $G$ a group of the Harish-Chandra class,                                                                                                                                                                                                                                                                                   $K A N_0$ an Iwasawa decomposition of $G$ and $\chi$ a regular unitary character of $N_0.$ 
An MS setting is a Whittaker datum as above together with a tuple $P,Q$ of parabolic subgroups
of $G$ containing $A,$ and with equal split components. We will say that such a setting satsfies he assertions of Thm. \ref{t: MS for B} if for all $\gs \in \dMPds$ the identity (\ref{e: Maass-Selberg for B}) is valid.
Finally, an MS setting $(G = KAN_0, \chi, P,Q)$ is called basic if $P,Q$ are maximal and adjacent. In particular, in this
case $P$ and $Q$ are opposite.

                                                                                                                                                                                                                                                                \begin{lemma}                                                                                                                                                                                                                                                              
\label{l: basic case MS} If the assertions of Thm. \ref{t: MS for B} hold for every basic $MS$-setting, 
then they hold in general.
\end{lemma}

\proof
Suppose the assertions of  Thm. \ref{t: MS for B} are valid for every basic setting, and let $(G= K A N_0, \chi)$ 
be a Whittaker datum. Let $P',Q'$ determine an associated MS setting. By Corollary \ref{c: reduction to cc} it suffices to prove the assertions of Theorem \ref{t: MS for B} for the setting $(G= K  A N_0, \chi, P',Q')$ under the assumption
that $G$ has compact center.  By Lemma \ref{l: reduction B to adjacent} it suffices to prove
the assertions under the additional condition that $P',Q'$ are 
adjacent. By Lemma \ref{l: MS and v conjugacy} we may further reduce to the case that
 $P'$ is opposite standard. If $P'$ is maximal, so is $Q'$ and hence $(G,P',Q')$ is basic,  and by hypothesis there is nothing 
 left to be proven. Thus, we may in addition assume that $P'$ is not maximal. 
We write $P' = \bar P,$ with $P$ standard. Then $Q := \bar Q'$ is adjacent to $P.$ 
It remains to prove the assertion of Theorem \ref{t: MS for B} for 
the setting $(G = K A N_0, \chi, \bar P, \bar Q).$ 
We now select  a subgroup $G^{(\ga)}$ of $G$  related to the pair $(P,Q)$ as in the previous section. Then 
 $(G^{(\ga)} = K^{(\ga)} A N_0^{(\ga)},  \chi^{(\ga)})$ is a Whittaker datum, and  $(P^{(\ga)}, Q^{(\ga)})$  determines an associated  MS-setting. By Lemma \ref{l: comparison B and B ga} it suffices to prove the assertions of Thm. \ref{t: MS for B}  for the latter setting. By Lemma \ref{l: comparison with cc} we see that it suffices to 
 verify the assertions of Thm. \ref{t: MS for B} for the setting
 $(G^{(\ga)} = K^{(\ga)} {}^{\circ}\!AN_0^{(\ga)},  \chi^{(\ga)}, {}^\circ P^{(\ga)}, {}^\circ Q^{(\ga)}).$ 
Since the latter setting is basic, the validity of the assertions is garanteed by the hypothesis.
\qed

\section{Smoothness of $J$}
\label{s: smoothness J}
In this section, we assume that $P, Q \in \cP$ are adjacent parabolic subgroups of $G,$ containing $A.$ In addition
we assume that $P$ is standard and not maximal. Furthermore, $\ga \in \gS(P)$ and $X \in \ker \ga$ are as in Section \ref{s: parabolics}. We retain the notation introduced in the text following Lemma \ref{l: root deco principle}.

Fix $\eta \in \ccH_{\gs,\chi_{\bar P}}^{-\infty};$ recall that $P$ is standard.
For $\nu \in \faPdc$ we define $\ge_\nu: N_0 \bar P \to H^{-\infty}$ by 
$$ 
\ge_\nu(n ma\bar n) = \chi(n) a^{-\nu +\rho_P} \gs^{-1}(m) \eta,
$$
for $n\in N_P , m\in M_P, a \in A_P$ and $\bar n \in \bar N_P.$
We view the restriction of $\ge_\nu$ to $K\cap N_0 \bar P$ 
as an almost everywhere defined function $K \to H_\gs^{-\infty}.$ 
This function satisfies $\ge_\nu(k m) = \gs(m)^{-1} \ge_\nu(k)$ 
for almost all $k \in K$ and $m \in K_P.$ 

From  \cite[Prop.\ 8.12]{Bunitemp} it follows that for $\Re \nu$ $P$-dominant 
the Whittaker vector $j_\nu := j(\bar P, \gs, \nu)\eta$ is represented by $\ge_\nu$ 
in the sense that in the compact picture one has, for all $\gf \in C^{\infty}(K/K_P:\gs_P),$ 
\begin{equation}
\label{e: repres j by ge}
\inp{j_\nu}{\gf} = \int_K \inp{\ge_\nu(k)}{\gf(k)}\; dk,
\end{equation}
with absolutely convergent integral.

Likewise, the element $j^{(\ga)}_\nu = j (\bar P^{(\ga)}, \gs, \nu)\eta$ associated
with $G^{(\ga)}$  is, for $\nu \in \faPdc$ with $\Re\nu(H_\ga) > 0,$ 
given by the almost everywhere defined function
$$ 
\ge^{(\ga)} (n ma\bar n) = \chi(n) a^{       -\nu +\rho_{ P^{(\ga)} }  } \gs^{-1}(m) \eta,
$$ 
for $n\in N_P^{(\ga)} , m\in M_P , a \in A_P $ and $\bar n \in \bar N_{P^{(\ga)} }.$
Here $\rho_{ P^{(\ga)}}$ is the rho of the standard parabolic subgroup $ P^{(\ga)}$
in $G^{(\ga)}.$ We note that the difference 
$\rho_{P} - \rho_{P^{(\ga)}}$ restricts to zero on $\ker \ga = {}^\circ \fa^{(\ga)}$ hence 
does not appear in the analysis on ${}^\circ G^{(\ga)}.$ 

For $\mu \in \faPdc $ we define $\gf_\mu: G \to \C$ by
$$ 
\gf_\mu(k m a \bar n) = a^{-\mu}, \qquad (k \in K, m\in M_P, a \in A_P, \bar n \in \bar N_P).
$$ 
We note that the operator 
$$
m_{\nu,\mu}: C^\infty(G/\bar P:\gs: \nu) \to C^\infty(G/\bar P: \gs: \nu),\;\;\; \psi \mapsto \gf_\mu\psi
$$ 
is given by the identity of $C^\infty(K/K_P: \gs_P)$ in the compact picture.
In particular, it follows that the operator has a unique continuous linear
extension to an operator $C^{-\infty}(G/\bar P:\gs: \nu) \to 
C^{-\infty}(G/\bar P: \gs: \nu + \mu).$ We denote this operator by 
$m_{\nu, \mu}$ again, and write $m_{\nu,\mu}: \psi \mapsto \gf_\nu\psi.$  
In the compact picture the extended operator $m_{\nu, \mu}$ is given by the identity of $C^{-\infty}(K/K_P:\gs_P).$

It follows from Cor.  \ref{c: holomorphy j Q} 
 that $\nu \mapsto j_\nu$ defines a holomorphic function
$\faPdc \to C^{-\infty}(K/K_P:\gs_P).$ Clearly, $(\nu,\mu)\nu,\mu \mapsto \gf_\mu j_{\nu - \mu}$ is given by 
$(\nu,\mu) \mapsto  j_{\nu- \mu}|_K$ in the compact picture, hence is holomorphic on $\faPdc \times \faPdc$.

It follows from  \cite[Prop. 8.14-15]{Bunitemp}  that for $\Re  \nu$ strictly $P$-dominant 
the generalized function  $\gf_\mu j_{\nu -\mu} \in C^{-\infty}(G/\bar P:\gs:\nu)$ 
is represented by $\gf_\mu \ge_{\nu - \mu}$ in the sense that for all $\gf\in C^\infty(K/K_P:\gs_P),$ 
$$ 
\inp{\gf_\mu  j_{\nu -\mu}}{ \gf} = \int_K \inp{\ge_{\nu - \mu}(k)}{ \gf(k)} \;dk.
$$ 

Put $A(\nu) = A(\bar Q , \bar P, \gs, \nu)$ and let $\cS_A\subset \faPdc$ denote the singular locus for 
$A(\dotvar).$ For $(\nu,\mu) \in (\faPdc \setminus \cS_A)\times \faPdc  , $ we define 
\begin{equation}
\label{e: first intertwiner on phi j}
J_{\nu,\mu} :=A( \nu) [\gf_{\mu} j_{ \nu - \mu}] \in C^{-\infty}(G/\bar Q:\gs :\nu).
\end{equation}
The extra parameter $\mu$ is introduced to allow the choice of 
pairs  $(\nu, \mu)$ where $A(\nu)$ and $j_{\nu -\mu}$ are simultaneously representable
by a convergent integral and a locally integrable function, respectively. This idea is inspired 
by an  argument of T. Oshima and J. Sekiguchi 
\cite[text prec. Lemma 4.13]{OS80},  developed further by  \cite[Lemma 7.4]{Bps1}, 
see also \cite[Prop. 6]{CD}. 

Our first main goal is to understand 
the dependence of $\pi_{\bar Q, \gs, \nu}(n) J_{\nu,\mu}$ 
on $(n, \nu, \mu) \in N_P \times \faPdc \times\faPdc.$  
 
From (\ref{e: defi fa Q R}) and (\ref{e: defi fa Q21}) we have, for $R \in \R,$
\begin{eqnarray*}
\fad(P, R): &=& \{\nu\in\fa_{P\iC}^*\mid
\;\Re\inp{\gb}{\nu} > R,\;\; (\gb \in \gS(P))\},\\
\fad(Q|P, R): &=&  \{\nu\in\fa_{P\iC}^*\mid
\Re\inp{\gb}{\nu} > R,\;\; (\gb \in \gS(\bar Q) \cap \gS(P))\}.
\end{eqnarray*}

From  \cite[Prop 14.8]{Bunitemp} we know that for every $R\in \R$ there exists a positive integer $s$
such that $\nu \mapsto j(\bar P ,\gs, \nu)$ is holomorphic as a map 
$\fa_{P}^*(P, R) \to C^{-s}(K/K_P:\gs_P).$  

Furthermore,  for $\gs \in \widehat M_{P, { \rm ds}}.$
we know by Lemma \ref{l: mero prop A} that for every $R \in \R$ there exists a polynomial function $q: \faPdc \to \C$ 
and a constant $r \in \N$ such that for every $t \in \N$ the assignment
$$
\nu \mapsto q(\nu)\,A(\bar Q,\bar P, \gs, \nu)
$$
is holomorphic as a function on $\fad(\bar Q|\bar P, R)$ with values in $B(C^{-t}, C^{-t- r})$, the Banach space 
 of bounded linear maps 
$C^{-t}(K/K_{P}:\gs_{P}) \to C^{-t -r}(K/K_{P}:\gs_{P})$.

\begin{lemma}
\label{l: first holomorphy J}
For every bounded open subset $\Omega$ of $\faPdc \times \faPdc$  there exists a $p \in \N$ and a polynomial function $q\in P(\faPd)$ 
such that the map $(n, \nu, \mu)\mapsto   q(\nu) \pi_{\bar Q, \gs, \nu}(n) J_{\nu,\mu}$ 
is a smooth map $N_P \times \Omega \to C^{-p}(K/K_P:\gs_P)$ 
which is holomorphic in the variable $(\nu, \mu) \in \Omega.$
\end{lemma}

\proof
Without loss of generality we may assume that $\Omega = \Omega_1 \times \Omega_2$ with
$\Omega_j$ bounded open in $\faPdc.$ 
We will write $C^s(K:\gs_P):= C^s(K/K_P : \gs_P),$ and keep in mind
that $K_P = K_Q.$
For every $s\in \Z$ we  will write $\psi \mapsto \psi|_K$ for the isomorphism 
$C^s(G/\bar P:\gs:\nu) \to C^s(K:\gs_P)$ induced by restriction to $K.$ 
Note that with this notation,
$$
 \pi_{\bar P, \gs, \nu} (n) [ \gf_{\mu} j(\bar P , \gs, \nu - \mu, \eta)]|_K 
 = \psi(n, \mu) [ j(\bar P , \gs, \nu - \mu, \eta)|_K]
 $$ 
 where 
 $$ 
 \psi(n, \mu)(k)  = \chi(n)^{-1} \gf_{-\mu}(n^{-1}k) = \chi(n)^{-1} a_{\bar P}(n^{-1} k)^{-\mu}
 $$
 is a smooth function $N_P \times \faPdc \to C^\infty(K/K_M),$ which is holomorphic 
 in the second variable.
 Let $r$ be the order of $j(P, \gs, \dotvar)$ over $\Omega_1 - \Omega_2 = \{\nu - \mu\mid \nu \in \Omega_1, \mu\in \Omega_2\}.$ Then  $(\nu,\mu) \mapsto j_{\nu - \mu} $ defines a holomorphic map $\Omega \to C^{-r}(K: \gs_P)$.
 It follows that  the map 
$$ 
(n, \nu, \mu) \mapsto \pi_{\bar P, \gs, \nu} (n) [\gf_{\mu} j_{\nu - \mu}|_K]
$$ 
is smooth $N_P \times \Omega  \to C^{-r}(K:\gs_P)$ 
and in addition holomorphic in the second variable.

Let $t$ be the order of the family $A(\bar Q,\bar P, \gs, \dotvar)$ over $\Omega_1. $
Then there exists a polynomial function $q: \faPdc \to \C$ 
such that for every positive integer $s$ the operator $A(\bar Q, \bar P, \gs, \dotvar)$ defines a holomorphic 
function $\Omega_1 \to  B(C^{-s}(K:\gs_P), C^{-s - t}(K:\gs_P)).$ 
 Since the natural map 
$$
B(C^{-r}(K:\gs_P), C^{-r - t}(K:\gs_P))\;\times \;C^{-r}(K:\gs_P) \to C^{-r-t}(K:\gs_P)
$$
 is a continuous bilinear map of Banach spaces, it follows from the usual 
rules for differentiation that 
%\end{document}
 \begin{eqnarray}
 \lefteqn{\!\!\!\!\!\!\!\!\!\!(n, \nu) \mapsto 
 q(\nu) \pi_{\bar Q, \gs, \nu}(n) J_{\nu,\mu}} \nonumber
 \\
 &=& 
 q(\nu) A(\nu)  
 \pi_{\bar P, \gs, \nu} (n) [\gf_{\mu} j_{\nu - \mu}|_K]
 \label{e: formula for J nu mu}
 \end{eqnarray}
 defines a smooth map from $N_P \times \Omega$ to 
 $C^{-r -t }(K: \gs_P)$ which is holomorphic in the second variable. 
 This proves the result with $p = r + t.$  
 \qed
 %\end{document}
 
 In view of the results of Section \ref{s: defi j} for the group $G^{(\ga)}$
 in place of $G$ 
  there exists an element $v \in N_{K^{(\ga)}}(\fa)$ such that $N_0^{(\ga)} v \bar Q^{(\ga)}$ is open in $G^{(\ga)}.$ 
 We may therefore choose 
 $$
 v_{\bar Q} = v^{(\ga)}_{\bar Q^{(\ga)}} := v
 $$
 for our maps $v:\cP \to N_K(\fa)$ and $v^{(\ga)}: \cP^{(\ga)} \to N_{K^{(\ga)}}(\fa)$
 for $G$ and $G^{(\ga)}$ as discussed in Definition \ref{d: choice of vQ}. 
 Then $\bar P = v\bar Qv^{-1}.$  
 It follows that  $G^\od := N_P v \bar Q$ is a
 right $\bar Q$-invariant (dense) open subset of
 $G.$ The action map $N_P  \to G/\bar Q$, $n \mapsto n v \bar Q,$ induces  
 an open embedding into $G/\bar Q$ with image $G^\circ/ \bar Q.$ 
 Composing the defined  embedding $N_P \to G/\bar Q$ with the inverse of the diffeomorphism  $K/K_Q \to G/\bar Q $ we obtain an embedding $N_P \to K/K_Q$ 
 with image $K^\od = [N_P v \bar Q] \cap K.$ The defined maps form a commutative diagram of diffeomorphisms
 $$
 \begin{array}{ccc}
 N_P & \longrightarrow & G^\od/\bar Q\\
 & \searrow& \uparrow\\
 && K^\od/K_Q.
 \end{array}
 $$ 
  By pull-back we then obtain for every $\nu \in \faPdc$ 
 a commutative diagram of topological linear isomorphisms
 $$
 \begin{array}{ccc}
 C^{r}(N_P, H_\gs) & \longleftarrow & C^{r}(G^\od /\bar Q: \gs : \nu)\\
 &\nwarrow& \uparrow\\
 && C^{r}(K^\od / K_Q: \gs_Q).
 \end{array}
 $$
 The diagram with
 arrows representing the inverted maps is still a commutative diagram of topological linear isomorphisms.
 For a given $u \in H_\gs^\infty$ we consider the embedding 
 $\iota_u: C^p(N_P) \to C^{p}(N_P, H_\gs),$ $f \mapsto f \otimes u.$ 
 Combining this with the inverted diagram, we obtain a diagram
 of continuous linear maps
 \begin{equation}
 \label{e: embedding of Cr NP} 
 \begin{array}{ccc}
 C^{r}(N_P) &{\buildrel {}_u^{ \bp} T_\nu \over \longrightarrow} & C^{r}(G^\od/\bar Q: \gs : \nu)\\
 &\!\!\!\!\!\!\!\!\!\! \scriptstyle{ {}_u T_\nu} \searrow& \downarrow r_\nu \\
 && C^{r}(K^\od/ K_Q: \gs_Q).
 \end{array}
 \end{equation}
 
 \begin{lemma}
For every $u \in H_\gs^\infty$ and $\nu \in \faQdc,$  
$$ 
\pi_{\bar Q, \gs, \nu}(n) \after {}_uT_{\nu} = {}_uT_{\nu} \after L_n, \qquad (n \in N_P).$$
\end{lemma}

\proof
This follows readily from the definitions.
\qed

 Given a compact subset $\cK \subset N_P$ we will denote the canonical image of 
$ \cK v$  in $G/\bar Q$ by $\cK'$ and the  image in $K/K_Q$ by $\cK''.$ 
 Note that $\cK'\subset G^\od/\bar Q$ and $\cK''\subset K^\od/K_Q.$
 Let $r$ be a  positive integer.
 Identifying $C^r_{\cK'}(G^\od/ \bar Q:\gs:\nu)$ with $C^r_{\cK'}(G/ \bar Q)$ 
 in the usual way, through extension by zero, and using the analogous identification
 for functions on $K,$ we infer that the diagram (\ref{e: embedding of Cr NP}) induces 
 a commutative diagram 
 \begin{equation}
 \label{e: diagram with specified compact supp}
 \begin{array}{ccc}
 C_\cK^{r}(N_P) &{\buildrel {}_u^{ \bp} T_\nu \over \longrightarrow} & C^{r}_{\cK'}(G/\bar Q: \gs : \nu)\\
 &\!\!\!\!\!\!\!\!\!\! \scriptstyle{ {}_u T_\nu} \searrow& \downarrow r_\nu \\
 && C^{r}_{\cK''}(K/ K_Q: \gs_Q)
 \end{array}
 \end{equation}
of bounded linear maps between Banach spaces. 
If $f \in C^r(N_P)$ we denote by ${}_u f_{\bar Q, \gs,  \nu}$ the function $G^\circ \to H_\gs^\infty$
given by 
\begin{equation}
\label{e: defi f bar Q nu}
{}_u f_{\bar Q, \gs,  \nu}(nv ma\bar n) = a^{-\nu + \rho_Q}f(n) \otimes \gs(m)^{-1}u,
\end{equation}
for $n \in N_P, m \in M_Q, a\in A_Q$ and $\bar n \in \bar N_Q.$ 
 Then 
$ {}_u^{ \bp} T_\nu $ and ${}_u T_\nu$ are given by 
$$
 {}_u^{ \bp} T_\nu (f) = {}_u f_{\bar Q, \gs, \nu}, \;\;\;{\rm and} \;\;\;{}_u T_\nu(f) =
  {}_u^{ \bp} T_\nu (f)|_{K^\circ},  \qquad (f \in C^r(N_P)).
$$ 
If $f \in C_\cK^r(N_P)$ we  view ${}_u^{ \bp} T_\nu (f) $ as an element
of $C^r_{\cK'}(G/\bar Q:\gs:\nu)$ as explained above. Then $  {}_u T_\nu (f) = {}_u^{ \backprime} T_\nu (f)|_{K}. $

\begin{lemma}
Let $\cK \subset N_P$ be compact, and $ \nu \in \faPdc.$  For every $r \in \N$ the map $(u, f) \mapsto 
{}_uT_\nu(f)$ is continuous bilinear $H_\gs^\infty \times C_\cK^{r}(N_P)\to C^r_{\cK''}(K/K_Q,\gs_Q).$
\end{lemma} 

\proof Straightforward. See also  \cite[Lemma 8.8]{Bunitemp} for a related
discussion.
\qed

We will now investigate the family $J_{\nu, \mu} \in C^{-r}(G/\bar Q: \gs :\nu)$ in more detail.  

First of all, if $J  \in C^{-r}(G/\bar Q: \gs : \nu)$ and $u \in \ccH_\gs^\infty$ we define the continuous 
linear functional
${}_u J : C^r_c(N_P) \to \C$ by 
\begin{equation}
\label{e: defi u J}
{}_u J (f) = \inp{J }{ { }_u T_{-\bar \nu}(\bar f)},\qquad (f \in C^r(N_P)).
\end{equation}
 
We denote by $a_{\bar Q}$ the function $K^\od \to A_Q$ (uniquely)  determined by 
$$
x \in N_P v M_Q a_{\bar Q}(x) \bar N_Q, \qquad (x \in K^\od).
$$ 
Then $a_{\bar Q}$ is real analytic $\Kcirc/K_Q \to A_Q.$ For $r \in \N$ and $\nu \in \faPdc$ we define the map 
$\rmm_\nu: C^r(K^\od /K_Q:\gs_Q) \to C^r(K^\od /K_Q:\gs_Q)$ by 
$$
\rmm_\nu(f)(k) = a_{\bar Q}(k)^{-\nu} f(k), \qquad (f \in C^r_c(\Kcirc /K_Q)).
$$

\begin{lemma}
\label{l: property map m nu}
Let $\cK$ be a compact subset of $N_P.$ Then for every $\nu\in \faQdc$  
the map ${\rm m}_\nu$ restricts to 
a bounded automorphism of the Banach space 
$
C^r_{\cK''}(K/K_Q:\gs_Q). 
$
The assignment $\nu \mapsto \rmm_\nu$ is holomorphic as a function on $\faPdc$ with values
in the space of bounded operators $B(C^r_{\cK''}(K/K_Q:\gs_Q)).$
\end{lemma}
\proof 
Straightforward.
\qed

\begin{lemma}
\label{l: key lemma for T}
Let $u \in H_\gs^\infty$. Let $\cK \subset N_P$ be a compact subset. 
For $f \in C_\cK^r(N_P)$ and 
$\nu \in \faQdc$ we have 
\begin{enumerate}
\itema 
${}_u T_\nu( f) \in C^r_{\cK'}(K/K_Q : \gs:  \nu);$
\itemb
${}_u T_\nu(f) =  \rmm_{\nu}  ( {}_u T_0(f) );$
\itemc
the map ${}_u T_\nu$  is a bounded linear map between the Banach spaces
$C^r_\cK(N_P)$ and $C^r_{\cK''}(K/K_Q: \gs_Q);$ 
\itemd 
the assignment  $\nu \, \mapsto {}_uT_\nu$ is holomorphic as a function on $\faQdc,$
with values in the Banach space $B$ of bounded linear maps from $C_{\cK}^r(N_P)$ to  $C^r_{\cK''}(K/K_Q: \gs_Q),$ equipped with the operator norm.
\end{enumerate}
\end{lemma} 

\proof Assertion (a) is true by definition. 
For (b), fix $k \in \cK'.$ Then $k$ has a unique decomposition
$k= n v m a \bar n \in \cK v \bar Q \simeq \cK v M_Q A_Q \bar N_Q
\subset N_P v \bar Q.$ 
It follows that 
\begin{eqnarray*}
{}_u T_\nu (f)(k)&= &  f(n) a^{-\nu + \rho_Q} \gs(m)^{-1}u\\
&=& 
a_{\bar Q}(k)^{-\nu} {}_u T_0(f)(k) = m_\mu({}_u T_0(f))(k).
\end{eqnarray*}

From (a) and (b) it follows that ${}_u T_\nu = \rmm_\nu \after {}_u T_0.$
Thus, view of Lemma \ref{l: property map m nu} it suffices to show (c) for $\nu = 0.$

Put $X = C^p_\cK(N_P),$  $Y = C^p_{\cK''}(K/K_Q: \gs_Q);$ these are
Banach spaces. We equip $B(X,Y)$ and $B(Y)$ with the operator norms
and consider the natural map $\gb: B(X, Y) \times B(Y) \to B(X,Y)$ 
$(\tau, \mu) \mapsto \mu \after \tau.$ Then $\gb$ is bilinear and 
$\|\gb(\tau, \mu)\|_{\rm op} = \|\mu \after \tau\|_{\rm op} \leq \|\mu\|\|T\|.$ Thus, $\gb$ is continuous. If $\nu \mapsto \tau_\nu$ and $\mu \mapsto \mu_\nu$ are holomorphic,
it follows readily that $\mu \mapsto \gb(Y_\nu, \mu_\nu)$ is holomorphic 
with values in $B(X,Y),$ Applying this to $\tau_\nu = {}_uT_0$   and $\mu_\nu = m_\nu,$ 
we find that $\nu \mapsto m_\mu \after {}_u T_\nu = \gb({}_u T_0, m_\nu)$ 
is holomorphic  $\faQdc \to B(X,Y),$ This  establishes (d).
\qed

Recall the definition of $J_{\nu,\mu}$ from (\ref{e: first intertwiner on phi j}), for $(\nu, \mu) \in (\faPdc\setminus \cS_A) \times \faPdc.$ 
For $u \in H_\gs^\infty$ 
we define ${}_u J_{\nu, \mu}$ as in (\ref{e: defi u J}). This is an element of $C^\infty_\cK(N_P) '$
for every compact $\cK \subset N_P.$ 

\begin{thm}
\label{t: tilde J smooth}
For every $u\in H_\gs^\infty$ and $(\nu,\mu) \in (\faPdc \setminus \cS_A )\times \faPdc$ 
there exists a unique smooth function ${}_u \widetilde J_{\nu, \mu}\in C^\infty(N_P)$ such that
${}_u J_{\nu,\mu}$ is represented by the density ${}_u \widetilde J_{\nu, \mu} dn_P$ in the sense  
that 
$$
{}_u J_{\nu,\mu}(f) = \int_{N_P} {}_u \widetilde J_{\nu, \mu} (n_P) f(n_P) d n_P, \qquad (f \in C^\infty_c(N_P)).
$$
If $\Omega \subset \faPdc \times \faPdc$ is a bounded open subset then there exists a polynomial function $q: \faPdc \to \C$ such that the map $(\nu, \mu)\mapsto q(\nu) {}_u \widetilde J_{\nu, \mu}$ is holomorphic $\Omega \to C^\infty(N_P).$
\end{thm}

\proof
Let $u\in H_\gs^\infty$ be fixed. The asserted uniqueness is clear. Therefore,
it suffices to prove the assertions that arise if we replace 
$(\faPdc \setminus \cS_A )\times \faPdc$ by $(\Omega_1\setminus \cS_A)\times \Omega_2$ 
and $\Omega$ by $\Omega_1 \times \Omega_2$ for an arbitrary pair  $\Omega_1, \Omega_2$  of bounded open subsets of $\faPdc.$ Suppose such a pair is fixed and let $q \in P(\faPdc)$ and $p \in \N_+$ be associated with 
$\Omega = \Omega_1 \times \Omega_2$ as in Lemma 
\ref{l: first holomorphy J}.

Let $\cK \subset N_P$ be a compact subset.
We fix an open neighborhood $V$ of $ e$ in $N_P$ whose closure in $N_P$
is compact. 
Then $\cK_e := \cl (V^{-1}) \cK$ is a compact subset of $N_P.$
We fix $\cK_e'$ and $\cK_e''$ as in the discussion
of the diagram (\ref{e: diagram with specified compact supp}) with $\cK_e$ in place of $\cK.$ 

For $f\in C^p_\cK(N_P)$ and $n \in V$ we have $L_n^{-1} f \in C_{\cK_e}^p(N_P).$  
We now note that  
\begin{eqnarray*}
 \inp{L_n[ {}_u J_{\nu, \mu}]}{\bar f} &=& \inp{{}_u J_{\nu, \mu}}{L_n^{-1} \bar f}
 \\
 &=&\inp{J_{\nu, \mu}}{{}_u T_{\mu, -\bar\nu} [L_n^{-1} \bar f]}
 \\
 &=& \inp{J_{\nu, \mu}}{\pi_{Q, \gs, -\bar\nu}(n)^{-1}[{}_u T_{\mu, -\bar\nu}\bar f]}
 \\
 &=& \inp{\pi_{Q, \gs, \nu}(n) J_{\nu, \mu}}{ {}_u T_{-\bar \nu} (\bar f)}.
\end{eqnarray*}

We write $X= C^p_\cK(N_P)$ and  $Y= C^p(K/K_Q:\gs_Q).$ Furthermore,
$\bar Y$ denotes the conjugate of $Y$ and $B(X, \bar Y)$ the space of bounded linear maps from $X$ to $Y,$ equipped with the operator norm.
  
It follows from Lemma \ref{l: key lemma for T} (d) that
$f \mapsto {}_u T_{-\bar \nu}(\bar f)$ is an element of 
$B(X, \bar Y),$ depending holomorphically on $\nu \in \faPdc.$ On the other hand, $(n, \nu, \mu) \mapsto q(\nu) \pi_{Q, \gs, \nu}(n) J_{\nu, \mu}$
is a smooth function $V \times \Omega \to \bar Y' = C^{-p}(K/K_Q:\gs_Q),$
which is holomorphic in $(\nu,\mu).$ 

We now consider the natural bilinear map $\gb: B(X, \bar Y) \times \bar Y' \to X' $ 
given by $(t, \eta) \mapsto \eta \after t.$  Note that $\|\gb(t,\eta)\|\leq \|\eta\|\|t\|;$ 
this shows that $\gb$ is continuous bilinear. 
We observe that 
$$
L_n[{}_u J_{\nu,\mu}]  = \gb({}_u T_{-\bar \nu}\;,\; \pi_{Q,\gs,\nu}(n)J_{\nu,\mu}).
$$ 
By the usual rules for differentiation it follows that $(n, \mu,\nu) \mapsto q(\nu) L_n[{}_u J_{\nu,\mu}]$
is a smooth map  $N_P \times \Omega_{R_1, R_2} \to C^p_\cK(N_P)'$, which 
is holomorphic in the variable from $\Omega.$
 
Let $\bp v_1, \ldots , \bp v_n$ be a basis of the Lie algebra of $N_P$ 
and let $v_j= L_{\bp v_j} $ be the associated right invariant vector fields
on $N_P.$ Then it follows for all $(\nu,\mu) \in \Omega$
and every multi-index $\ga \in \N^n$ that the 
distribution $q(\nu) v^\ga({}_u J_{\nu,\mu})$ belongs to $C^p_\cK(N_P)'.$ 
Furthermore, 
$$
(\nu, \mu) \mapsto q(\nu) v^\ga({}_u J_{\nu,\mu}): \Omega_{R_1, R_2} \to C_\cK^p(N_P)'
$$ 
is holomorphic. By application of the lemma of the appendix we now 
conclude that $q(\nu) {}_u J_{\nu, \mu}$ is a smooth density of the form asserted,
with holomorphic dependence on $(\nu,\mu).$ 
\qed

\section{A useful integral formula}
\label{s: integral formula}
We keep working under the hypothesis of Section \ref{s: parabolics}.
Thus, $P\in \cP$ is standard non-maximal, $Q$ is adjacent to $P$ 
and $\ga \in \gS(P) \cap [- \gS(Q)].$ The element $X \in \ker \ga\cap \cl(\faPp)$ is generic.
$P_X$ is the unique parabolic subgroup having $X$ in its positive chamber,
and $G^{(\ga)} = M_{1X} = P_X \cap \bar P_X.$

We define smooth maps $k_{\bar P}: G \to K$, $m_{\bar P}: G \to M_P\cap \exp \fs,$ $a_{\bar P}: G \to A_P$ and $n_{\bar P}: G \to \bar N_P$ by 
$$
x = k_{\bar P}(x) m_{\bar P}(x) a_{\bar P}(x) n_{\bar P}(x), \qquad (x \in G).
$$
The multiplication map $K\times \bar P \to G$ 
factors through  a diffeomorphism $K \times_{K_P} \bar P \to G.$ 
Since $a^{- 2\rho_P}dm da d\bar n$ defines a right-invariant measure on the
group $\bar P,$ it follows that for a Lebesgue integrable function 
$f: G \to \C$ we have
$$ 
\int_G f(x) dx = \int_{K \times M_P \times A_P \times \bar N_P} f(kma\bar n)\;
dk dm da d\bar n.
$$ 

\begin{lemma}
Let $\gf: K \to \C$ be Lebesgue integrable.
Then 
$$
\int_K \gf(k)\; dk = \int_{N_X \times K^{(\ga)}}\gf(k_{\bar P}(n_X k_\ga))\; a_{\bar P}(n_X k_\ga)^{2\rho_P}\;
dn_X dk_\ga.
$$ 
\end{lemma}

\proof
We fix $\psi\in C_c(\bar P)$ left $K_P$-invariant, such that 
\begin{equation}
\label{e: int psi with rho is 1} 
\int_{M_P \times A_P\times  \bar N_P} \psi(ma \bar n)  a^{- 2\rho_P} dm da d\bar n = 1.
\end{equation}
Furthermore, we extend $\gf$ to $G$ by the formula 
$$ 
\gf(k m a \bar n) = \gf(k) \psi(m a n),
$$ 
for $k \in K, m\in M_P, a\in A_P$ and $\bar n \in \bar N_P.$ 
Then $\gf$ is Lebesgue integrable on $G$ and 
\begin{eqnarray*}
\int_G \gf(x)\; dx  &=& \int_{K\times M_P \times A \times \bar N_P } \gf(k)\psi(ma\bar n) a^{-2\rho_P} dk dmdad\bar n
 \\
&=& \int_K \gf(k)\; dk.
\end{eqnarray*} 
Since $K^{(\ga)} $ normalizes $N_X,$ the density $dn_X d k_\ga$ on $N_X \times K^{(\ga)}$ 
is left invariant.  Therefore, 
\begin{equation}
\int_G \gf(x)\; dx = \int_{N_X \times K^{(\ga)} \times M_P \times A_P \times \bar N_P} \gf(n_X k_\ga m a \bar n)\;
a^{-2\rho_P} \; dn_X dk^{(\ga)}  dm da d\bar n,
\label{e: int gf over G}
\end{equation}
provided the measures are suitably normalized. The pull-back of $a^{-2\rho_P} dm da dn$ on $\bar P$
under left multiplication by $m_1 a_1 \bar n_1$ equals $a^{2\rho_P} _1  a^{-2\rho_P} dm da dn.$
Therefore, the second integral in (\ref{e: int gf over G}) equals
\begin{eqnarray*}
\!\!\!\!\!\!\!\!\!\!\!\!\!\!\!\!\!\!\!\!\!\!\!\!\!\!\!\!\!\!\!\!\!\!\!\!\!\!\!\!\!\!\!\!\!\!\!\!\!\!\!\!\!\!\!\!\!\!\!\!\!\!\!\!\!\!\!\!\!\!\!\!\!\!\!\!\!\!\!\!\!\!\!\!\!\!\!\!\!\!
\lefteqn{\!\!\!\!\!\!\!\!\!\!\!\!
\int_{N_X \times K^{(\ga)} \times M_P \times A_P \times \bar N_P} \gf(k_{\bar P}(n_X k_\ga))
a_{\bar P}(n_X k_\ga)^{2\rho_P} \; \psi( m a \bar n)a_{\bar P}^{-2\rho_P}   dm da d\bar n dn_X dk_\ga\;\;\;\;\;\;\;\;}
\\
\!\!\!\!\!\!\!\!\!\!\!\!\!\!\!\!\!\!\!\!\!\!\!\!\!\!\!\!\!\!\!\!\!\!\!\!\!\!\!\!\!\!\!\!
&=& \int_{N_X \times K^{(\ga)} } \gf(k_{\bar P}(n_X k_\ga))\; a(n_X k_\ga)^{2\rho_P} dn_X dk_\ga. 
\end{eqnarray*}
\qed

\begin{lemma}
\label{l: integral formula ga X}
Let $f: G \to \C$ be such that $f|_K$ is Lebesgue integrable and 
$$
f(xma\bar n) = a^{2\rho_P}  f(x),
$$
for all $x \in G, (m ,a, \bar n) \in M_P \times A_P \times \bar N_P.$ Then
$$
\int_K f(k) dk = \int_{N_X \times K^{(\ga)}}  f(n_X k_\ga)\; dn_X dk_\ga.
$$
\end{lemma}

\proof
\begin{eqnarray*}
\int_K f(k) dk & = &  \int_{N_X \times K^{(\ga)}}  f(k_{\bar P}(n_X k_\ga))  \; a_{\bar P}(n_X k_\ga)^{2\rho_P}\; dn_X dk_\ga.
\\
&  =  & \int_{N_X \times K^{(\ga)}}  f(n_X k_\ga)\; dn_X dk_\ga.
\end{eqnarray*}
\vspace{-24pt}

\qed{\ }
\begin{rem} 
\label{r: integral formula}
In proof given above we have not used the particular definitions of $N_X$ and $K^{(\ga)}.$ 
The proof works under the assumptions that $\fn_X , \fn_\ga$ are sums of $\fa$-root spaces such that $\fn_X \oplus \fn_\ga = \fn_P,$ $N_X = \exp \fn_X,$ 
$\fk^{(\ga)} = \fk\cap (\fn_\ga + \bar\fn_\ga) + \fk_P$ and $K^{(\ga)}$ is the
group generated by $\exp (\fk^{(\ga)}) K_P.$  In particular the proof works for the case 
$\fn_\ga = 0,$ so that in particular $N_X = N_P$ and $K^{(\ga)} = K_P.$ In this setting 
the above result is well known.
\end{rem}

In the sequel we will also need the following result, for $P\in \cPst$, and $Q= v^{-1}P v,$ $v \in N_K(\fa).$

\begin{lemma}
\label{l: int over K is int over NPv}
Let $f: G\to C$ be right $M_Q\bar N_Q$-invariant, and let $R_a f = a^{2\rho_Q}f$ 
for all $a \in A_Q.$ If $f|_K$ is Lebesgue integrable, then 
$$
\int_K f(k)\; dk = \int_{N_P} f(nv) \; dn.
$$
\end{lemma}

\proof The function $L_v f: G\to \C$ has the same $\bar Q$-equivariance on the right at $f.$
In view of Remark \ref{r: integral formula} we obtain
$$
\int_K  L_{v} f(k)\; dk = \int_{N_Q} f(vn_Qv^{-1} v) dn_Q = \int_{N_P}  f(nv) dv.
$$
Since $dk$ is left invariant, the desired result follows.
\qed
%%%%%%%%%%%%%%%%%%%%%%%%%%%%%%%%%%%%%%

\section{Comparison of $J$ with $J^{(\ga)}$}
\label{s: comparison J}
The preceding discussion applies to any Whittaker datum $(G = K A N_0, \chi).$  In particular it applies to the group $G^{(\ga)} = K^{(\ga)} A N_0^{(\ga)},$ with the 
character $\chi^{(\ga)} = \chi|_{N_0^{(\ga)}};$ see Lemma \ref{l: char chi ga reg}
 and its adjacent parabolic
subgroups $P^{(\ga)} = P \cap G^{(\ga)}$ and $Q^{(\ga)} = Q \cap  G^{(\ga)}.$
Their respective nilpotent radicals are $ N_\ga$ and $\bar  N_\ga.$ Accordingly,
$P^{(\ga)}$ and $Q^{(\ga)}$ are opposite parabolic subgroups of $G^{(\ga)}$
with split components $A_P$. Since they are maximal parabolic subgroups of $G^{(\ga)}$ they are adjacent.
The element $v$ belongs to the normalizer in $K^{(\ga)}$ of  $\fa.$ In particular, from $P = vQv^{-1}$
it follows that $P^{(\ga)} = v Q^{(\ga)} v^{-1}.$ Thus, $N_0 v Q$ is open in $G$ and 
$N_0^{(\ga)} v Q^{(\ga)}$ is open in $G^{(\ga)}.$ Moreover, $N_0 v \bar Q = N_P v \bar Q$
and $N_0^{(\ga)} v \bar Q^{(\ga)} = N_\ga v \bar Q^{(\ga)}.$ 

We recall that $\fn_P = \fn_\ga \oplus \fn_X.$ Since both $\fn_\ga$ and $\fn_X$ are subalgebras of $\fn_P$
and each of them is a direct sum of root spaces $\fg_\gb$ with $\gb \in \gS(\fa),$ it follows that 
the multiplication map $N_X \times N_\ga \to N_P$ is a diffeomorphism. Since $G^{(\ga)}$  normalizes $N_X,$
it follows that $N_X$ is a normal subgroup of $N_P$ and 
$$
N_P = N_X N_\ga = N_X  \rtimes N_\ga \qquad \mbox{\rm  (semidirect product)}.
$$
We note that also $N_X K^{(\ga)}$ is a closed subgroup of $G$; clearly 
$$
N_X K^{(\ga)} = N_X \rtimes K^{(\ga)}.
$$

For $\nu \in \faPdc$ with $\Re \nu$ strictly $P^{(\ga)}$-dominant the Whittaker vector
$j^{(\ga)}(\bar P^{(\ga)}, \gs, \nu, \eta)$ is represented by the function
$\ge^{(\ga)}: G^{(\ga)} \to \C$ defined by 
$$ 
\ge^{(\ga)}_\nu(n_\ga m a \bar n_\ga) = \chi(n_\ga)\gs(m)^{-1} a^{-\nu + \rho_{P^{(\ga)}}}\eta
$$ 
for $(n_\ga, m, a, \bar n_\ga) \in N_\ga \times M_P \times A_P\times \bar N_\ga$ and by 
zero on the complement of $N_\ga M_P A_P \bar N_\ga.$ 
For $\mu \in \faPdc$ the function $\gf^{(\ga)}_\mu: G^{(\ga)} \to \C$ is defined  by 
$$ 
\gf^{(\ga)}_\mu(k_\ga m a n_\ga) = a^\mu,
$$  
for $(k_{\ga}, m,  a , n_\ga) \in K^{(\ga)}\times M_P\times A_P \times N_P.$ 
We define the character $\xi: A_P \to \R$ by 
$$
\xi(a) =a^{\rho_P - \rho_{P^{(\ga)}}} 
$$
and note that $\xi = 1$ on $\exp \R H_\ga.$ Now $\exp \R H_\ga \subset M_X=
{}^\circ G^{(\ga)}$
commutes with $A_X = \exp \ker\ga.$ It follows that $\xi$ uniquely extends to a character of
$G^{(\ga)}$ which is $1$ on $M_X.$ This extension is also denoted by $\xi.$ 
The following result is straightforward.

\begin{lemma}
For all $\nu, \mu \in \faPdc,$ 
$$
\ge_\nu|_{G^{(\ga)}} = \xi \ge^{(\ga)}_\nu \quad{\rm and}\;\; \gf_\mu|_{G^{(\ga)}} = \gf^{(\ga)}_\mu.
$$
\end{lemma}

As the character $\xi$ is only non-trivial on the center of $G^{(\ga)}$ its role 
is easily understood in the calculations that follow. From the lemma it follows
that $\gf_\mu\ge_{\nu - \mu}$ restricts to $\xi \gf_\mu^{(\ga)}\ge^{(\ga)}_{\nu - \mu}.$ 
This suggests that $J_{\nu,\mu}$ and $J^{(\ga)}_{\nu,\mu}$ might be related.

In fact, we will show that the following is valid.

\begin{prop} 
\label{p: J and J ga}
Let $u \in H_\gs^\infty$. Then, for $(\nu,\mu) \in \faPdc\times \faPdc,$ 
$$
{}_u\widetilde J_{\nu,\mu}|_{N_\ga} = {}_u\widetilde J^{(\ga)}_{\nu,\mu}.
$$ 
\end{prop}

For the proof we need some preparation. First we will describe a direct relationship 
between $\gf_\mu j_{\nu - \mu}$ and $\gf^{(\ga)}_\mu j^{(\ga)}_{\nu - \mu}.$

\begin{lemma}
\label{l: long lemma j nu mu}
Let $\Omega \subset \faPdc \times \faPdc$ be a bounded connected open subset which contains 
a point $(\mu, \nu)$ such that $\Re (\nu - \mu)$ is $\bar P$-dominant.  There
exists a positive integer $r$ such that the following assertions are valid.
\begin{enumerate}
\itema
The map $(\mu,\nu) \mapsto j_{\nu - \mu}$ is holomorphic  $\Omega \to 
C^{-r}(K/K_P:\gs_P).$ 
\itemb
For every $f \in C^{r}(K/K_P:\gs_P)$ 
the function
$N_X \times \Omega \to C^r(K^{(\ga)}/K_P:\gs_P),$ 
$$
(n_X, \nu,\mu) \mapsto  L_{n_X}^{-1}( \gf_{\bar \mu} ) 
L_{n_X}^{-1} (f_{P, \gs, -\bar \nu} )|_{K^{(\ga)}} 
$$ 
is smooth, and holomorphic in the variable $(\nu,\mu)$ from $\Omega$.
\itemc
For every $f \in C^{r}(K/K_P:\gs_P)$ with support contained in $K \cap N_P \bar P,$
$$
\inp{\gf_\mu j_{\nu - \mu}  }{  f  } 
 =
\int_{N_X} \chi(n_X)
 \inp{ j^{(\ga)}_{\nu - \mu}|_{K^{(\ga)}} }{  L_{n_X}^{-1} (\gf_{\bar \mu})
L_{n_X}^{-1} (f_{\bar P, \gs, -\bar \nu} )|_{K^{(\ga)}}  }\; \;d n_X
$$ 
for all $(\nu, \mu) \in \Omega.$ 
\end{enumerate}
\end{lemma}

\proof
Assertion (a) follows from Corollary \ref{c: holomorphy j Q}.
Assertion (b) is obvious. We address (c). 
The set  $\Omega_0$ of points $(\mu,\nu) \in \Omega$ such  
the real part $\Re (\nu - \mu)$ is strictly $\bar P$-dominant is non-empty and open.
Let $(\nu,\mu) \in \Omega_0;$ then it 
follows that for $f \in C^r(K/K_P:\gs_P)$ the function
$\inp{\gf_\mu j_{\nu - \mu}}{f}_\gs$ is integrable over $K/K_P.$  
We now observe that the function $F: x \mapsto \inp{\gf_\mu(x) j_{\nu - \mu}(x)}{f_{Q,\gs, -\bar \nu}(x)}_\gs$ on $G$ 
is right $M_P \bar N_P$-invariant, and satisfies $R_a F = a^{2 \rho_P} F$ 
for all $a \in A_P.$In view of Lemma \ref{l: integral formula ga X} it follows that
\begin{eqnarray*}
\inp{\gf_\mu j_{\nu - \mu}  }{  f  } &=&
\int_{N_X \times K^{(\ga)}}
 \inp{ \ge_{\nu - \mu}(n_X k_\ga)
 }{  
 ( \gf_{\bar \mu} f_{\bar P,\gs,-\bar \nu})(n_X k_\ga)} \; d n_X \,d k_\ga = 
\\ 
& =  &
 \int_{N_X} \chi(n_X) 
 \int_{K^{(\ga)}}
 \inp{
\ge^{(\ga)}_{\nu -\mu}(k_\ga)}{
 \gf_{\bar \mu}(n_X k_\ga) f_{\bar P,\gs,-\bar  \nu})(n_X k_\ga)}\; dk_\ga d n_X 
\\
&=&
 \int_{N_X} \chi(n_X) 
\inp{
j_{\nu -\mu}^{(\ga)}|_{K^{(\ga)} }}{  L_{n_X^{-1}}
( \gf_{\bar \mu}) L_{n_X^{-1}} (f_{\bar P,\gs,-\bar  \nu})|_{K^{(\ga)}}
}\;  d n_X.
\end{eqnarray*}
This establishes the identity of (c) for $(\nu, \mu) \in \Omega_0.$ 

If $f$ satisfies the mentioned support condition, it follows
that there exists a compact set $\cK\subset N_X$ such that 
$$ 
L_{n_X^{-1}}  (f_{\bar P,\gs,-\bar  \nu})|_{K^{(\ga)}} = 0
$$ 
for all $\nu \in \faPdc$ and all $n_X\in N_X\setminus \cK.$ From this it
is readily seen that the expressions on both sides of the equation in (c) 
are holomorphic functions of $(\nu,\mu) \in \Omega.$ The full result now follows by analytic continuation.
\qed

For $\gf_X\in C^{\infty}_c(N_X)$ and 
$\gf_\ga \in C^{\infty}_c(N_\ga),$ 
define $\gf_X \otimes \gf_\ga \in C_c^\infty(N_P)$ by 
$$ 
\gf_X \otimes \gf_\ga (n_X n_\ga) = \gf_X(n_X)\gf_\ga(n_\ga).
$$ 
We note that for $\nu \in \faPdc.$ According to the definition,
$$ 
{}_u^\backprime  T_\nu( \gf_X \otimes \gf_\ga) \in C^r(G/\bar Q:\gs :\nu)
$$ 
is given by
$$
  {}_u^\backprime T_\nu( \gf_X \otimes \gf_\ga)(n_X n_\ga v man) = a^{-\nu+\rho_Q}
 \gf_X(n_X)\gf_\ga(n_\ga) \;\gs(m)^{-1} u
 $$ 
 for $(n_X, n_\ga , m a \bar n ) \in N_X\times N_\ga  \times M_P A_P \bar N_Q$
 and by  $ {}_u^\backprime T_\nu ( \gf_X \otimes \gf_\ga)= 0$ on $G \setminus N_P v \bar Q.$  
 We recall that 
 $$ 
 T_\nu( \gf_X \otimes \gf_\ga) = {}_u^\backprime T_\nu( \gf_X \otimes \gf_\ga)|_K
 \in C^r(K/K_P: \gs_P).
 $$
 
 The map ${}_u T_\nu^{(\ga)}$ is defined similarly for the group $G^{(\ga)};$ 
 note that $N_{P^{(\ga)} } = N_\ga.$ For $\gf \in C^r_c(N_\ga)$, $u \in H_\gs^\infty$ 
 and $\nu \in \faPdc$ 
 we define  $ {}_u^\backprime T^{(\ga)}_\nu (\gf)\in C^r(G^{(\ga)}/\bar Q^{(\ga)}:\gs: \nu)$
 by
 $$
 {}_u^\backprime T^{(\ga)}_\nu (\gf)(n_\ga v m a \bar n) = 
  a^{-\nu + \rho_{Q^{(\ga)}}} \gs(m)^{-1} \gf(n_\ga) u
  $$ 
  for $n_\ga v m a \bar n \in N_\ga v M_P A_P \bar N_{Q^{(\ga)}}$ and by 
 ${}_u^\backprime T^{(\ga)}_\nu (\gf) = 0$ on $G^{(\ga)}\setminus N_\ga v \bar Q^{(\ga)}.$ 
Finally, ${}_u T^{(\ga)}_\nu: C^r_c(N_\ga) \to C^r(K^{(\ga)}/K_P : \gs_P)$ 
is defined by 
$$ 
{}_u T_\nu^{(\ga)} (\gf) = {}_u^\backprime T_\nu (\gf)|_{K^{(\ga)}},\qquad (\gf \in C_c^r(N_\ga)).
$$

 \medno
 Suppose now that $R$ is either of the parabolic subgroups $P$ and $Q$.
 In this situation, the natural multiplication map $N_X \times K^{(\ga)} \times_{K_P} \bar R \to G$ is an open embedding.
If  $\gf_X \in C_c^r (N_X)$, $\psi_\ga \in C^r(K^{(\ga)} / K_P:\gs_P)$ and $\nu \in \faPdc,$
we define the $C^r$-function $S_{R,  \nu}( \gf_X \otimes  \psi_\ga):  G \to H_\gs$ by 
$$
{}^\backprime S_{R, \nu}( \gf_X \otimes  \psi_\ga) 
(n_X k_\ga ma \bar n ) = a^{-\nu +\rho_R}
\gf_X(n_X) \psi_\ga(k_\ga ) \gs(m)^{-1} u,
 $$
for $(n_X, k_\ga, ma \bar n) \in  N_X \times K_\ga \times \bar R,$
and by $ {}^\backprime  S_{R,\nu}( \gf_X \otimes  \psi_\ga) = 0$ on $G\setminus N_X K_\ga \bar R.$
Thus, $S_{R, \nu} ( \gf_X \otimes  \psi_\ga) \in C^r(G/R:\gs: \nu).$
As  before, we define 
$$
S_{R, \nu} ( \gf_X \otimes  \psi_\ga) =  {}^\backprime S_{R, \nu} ( \gf_X \otimes  \psi_\ga)|_K
$$

\begin{lemma}
Let $u \in H_\gs^\infty.$ Then for $\gf_X \in C^r_c(N_X)$ and
$\gf_\ga \in C^r_c(N_\ga),$ 
$$ 
{}_uT_{\nu}(\gf_X \otimes \gf_\ga)
= S_{Q,\nu} \left 
(\gf_X \otimes {}_uT^{(\ga)}_\nu(\gf_\ga)
\right),
$$
for $\nu\in \faPdc.$
\end{lemma}

\proof 
We will use the notation $\gs_\nu$ for the character of $\bar Q$ given by 
$\gs_\nu(ma\bar n) = a^{\nu - \rho_Q}\gs(m).$ 

Let $k \in K$ and suppose ${}_u T_{\nu}(\gf_X \otimes \gf_\ga) (k) \neq 0.$ 
Then there exist $n_X \in N_X$ and $n_\ga \in N_\ga$ such that
$k = n_X n_\ga v  \bar q$ with $\bar q \in \bar Q$ and 
$$
{}_u T_{\nu}(\gf_X \otimes \gf_\ga) (k) = \gs_\nu(\bar q)^{-1}(\gf_X(n_X)\gf_\ga(n_\ga) u ).
$$
In particular, $\gf_\ga(n_\ga)  \neq 0.$ Write $n_\ga v = k_\ga \bar q_\ga,$ with 
$k_\ga\in K^{(\ga)}$ and $ \bar q_\ga \in \bar  \fq^{(\ga)}.$
Then 
$k = n_X k_\ga \bar q_\ga  \bar q.$ Hence,
\begin{eqnarray}
\nonumber
\lefteqn{S_{Q,\nu}(\gf_X \otimes {}_u T_\nu^{(\ga)} \gf_\ga)(k) }\\
&= & \gs_\nu(\bar q_\ga \bar q)^{-1} \gf_X( n_X)   \left( {}_u T^{(\ga)}_\nu \gf_\ga(k_\ga)\right)\nonumber\\
&= & \gs_\nu(\bar q_\ga \bar q)^{-1} \gf_X( n_X)\gs_\nu(\bar q_\ga)  \gf_\ga(n_\ga) u \nonumber\\
&=& \gs_\nu(\bar q)^{-1} \gf_X( n_X)\;  \gf_\ga(n_\ga)u\nonumber\\
&=& {}_u T_\nu(\gf_X \otimes \gf_\ga)(k). \label{r: long display} \end{eqnarray}
Conversely, suppose that 
$$
S_{Q,\nu}(\gf_X \otimes {}_u T_\nu^{(\ga)} \gf_\ga)(k) \neq 0.
$$
Then there exist $n_X \in N_X, k_\ga \in K_\ga$ such that $ k= n_X  k_\ga  \bar q$ with $\bar q \in \bar Q.$ 
Moreover,
$$
S_{Q,\nu}(\gf_X \otimes {}_u T_\nu^{(\ga)} \gf_\ga)(k) = \gs(\bar q)^{-1}\gf_X(n_X)\; {}_u T_\nu^{(\ga)}(\gf_\ga)(k_\ga).
$$
In particular, $ {}_u T_\nu^{(\ga)}(\gf_\ga)(k_\ga) \neq 0.$ 
Hence, there exist $n_\ga \in N_\ga$ and $\bar q_\ga  \in \bar Q^{(\ga)}$
such that $k_\ga = n_\ga v \bar \fq^{(\ga)}.$ 
Now
$$
{}_u T_\nu^{(\ga)}\gf_\ga(k_\ga) = \gs_\nu(\bar q_\ga)^{-1} \gf_\ga(n_\ga).
$$
We now have $k = n_X k_\ga \bar q$ and $k_\ga = n_\ga v \bar q_\ga$ so that
$k =n_X n_\ga v \bar q_\ga \bar q.$ It follows that
\begin{eqnarray*}
\lefteqn{{}_u T_\nu (\gf_X \otimes \gf_\ga)  (k)}\\
&=& \gs_\nu(\bar q_\ga \bar q)^{-1} \gf_X(n_X) \gf_\ga(n_\ga) u\\
&=& \gs_\nu(\bar q)^{-1}
\gf_X(n_X) ({}_uT_\nu^{(\ga)}\gf_\ga)(k_\ga)  \\
&=& 
 S_{Q,\nu} (\gf_X \otimes {}_u T_\nu^{(\ga)}\gf_\ga)(k) .
\end{eqnarray*}
Thus,  we have shown that for $k \in K,$ ${}_u T_\nu(\gf_X\otimes \gf_\ga)(k)$ is non-zero
if and only if $S_{Q,\nu}(\gf_X \otimes  {}_u T_\nu^{(\ga)}\gf_\ga)(k)$ is non-zero,
and that the desired  equality is valid at such $k.$ In the remaining points $k$, both functions are zero, hence also equal.
\qed

\begin{lemma}
\label{l: second lemma cond Omega}
Suppose that $\Omega$ is a bounded connected open subset of $\faPdc,$
containing a point $\nu$ such that $\Re\nu $ is $\bar Q$-dominant. Let $t
$ dominate both the order
of  $A(\bar P , \bar Q, \gs, \dotvar)$ over $\Omega$ and the order of $A^{(\ga)}(\bar P^{(\ga)}, \bar Q^{(\ga)}, \gs, \nu)$
over $\Omega.$ 
Then for $\gf_X \in C_c^{r+t} (N_X)$ 
and $\psi \in C_c^{r+t}(K_\ga /K_P:\gs_P),$ 
$$
A(\bar P, \bar Q, \gs, \nu) S_{Q, \nu} (\gf_X \otimes \psi_\ga) 
=
S_{P, \nu} (\gf_X  \otimes  A^{(\ga)} (\bar P^{(\ga)},  \bar Q^{(\ga)}, \gs, \nu)\psi_\ga).
$$
as meromorphic functions of $\nu\in \Omega$ with values in 
$C^r(K/K_P:\gs).$ 
\end{lemma}
\proof
The statements about meromorphy are well-known, and serve here to allow meromorphic
continuation of identities. Let $\Omega_0$ be the set of $\nu \in \Omega$ 
such that $\Re \nu$ is $\bar Q$-dominant.Then $\Omega_0$ is open and non-emtpy.
For $\nu \in \Omega_0$ the intertwining operators are given by the familiar integral formulas.
We use the abbreviated notation $A(\nu)$ and $A^{(\ga)}(\nu)$ for the above mentioned intertwining operators.
Then it suffices to show that 
$$ 
A(\nu) \bp S_{Q,\nu} (\gf_X \otimes  \psi_\ga) = \bp S_{P,\nu} (\gf_X \otimes A^{(\ga)}(\nu) \psi_\ga)
$$ 
at each point $g: = n_X k_\ga \bar p \in N_X K^{(\ga)} \bar P.$  Since the elements
on both sides of the equation belong to $C^r(G/\bar P:\gs:\nu),$ we may as well
assume that $g = n_X k_\ga.$ Then 
$$
A(\nu) \bp S_{P,\nu} (\gf_X \otimes  \psi_\ga) (g) =  
 \int_{\bar N_\ga}  \bp S_{P,\nu} (\gf_X \otimes  \psi_\ga)(n_X k_\ga \bar n_\ga) d \bar n_\ga.
$$ 
Since $\bar n_\ga \in G^{(\ga)} $ we may write 
$\bar n_\ga = \kappa(\bar n_\ga) \bar p(n_\ga)$ with $\kappa(\bar n_\ga) \in K^{(\ga)}$
and $\bar p(\bar n_\ga)\in \bar P^{(\ga)}$ smoothly depending on $\bar n_\ga.$ 
Therefore,
\begin{eqnarray*} 
A(\nu) \bp S_{P,\nu} (\gf_X \otimes  \psi_\ga) (g) 
& = & 
 \int_{\bar N_\ga}  \gs_\nu(\bar p(\bar n_\ga))^{-1} 
 \bp S_{P,\nu} (\gf_X \otimes  \psi_\ga)(n_X k_\ga \kappa(\bar n_\ga) )
 d \bar n_\ga\\
 &=& 
  \int_{\bar N_\ga}  \gs_\nu(\bar p(\bar n_\ga))^{-1} \gf_X(n_X) \psi_\ga(k_\ga \kappa(\bar n_\ga)) \; d \bar n_\ga\\
&=&  
\gf_X (n_X) \int_{\bar N_\ga} \gs_\nu(\bar p(\bar n_\ga))^{-1}\psi_\ga(k_\ga \kappa(\bar n_\ga)) \; d \bar n_\ga
\\
&=& 
 \gf_X (n_X) \int_{\bar N_\ga} \psi_\ga(k_\ga \bar n_\ga) \; d \bar n_\ga
= 
 \gf_X (n_X) [A^{(\ga)}(\nu) \psi_\ga](k_\ga) \\
 &=&  
 \bp S_{Q, \nu} (\gf_X \otimes A^{(\ga)}(\nu)\psi_\ga)(g).
 \end{eqnarray*}
 \qed
 
 {\em Proof of Prop.~\ref{p: J and J ga}.\ } Suppose that $u \in H_\gs^\infty$ is fixed.
 It suffices to prove the identity for $(\nu,\mu) \in \Omega$, 
 where $\Omega$ is an open subset of $\faPdc \times \faPdc$ satisfying the conditions of Lemmas \ref{l: long lemma j nu mu} and \ref{l: second lemma cond Omega}.
 
Let $\gf_X\in C_c^\infty(N_X)$ and $\gf_\ga \in C^\infty_c(N_\ga);$ 
then 
\begin{eqnarray}
\lefteqn{
\int_{N_X\times N_\ga}  {}_u \widetilde J_{\nu,\mu}  (n_X, n_\ga) \;\gf_X(n_X)\; \gf_\ga(n_\ga)\; dn_X \,dn_\ga
}
\label{e: first expression int nX nga}
\\
&=& 
{}_u J_{\nu,\mu} (\gf_X \otimes \gf_\ga)
\nonumber
\\
&=&
\inp{J_{\nu,\mu}}{
{}_uT_{-\bar \nu}(\bar \gf_X \otimes \bar \gf_\ga)}
\nonumber
\\
&=& 
\inp{J_{\nu,\mu}}{S_{Q,-\bar\nu}(\bar  \gf_X\otimes {}_u T^{(\ga)}_{-\bar \nu}(\bar \gf_\ga)}
\nonumber
\\
&=& 
\inp{A(\bar Q, \bar P, \gs, \nu) (\gf_\mu j_{\nu-\mu})} {S_{Q, -\bar \nu}(\bar  \gf_X\otimes {}_u T^{(\ga)}_{-\bar \nu}(\bar \gf_\ga)}\nonumber
\\
&=&  \inp{\gf_\mu j_{\nu-\mu}}{A(-\bar \nu) S_{Q, -\bar\nu}(\bar \gf_X \otimes \bar \psi_{\ga, -\bar \nu})}
\nonumber
\\
&=&  \inp{\gf_\mu j_{\nu-\mu}}{S_{P, -\bar\nu}(\bar \gf_X \otimes A^{(\ga)}(-\bar \nu) \bar \psi_{\ga, -\bar \nu})}
\label{e: last expression gf mu j}
\end{eqnarray}
where we have written $A(-\bar \nu) = A(\bar P, \bar Q, \gs, -\bar \nu)$, $A^{(\ga)}(-\bar \nu) = A^{(\ga)}(\bar P^{(\ga)} ,
\bar Q^{(\ga)}, \gs, -\bar \nu)$ and $\bar \psi_{\ga, -\bar \nu} = T^{(\ga)}_{-\bar \nu}(\bar \gf_\ga).$
Put 
$$
F = S_{P, -\bar\nu}(\bar \gf_X \otimes A^{(\ga)}(-\bar \nu) \bar \psi_{\ga, -\bar \nu}).
$$
Then $F \in C^r(K/K_P:\gs_P).$  Applying Lemma \ref{l: long lemma j nu mu} with $F$ in place of $f,$ 
we find that (\ref{e: last expression gf mu j}) equals
\begin{equation}
\label{e: new expression gf mu j with F}
\inp{\gf_\mu j_{\nu-\mu}}{F}=
\int_{N_X} \chi(n_X)
 \inp{ j^{(\ga)}_{\nu - \mu}|_{K^{(\ga)}} }{  L_{n_X}^{-1} (\gf_{\bar \mu})
L_{n_X}^{-1} (F_{P, \gs, -\bar \nu} )|_{K^{(\ga)}}  }\; \;d n_X.
\end{equation}
Now 
$$
F_{P,\gs, -\bar \nu}(n_X k_\ga) = \bar \gf_X(n_X)  A^{(\ga)}(-\bar \nu) \bar \psi_{\ga, -\bar \nu}(k_\ga),
$$
so that (\ref{e: new expression gf mu j with F}) equals
\begin{equation}
\label{e: last expression inp gf mu j and F}
\int_{N_X} \chi(n_X) \gf_X(n_X) 
 \inp{ j^{(\ga)}_{\nu - \mu}|_{K^{(\ga)}} }{  L_{n_X}^{-1} (\gf_{\bar \mu})
  A^{(\ga)}(-\bar \nu) \bar \psi_{\ga, -\bar \nu}|_{K^{(\ga)}} }\; \;d n_X.
 \end{equation}
The equality of (\ref{e: first expression int nX nga}) with (\ref{e: last expression inp gf mu j and F})
for arbitrary $\gf_X \in C^\infty_c(N_X)$ implies that for every $n_X \in N_X,$ 
$$
\int_{N_\ga}  {}_u \widetilde J_{\nu,\mu}  (n_X, n_\ga) \; \gf_\ga(n_\ga)\; dn_\ga 
=
\chi(n_X) 
 \inp{ j^{(\ga)}_{\nu - \mu}|_{K^{(\ga)}} }{  L_{n_X}^{-1} (\gf_{\bar \mu})
  A^{(\ga)}(-\bar \nu) \bar \psi_{\ga, -\bar \nu}|_{K^{(\ga)}} }.
$$
Substituting $n_X = e,$ we find that
\begin{eqnarray*}
\lefteqn{\int_{N_\ga}  {}_u \widetilde J_{\nu,\mu}  (e, n_\ga) \; \gf_\ga(n_\ga)\; dn_\ga }\\
&=&
 \inp{ j^{(\ga)}_{\nu - \mu}|_{K^{(\ga)}} }{ (\gf_{\bar \mu})
  A^{(\ga)}(-\bar \nu) \bar \psi_{\ga, -\bar \nu}|_{K^{(\ga)}} }\\
&=&
 \inp{A^{(\ga)}(\bar Q^{(\ga)}, \bar P^{(\ga)}, \gs, \nu) \gf^{(\ga)}_\mu j^{(\ga)}_{\nu - \mu}}{ {}_u T^{(\ga)}_{-\bar \nu}(\bar \gf_\ga)}
 \\
&=& 
{}_u J^{(\ga)}_{\nu,\mu} (\gf_\ga) = \int_{N_\ga} {}_u \widetilde J_{\nu,\mu}(n_\ga) \gf_\ga(n_\ga ) dn_\ga.
\end{eqnarray*}
As this is valid for any $\gf_\ga \in C^{\infty}_c(N_\ga),$ we conclude that
$$ 
 {}_u \widetilde J_{\nu,\mu}  (e, n_\ga) = {}_u \widetilde J^{(\ga)}_{\nu,\mu}(n_\ga), 
$$
for all $n_\ga \in N_\ga.$ 
\qed

\section{Comparison of $B$ with $B^{(\ga)}$, proof of Lemma \ref{l: comparison B  and B ga}}
\label{s: comparison B}
We retain the notation of the previous section. 
\begin{lemma}
For every $u \in H_\gs^\infty$ the assignment  $\nu \mapsto {}_u \widetilde J_{\nu,0} $ 
is meromorphic $\faPdc \to C^\infty(N_P)$ and given by 
\begin{equation}
\label{e: formula for u J nu 0}
{}_u \widetilde J_ {\nu,0}(n) =  \chi(n)\inp{B(\bar Q, \bar P, \gs, \nu)\eta}{u}_\gs, \qquad (n\in N_P),
\end{equation}
as an identity of meromorphic functions of the variable $\nu.$ 
\end{lemma}

It follows from the assertion about the singular set in Theorem \ref{t: tilde J smooth}
that $\nu \mapsto {}_u \widetilde J_{\nu,0} $ is a genuine meromorphic function
of $\nu$ with values in $C^\infty(N_P).$
On the other hand, by definition, for regular values of $\nu\in \faPdc$ 
the following identity of elements of $C^{-\infty}(G/\bar Q:\gs:\nu)_\chi$ is valid:
$$
J_{\nu,0} = A(\bar Q, \bar P, \gs, \nu) \gf_0 j(\bar P, \gs, \eta) = j(\bar Q , \gs, \nu) B(\nu) \eta,
$$
where $B(\nu) = B(\bar Q,\bar P, \gs, \nu).$
From this equality combined with Lemma \ref{l: chi invariance implies continuous} and 
\cite[Thm.~ 8.6]{Bunitemp} it follows that for 
regular $\nu,$ the generalized function  $J_{\nu,0} \in C^{-\infty}(G/\bar Q:\gs:\nu)$ is continuous
$H_\gs^{-\infty}$-valued on the open set $N_P v \bar Q$  where it is given by 
$$
J_{\nu, 0}|_{N_P v \bar Q}(n v m a \bar n) = \chi(n) a^{-\nu + \rho_Q} \gs^{-1}(m) B(\nu) \eta, \qquad (n \in N_P, ma \bar n \in \bar Q),
$$ 
in the sense that the identity is valid after testing with any function $\gf$ 
from $C^\infty(G/\bar Q:\gs: -\bar \nu)$ whose support is contained
in $N_P v \bar Q, $ i.e., 
$$
\inp{J_{\nu, 0}}{\gf} = \int_{K} \inp{J_{\nu, 0}|_{N_P v \bar Q}}{\gf}_\gs(k) \; dk.
$$

Let $u \in H_\gs^\infty;$ then for any  $f \in C^\infty_c(N_P)$ 
the function $\gf = {}_uT_{-\bar \nu}(\bar f)$ is of this type. The function 
$F:= \inp{J_{\nu, 0}}{{}_uT_{-\bar \nu}(\bar f)}_\gs$ is a continuous function 
$G \to \C$ which is right $M_Q \bar N_Q$-invariant and transforming
according to the rule $R_a F = e^{2\rho_Q} F.$ Hence, by Lemma \ref{l: int over K is int over NPv} its integral over
$K$ is given by 
\begin{eqnarray} 
\label{e: first expression J nu 0 in array}
\inp{J_{\nu, 0}}{{}_uT_{-\bar \nu}(\bar f)} &= &
\int_{N_P} \inp{J_{\nu, 0}(nv)}{{}_uT_{-\bar \nu}(\bar f) (nv)}\; dn \\
&=& 
\int_{N_P} \inp{J_{\nu, 0}(nv)}{\bar f(n) u }\; dn \nonumber \\
&=&  
\int_{N_P} \chi(n) \inp{B(\nu)\eta}{u} f(n)\; dn.\label{e: last expression J nu 0 in array}
\end{eqnarray}
The expression on the left-hand side of 
(\ref{e: first expression J nu 0 in array}) is equal to 
\begin{equation}\label{e: second array with J nu 0}
\inp{J_{\nu, 0}}{{}_uT_{-\bar \nu}(\bar f)} =
{}_u J_{\nu, 0}(f) = 
\int_{N_P} f(n) {}_u \widetilde J_{\nu, 0}(n) \; dn
\end{equation}
see the text preceding Theorem \ref{t: tilde J smooth}. It follows  that the integral in \ref{e: second array with J nu 0} is equal to the integral
in (\ref{e: last expression J nu 0 in array})  for all $f\in C^\infty_c(N_P).$
Since $ {}_u \widetilde J_{\nu, 0}$ and $\chi$ are continuous functions on $N_P$ the desired identity follows.
\qed
\medbreak 

{\em End of proof  Lemma \ref{l: comparison B and B ga}.\ }
All arguments presented so far in this section are valid for the triple 
$G^{(\ga)}, P^{(\ga)}, Q^{(\ga)} = \bar P^{(\ga)}$ in place of $G, P, Q.$   
In particular, if
$u\in H_\gs^\infty,$ then the function ${}_u\widetilde J^{(\ga)}_{\nu,0} \in C^\infty(N_\ga)$ 
depends meromorphically on $\nu \in \faPdc$ and is given by
\begin{equation}
\label{e: formula for tilde J ga}
{}_u\widetilde J^{(\ga)}_{\nu,0}(n) = \chi^{(\ga)}(n) \inp{B^{(\ga)}(\bar Q^{(\ga)}, \bar P^{(\ga)}, \gs, \nu)\eta}{u}
\end{equation}
for $n \in N_{P^{(\ga)}} = N_\ga.$ From Proposition \ref{p: J and J ga} it follows that
$$
{}_u\widetilde J_{\nu,0}(e)  =  {}_u\widetilde J^{(\ga)}_{\nu,0}(e)
$$
as meromorphic functions of $(\nu, \mu) \in \faPdc \times \faPdc.$  Combining this with
(\ref{e: formula for tilde J ga}) and (\ref{e: formula for u J nu 0}) we obtain that
$$
\inp{B(\bar Q, \bar P, \gs, \nu)\eta}{u} = 
\inp{B^{(\ga)}(\bar Q^{(\ga)}, \bar P^{(\ga)}, \gs, \nu)\eta}{u}.
$$
Since this holds for every $u \in H_\gs^\infty$ the proof is complete.
\qed

\section{The C-functions and the Maass--Selberg relations}
\label{s: c functions MS}
From now on we assume that $\ft \subset \fm$ is maximal abelian, so that
$\fh = \ft \oplus \fa$ is a maximally split Cartan subalgebra of $\fg.$ 
Let $R \in \cP.$ We write ${}^*\fh_R$ for the orthocomplement of $\fa_R$ in $\fh.$ This is a maximally split Cartan subalgebra of $\fm_R,$ which decomposes as 
$$
{}^*\fh_R = \ft \oplus {}^*\fa_R.
$$

We consider the ($\tau$-spherical) Whittaker integral 
$\Wh(R, \psi, \nu),$ for $\psi \in \cA_{2,R} = \cA_2(\tau: M_R/M_R \cap v_R^{-1} N_0 v_R:\chi_R).$ 
If $R$ is non-cuspidal, then $\cA_{2,R} = 0$ so that the Whittaker integral is trivial. 
Therefore, we assume  $R$ to be cuspidal.

 Let $\gL \in {}^*\fhRdc$ be the infinitesimal
character of a representation of the discrete series of $M_R.$ For  $\ge > 0$ we define $\faRd(\ge) = \{\nu \in \faRdc\mid |\Re(\nu)| < \ge\}.$ For $r, \ge > 0$ 
we consider the set $\twohol(\gL, \faR, \ge , r, \tau)$ of families of type 
$\twohol$ as defined in  \cite[\S 7]{Bcterm}. This set consists of  
families $(f_\nu)_{\nu \in \faRd(\ge)}$ of functions $f_\nu \in C^\infty(\tau: G/N_0: \chi)$ 
such that 
\begin{enumerate}
\itema the function $\nu \mapsto f_\nu$ is holomorphic $\faRd(\ge) \to C^\infty(\tau: G/N_0: \chi);$
\itemb $Z f_\nu = \gg(Z,\gL + \nu) f_\nu \quad
 (\nu \in \fa_R^*(\ge), 
 \;Z \in\fZ)$;
\itemc  for every $u \in U(\fg)$ there exist $C > 0$ and $N\in \N$ such that
$$
|L_u f_\nu(x)|\leq C_N |(x,\nu)|^N e^{-\rho H(x) + r|\Re\nu||H(x)|}, \qquad ((\nu, x) \in \faRd(\ge) \times G).
$$ 
\end{enumerate}
Here $|(x,\nu)|:= (1 +|H(x)|)(1 + |\nu|).$ 

\begin{lemma}
\label{l: Wh is twohol}
Let for $\gs \in \dMRds.$ There exist constants $\ge > 0, r >0$ such that
for every $\psi_{2,R,\gs}\in \cA_{2,R,\gs}$
the family $\faRdc \to C^\infty(\tau: G/N_0: \chi),$ 
$$
\nu \mapsto \Wh(R, \psi, \nu)
$$  
belongs to $\twohol(\gL, \faR, \ge , r, \tau);$ here $\gL$ denotes the infinitesimal character
of $\gs.$ 
\end{lemma}
\proof
See  \cite{Bcterm}.
\qed

In particular, the function $\Wh(R, \psi, \nu)$ belongs to 
the space $\cA(\tau :  G/N_0:\chi)$ of $\tau$-spherical tempered Whittaker functions. 
More precisely, it follows from Lemma \ref{l: Wh is twohol} and  from the theory developed in  \cite{Bcterm} that for $\nu \in i\fad,$ 
its constant term along a parabolic subgroup $Q\in \cP,$ defined as in 
 \cite[\S 1.4]{HCwhit},  is denoted by $\Wh_Q(R, \psi, \nu).$  

It follows from the theory in  \cite{Bcterm} that for $\ge > 0$ 
sufficiently close to zero this constant term extends to a holomorphic 
function  $\Wh_Q(R,\psi, \dotvar)$ on $\faRd(\ge)$ with values in $C^\infty(\tau: G/N_0: \chi).$

For this constant term to be non-zero for any particular value of $\nu \in i\faRd,$ 
the parabolic subgroup $Q$ needs to be
 standard  (see  \cite{HCwhit} and  \cite{Bcterm}), and there needs to be a  standard  parabolic subgroup $P$ contained in $Q$ such that 
$P \sim R$ (meaning that $\faR$ and $\faP$ are conjugate under $W(\fa)).$ 
In this case, if $Q\not\sim R$ the function $m \mapsto  R_a[\Wh_Q(R,\psi,\nu)]$
is perpendicular to $L^2_{\ds}(\tau: M_Q/M_Q\cap N_0:\chi)$ for all $a \in A_Q.$ 
If $Q\sim R,$ then the function $m \mapsto  R_a[\Wh_Q(R,\psi,\nu)]$
belongs to  $L^2_{\ds}(\tau: M_Q/M_Q\cap N_0:\chi)$ for all $a \in A_Q.$ In this
case, the precise form of the constant term is given in the following result.

If $Q,R \in \cP$ then $W(\faQ|\faR)$ denotes the set of $s\in \Hom(\faR, \faQ)$ for which there exists
a $w \in W(\fa)$ such that $s = w|_{\faR}.$ 

\begin{thm} 
\label{t: asymp char C functions} Let $R \in \cP.$ Then for $\ge> 0$ sufficiently small and 
for every $Q\in \cPst$ with $Q\sim  R$ there exist unique meromorphic
functions $C_{Q|R}(s, \dotvar)$ on $\faRdc(\ge)$ with values in $\Hom(\cA_{2,R}, \cA_{2,Q})$ 
such that for all $\nu \in \faRdc(\ge)$ and $\psi \in \cA_{2,R},$ we have 
\begin{equation}
\label{e: characterization C functions}
\Wh_Q(R, \psi, \nu) (ma) = \sum_{s \in W(\faQ|\faR)} a^{s\nu} C_{Q|R}(s,\nu)(\psi)(m),
\end{equation}
for $m\in M_Q$ and $a \in A_Q,$ as meromorphic functions of $\nu.$ 
\end{thm}

 From now on we will assume that $\ge > 0$ is sufficiently small.
 We proceed to obtain more detailed information on the $C$-functions from
 their characterization through Theorem \ref{t: asymp char C functions}.

 \begin{lemma}\label{l: C Q Q}
 Let $Q\in \cPst.$ Then for each $\gs\in \dMQds$
 appearing as an isotype in $\cA_2(\tau: M_Q/M_Q \cap N_0:\chi_Q),$ 
 and for all $T \in C^{\infty}(\tau: K/K_Q: \gs_Q) \otimes \ccH_{\gs, \chi_Q}^{-\infty},$ 
 we have 
\begin{equation}
C_{Q|Q}(1 : \nu) \psi_T = \psi_{[A(Q, \bar Q, \gs, -\nu) \otimes I]T}
\end{equation}
as an identity of meromorphic functions of $\nu \in \faQd(\ge).$  
\end{lemma}
\proof 
By linearity we may assume that $T = \gf \otimes \eta,$ with 
$\gf \in C^\infty(\tau: K/K_Q: \gs_Q)$ and $\eta \in \ccH^{-\infty}_{\gs, \chi_Q}.$ 
Then by Definition \ref{d: defi j Q} we have, for all $\nu \in \faQdc$ with $\inp{\Re \nu}{\ga}> 0$ 
for all $\ga \in \gS(\faQ, \bar\fn_Q),$ and for all $(m,a)\in M_Q \times A_Q,$ 
\begin{eqnarray*}
\lefteqn{\Wh(Q, \psi_T ,  \nu)(ma)
= \int_{K/K_Q} \inp{\gf_{\bar Q,\gs, -\nu}(ma k )}{j(\bar Q, \gs, \bar \nu, \eta)(k)}\; dk}\\
&=& \int_{N_Q}  \inp{\gf_{\bar Q,\gs, -\nu}(ma n )}{j(\bar Q, \gs, \bar \nu, \eta)(n)}\; dn \\
&=& a^{\nu - \rho_Q}\; \int_{N_Q}  \inp{\gf_{\bar Q,\gs, -\nu}(m n )}{j(\bar Q, \gs, \bar \nu, \eta)(a^{-1}na )} \;dn\\
&=&
a^{\nu - \rho_Q}\; \int_{N_Q} \chi(a^{-1}n^{-1} a ) \inp{\gs(m)^{-1} \gf_{\bar Q,\gs, -\nu}(n)}{\eta} \; dn.
\end{eqnarray*}

The integrand of the final integral may be estimated by 
$\ge(n):= C e^{(\nu - \rho_{\bar Q})H_{\bar Q} (n)}$, with $C>0$ uniform in $n \in N_Q$ and $a \in A_Q.$ Since the mentioned function
$\ge$ is absolutely integrable over $N_Q,$  it follows by dominated convergence
that 
\begin{eqnarray*}
\lefteqn{\lim_{a {\buildrel \bar Q \over \to} \infty} a^{-(\nu -\rho_Q)} \Wh(Q, \psi_T ,  \nu)(ma) =  \int_{N_Q} 
\inp{\gs(m)^{-1} \gf_{\bar Q,\gs, -\nu}(n)}{\eta} \; dn}\\
&  =& 
\inp{\gs(m)^{-1} A(Q, \bar Q, \gs, -\nu) \gf}{\eta} \; dn
= \psi_{[ A(Q, \bar Q, \gs, -\nu) \gf  \otimes \eta]}(m).
\end{eqnarray*}
Here the limit means that $a^\ga\to \infty$ for each $\bar Q$-root $\ga.$ 
On the other hand, it follows from (\ref{e: characterization C functions})
that for $\nu \in \faQd(\ge)$ with $\Re\nu$ $\bar Q$-dominant the limit is given 
by $C_{Q|Q}(1: \nu)\psi_T.$ This establishes (4.2) for all $\nu$ in a non-empty open subset of $\faQd(\ge).$ The validity of (4.2) for all generic $\nu \in \faQd(\ge)$ follows by application of analytic continuation.
\qed

 \begin{lemma}
\label{l: standard int ops on Whit int}
Let $P, P' \in \cP$ have the same split  component and suppose $\gs \in \dMPds.$ If 
$T$ is an element of $ C^\infty(\tau: K/K_P : \gs_P) \otimes \ccH^{-\infty}_{P, \chi_P},$ then
for generic $\nu \in \faPd(\ge),$ 
\begin{equation}
\label{e: comp Wh and A}
\Wh(P, \psi_{T} , \nu) =  \Wh(P', \psi_{[A(\bar P, \bar P', \gs, -\nu)^{-1}  
\otimes B(\bar P', \bar P, \gs, \bar \nu) ]T }, \nu).
\end{equation} 
\end{lemma}
\proof
We may assume that $T = \gf \otimes\eta.$ 
It follows from Lemma \ref{l: A B and j} that 
$$
A(\bar P', \bar P, \gs, \nu) j(\bar P,\gs, \nu)\eta = j(\bar P', \gs, \nu) B(\bar P', \bar P, \gs, \nu)\eta.
$$ 
Using this and the identity $A(\bar P, \bar P', \gs, \bar \nu)^* = A(\bar P', \bar P, \gs, -\nu),$ 
we infer that
\begin{eqnarray*}
\Wh(P, \psi_T, \nu)(x)  & = & \inp{ \pi_{\bar P, \gs, -\nu}(x)^{-1} \gf}{j(\bar P, \gs, \bar \nu)\eta}
\\
&=&  \inp{  \pi_{\bar P, \gs, -\nu}(x)^{-1} \gf}
 { A(\bar P', \bar P, \gs, \bar \nu)^{-1}j(\bar P', \gs, \bar \nu)B(\bar P', \bar P, \gs, \bar \nu)\eta} \\
 &=& 
  \inp{ 
 \pi_{\bar P', \gs, -\nu}(x)^{-1} A(\bar P, \bar P', \gs, -\nu)^{-1} \gf}
 {j(\bar P', \gs, \bar \nu)
 B(\bar P', \bar P, \gs, \bar \nu)\eta}.
\end{eqnarray*}
The required identity now follows.
\qed

\begin{lemma}
\label{l: C func and A B}
Let $P \in \cP$ and  $Q\in \cPst$ have the same split component and suppose $\gs \in \dMPds.$   If 
$T$ is an element of $ C^\infty(\tau: K/K_P : \gs_P) \otimes \ccH^{-\infty}_{P, \chi_P},$ then
$$
C_{Q|P}(1, \nu)  \psi_T  = \psi_{[A(Q, \bar P, \gs, -\nu)  \otimes 
B(\bar Q, \bar P, \gs, \bar \nu)]T} 
$$ 
In particular, $\nu \mapsto C_{Q|P}(1, \nu)$ extends to a meromorphic
 $\Hom(\cA_{2,P}, \cA_{2,Q})$-valued function on $\faPdc.$ 
\end{lemma}
\proof 
For any $P' \in \cP$ with the same split component as $P$ we obtain,
by taking the  constant terms of the Whittaker integrals in (\ref{e: comp Wh and A}) along $Q$ and comparing coefficients of exponents,
$$
C_{Q|P}(1, \nu)  \psi_T = C_{Q|P'}(1, \nu)  \psi_{[A(\bar P ,\bar P', \gs, -\nu)^{-1} \otimes 
B(\bar P', \bar P, \gs, \bar \nu)]T}
$$ 
In particular, substituting $P' = Q$ and using Lemma \ref{l: C Q Q} we obtain
\begin{eqnarray*}
C_{Q|P}(1, \nu)  \psi_T & = & \psi_{[A(Q, \bar Q, \gs, -\nu) A(\bar P, \bar Q, \gs, -\nu)^{-1} \otimes 
B(\bar Q, \bar P, \gs, \bar \nu)]T}\\
&=& \psi_{[A(Q, \bar P, \gs, -\nu)  \otimes 
B(\bar Q, \bar P, \gs, \bar \nu)]T}.
\end{eqnarray*}
\qed

\begin{cor} Let  $Q\in \cPst$ and suppose $\gs \in \dMQds.$   If 
$T$ is an element of $ C^\infty(\tau: K/K_Q : \gs_Q) \otimes \ccH^{-\infty}_{Q, \chi_Q},$ then 
$$
C_{Q|\bar Q}(1, \nu) \psi_T = \psi_{[I \otimes B(\bar Q, Q, \gs, \bar \nu)] T},
\qquad (\nu\in \faQdc).
$$ 
\end{cor}
 \proof
This follows from the previous lemma by taking $P = \bar Q,$
 \qed

The next step is to obtain a formula for $C_{Q\mid P}(s, \nu), $ for $s \in N_K(\fa).$ Let $Q\in \cPst$, $P \in \cP$ and let $s \in N_K(\fa)$ such that $s P s^{-1} = Q.$
Then the right regular action of $s$  defines an intertwining operator
$
R_{s}: C^\infty(G/P: \gs:\nu) \to C^\infty(G/Q: s \gs: s \nu).
$
According to Cor. \ref{c: diagram composition Rw and j} applied with $P, s$ in place of $Q,w$
there exists an isometric isomorphism
$$
\cR_{s, P}: \ccH^{-\infty}_{\chi_P} \to  \ccH^{-\infty}_{s\gs, \chi_{s P s^{-1}}}
$$
such that 
$ R_s\after j(P,\gs,\nu) = j(sPs^{-1}, s\gs, s\nu)\after \cR_{s,P} .$ 
Furthermore, $R_s$ induces an isometric isomorphism
 $
 \cR_{s}: C^\infty(\tau: K/K_P: \gs_P) \to C^\infty(\tau: K/K_Q: s\gs_Q).
 $
As in Lemma \ref{l: existence cR Q} it follows that there exists a unique  isometric isomorphism 
$\underline \cR_s: \cA_{2,P}\to \cA_{2, sPs^{-1}},$ such that for every $\gs \in \dMPds$ 
and every  $T \in C^\infty(\tau: K/K_P: \gs_P) \otimes \ccH^{-\infty}_{\gs, \chi_P}$
$$
\underline \cR_s( \psi_T )=  \psi_{(\cR_s \otimes \cR_{s,P})T} \in  \cA_{2, sPs^{-1}, s\gs}.
$$

\begin{lemma}
\label{l: action of s on Whittaker int}
Let $P \in \cP$ and $s\in N_K(\fa)$. Then, for $x\in G,$
$$ 
 \Wh(P, \psi, \nu)(x) = \Wh(s P s^{-1},\underline \cR_s \psi ,  s \nu), 
 \qquad (\nu \in \faPdc).
 $$ 
 \end{lemma} 
 
 \proof
 This is derived from the intertwining property of $R_{s}$ as follows.
 Let $\gs \in \dMPds$ be a type appearing in $\cA_{2, P}.$ By linearity it suffices
to check  the identity for $\psi \in \cA_{2, P, \gs}.$ Then $\psi = \psi_T$ 
with $T \in C^\infty(\tau: K/K_P: \gs_P) \otimes \ccH^{-\infty}_{\gs, \chi_P}.$
By linearity we may assume that $T = \gf \otimes \eta.$ 
Then
\begin{eqnarray*}
\Wh(P, \psi, \nu)(x) &= & \inp{\pi_{\bar P, \gs, -\nu}(x)^{-1} \gf}{j(\bar P, \gs, \bar\nu)\eta}\\
&=& 
\inp{R_{s}\; \pi_{\bar P, \gs, -\nu}(x)^{-1} \gf}{R_{s}\; j(\bar P, \gs, \nu)\eta}
\\
&=& \inp{\pi_{s\bar Ps^{-1}, s\gs,  -s \nu}(x)^{-1} 
\cR_{s} \gf}{j(s\bar Ps^{-1}, s \gs, s\bar \nu) \cR_{s, P}\eta}
\\
&=& \Wh (sPs^{-1}, \psi_{ \cR_{ s} \gf \otimes \cR_{s, P}\eta} ,  s \nu)( x)
\\
&=& \Wh(sPs^{-1}, \underline \cR_{s} \psi_{T}, s \nu  )(x).
\end{eqnarray*}
\vspace{-15pt}
\qed
 
 \begin{cor}
 \label{c: C func s on P}
 Let $s \in N_K(\fa)$ be such that $s(\fa_P) = \faQ.$  Then 
 $$ 
C_{Q|P} (s, \nu) \psi_T =  C_{Q|sP s^{-1}}(1 , s\nu)  \underline \cR_{s} \psi_{T}
$$
for generic $\nu\in i \faPd.$
 \end{cor}
 
\begin{cor}
 \label{c: Cstar C and eta}
 Let $s \in N_K(\fa)$ be  such that $sP s^{-1}= Q,$ and let $\gs \in \dMPds.$ Then for $\psi \in \cA_{P,2,\gs},$
 \begin{equation}
 \label{e: Cstar C and eta}
 C_{Q|P} (s,-\bar \nu)^* C_{Q|P}(s,  \nu) \psi = \eta(P, \bar P, \gs, -\nu) \psi.
 \end{equation}
 \end{cor} 

\proof 
It suffices to prove this for $\psi = \psi_T$ with $T = f \otimes v.$ Combining  Cor. \ref{c: C func s on P} and
Lemma \ref{l: C Q Q} we find
$$
C_{Q|P}(s,  \nu) \psi_T = C_{Q\mid Q}(1, s\nu) \psi_T = \psi_{(A(Q,\bar Q, s\gs, -s\nu)\otimes I)T}.
$$ 
Hence,
$$
C_{Q|P}(s,  -\bar \nu)^*C_{Q|P}(s,  \nu)  \psi_T = \psi_{\eta(Q,\bar Q, s\gs, - s\nu) \otimes I)T}
= \eta(s^{-1} Q s, s^{-1}\bar Q s , \gs, - \nu) \psi_T
$$
$$
= \eta(P, \bar P, \gs, -\nu) \psi_T.
$$
\qed

Let $P \in \cP$ be a cuspidal parabolic subgroup. We denote by $[P]_\ist$ the set of $Q \in \cPst$ that are associated
with $P$.  

\begin{defi}
The following relations MSC(P) will be called Maass-Selberg relations for the 
$C$-functions of the Whittaker integral $\Wh(P)$
\medno
MSC(P): 
for all $Q_1, Q_2 \in [P]_\ist$ and all $s_j \in W(\fa_{Q_j}, \fa_P),$ $(j=1,2)$,
\begin{equation}
\|C_{Q_1|P}(s_1, \nu)\psi\| = \|C_{Q_2|P}(s_2, \nu)\psi\|, \qquad (\psi \in \cA_{2, P}),
\end{equation} 
for generic $\nu \in i \faPd.$ 
\end{defi}
\begin{lemma}
\label{l: MSrels C P}
The Maass-Selberg relations {\rm MSC(P)} for $P$ as stated above are equivalent to the following.
\medno
{\rm MSC(P)':}
for all $Q \in [P]_\ist$, all $s\in W(\fa_{Q}\mid \fa_P)$, all $\gs \in \dMPds$
and all generic $\nu \in i\faPd,$ 
$$
C_{Q|P}( s , \nu)^*  C_{Q|P}(s,  \nu)  = \eta(P, \bar P, \gs, -\nu)\;\;{\rm on\;\;} \cA_{2,P, \gs} .
$$
\end{lemma}

\proof
Assume that MSC(P) hold.
There exists a unique $Q_1 \in \cPst$ which is $W$-conjugate to $P$. It is given by $Q_1 = s_1P s_1^{-1},$ 
where $s_1 = v_P.$ In particular, $Q_1 \in [P]_\ist, s_1\in W(\faQ\mid \faP)$ and it follows from Corollary \ref{c: Cstar C and eta} that (\ref{e: Cstar C and eta}) is valid with
$Q_1, s_1$ in place of $Q, s$, for all $\psi \in \cA_{P,2,\gs}.$ 
For such $\psi$ we find by application of MSC(P) that 
\begin{eqnarray*}
\inp{C_{Q|P}(s , \nu)^*  C_{Q|P}(s,  \nu) \psi}{\psi} & = & \inp{C_{Q_1|P}(s_1, \nu)^*  C_{Q_1| P}(s_1,  \nu) \psi}{\psi}\\
&= &\eta(P, \bar P, \gs, -\nu)\inp{\psi}{\psi}.
\end{eqnarray*} 
Now $C_{Q_|P}(s , \nu)^*  C_{Q|P}(s,  \nu)
$ is Hermitian, and the only eigenvalue of its restriction to $\cA_{P,2,\gs}$ can be $\eta(P, \bar P, \gs, - \nu).$
It follows that the latter Hermitian map is the scalar $\eta(P, \bar P, \gs, - \nu)$ on $\cA_{P,2, \gs}.$ Therefore, MSC(P)' holds.

The converse implication is straightforward.
\qed

We will now compare the Maass--Selberg relations formulated above with those for the $B$-matrix.
Recall from the text following (\ref{e: definition eta}) that 
$\eta(P,Q, \gs, \nu) = \eta(Q, P, \gs, \nu).$

\begin{prop} 
\label{p: comp MS rels B and C}
Let $P\in \cP$ and $\gs \in \dMPds.$ Then the following assertions are equivalent, for each $Q \in [P]_\ist$ and all
$s \in W(\faQ\mid \faP).$
\begin{enumerate}
\itema $C_{Q\mid P}(s , \nu)^*  C_{Q\mid P}(s,  \nu)  = \eta(P, \bar P, \gs, -\nu)$ on $\cA_{P,2,\gs}$ for generic $\nu \in i\faPd;$
\itemb  $B(s^{-1}\bar Q s, \bar P, \gs, - \nu)^* B( s^{-1} \bar Q s , \bar P, \gs, -  \nu) = \eta(s^{-1}  Q s , P, \gs, -  \nu)$ on $\ccH^{-\infty}_{\gs, \chi_{\bar P}},$
\\for generic $\nu \in i\faPd$.
\end{enumerate}
\end{prop}

\proof
Let $s \in W(\faQ\mid \faP).$ For $T \in C^{\infty}(\tau: K/K_P : \gs_P) \otimes \ccH^{-\infty}_{\gs, \chi_P}$ 
we have, by Corollary \ref{c: C func s on P}  and Lemma \ref{l: C func and A B} that
\begin{eqnarray*}
C_{Q|P}(s: \nu)\psi_T &=& C_{Q  |sP s^{-1}}  (1, s\nu) \psi_{-}   [  (\cR_s \otimes \cR_{s,P}) T ] \\
&=&
\psi_{-} [(A(Q,s \bar P s^{-1}, s\gs, - s\nu ) \otimes B(\bar Q, s\bar P s^{-1}, s\gs,-s \nu) R_s T ].
\end{eqnarray*}
Since $T \mapsto \psi_T$ is unitary from $C^\infty(\tau: K/K_P: \gs_P) \otimes \ccH^{-\infty}_{\gs, \chi_P}$  onto 
$\cA_{P,2, \gs},$  and unitary from $C^\infty(\tau: K/K_P: s\gs_P) \otimes \ccH^{-\infty}_{s\gs, \chi_P}$
onto $\cA_{2,P,  s\gs},$  it follows from the above that
$$
C_{Q|P}(s: \nu)^* \psi_S = \psi_{-}[ R_{s^{-1}} (A(Q,s \bar P s^{-1}, s\gs, - s\nu )^*\ \otimes 
B(\bar Q, s\bar P s^{-1}, s\gs,s \nu)^*)  S].
$$
for $S \in C^\infty(\tau: K/K_P: s \gs_P) \otimes \ccH^{-\infty}_{s \gs, \chi_P}.$ Combining the above, and using
that  $A(Q,s \bar P s^{-1}, s\gs, - s\nu )^* A(Q,s \bar P s^{-1}, s\gs, - s\nu ) = \eta(Q, s\bar P s^{-1}, s\gs, - s\nu),$
we infer
\begin{equation}
\label{e: formula for CstarC}
C_{Q|P}(s: \nu)^*C_{Q|P}(s: \nu) \psi_T = \psi_-[\eta(Q, s\bar P s^{-1}, s\gs, - s\nu) \otimes b(\nu))  T]
\end{equation} 
where 
\begin{eqnarray*}
b(\nu) &= &R_{s^{-1}} B(\bar Q, s\bar P s^{-1}, s\gs,- s \nu)^* B(\bar Q, s\bar P s^{-1}, s\gs,- s\nu)R_s \\
&=& B(s^{-1}\bar Q s, \bar P, \gs,-  \nu)^* B(s^{-1} \bar Q s, \bar P , \gs,-\nu)
\end{eqnarray*}
and where we have used the notation $\psi_{-} T :=\psi_T.$
Suppose now that (a) is valid. Then it follows from (\ref{e: formula for CstarC}) that
$b(\nu)$ must be a multiple of the identity map by the 
 (non-negative real) factor
\begin{eqnarray*}
\lefteqn{\!\!\!\!\!\!\!\!\!\!\!\!\!\!\!\!\!\!\!\!\!\!\!\!\!\!\eta(P, \bar P, \gs, -\nu) \eta(Q, s\bar P s^{-1}, s\gs, - s\nu)^{-1}}\\
 &=& \eta(P, \bar P, \gs, -\nu) \eta(s^{-1} Q s, \bar P, \gs, -\nu)^{-1} \\
&=& \eta(s^{-1} Q s,  P , \gs, -\nu),
\end{eqnarray*}
and (b) follows.

Conversely, suppose that (b) is valid.
Then it follows that $b(\nu)$ is the scalar $\eta(s^{-1} Qs, P, \gs, -\nu).$  From  (\ref{e: formula for CstarC})
we  now see that $C_{Q|P}(s,\nu)^*C_{Q|P}(s,\nu)$ is the scalar
\begin{eqnarray*}
\lefteqn{\!\!\!\!\!\!\!\!\!\!\!\!\!\!\!\!\!\!\!\!\!\!\!\!\!\!\!   \eta(Q, s \bar P s^{-1}, s\gs, - s\nu) \eta(s^{-1}Q s, P, \gs, -\nu) }\\
&=& \eta(s^{-1}Q s, \bar P , \gs, - \nu) \eta(s^{-1} Q s,  P, \gs, -\nu)\\
&=& \eta(\bar P, P, \gs, -\nu).
\end{eqnarray*}
\qed

\begin{defi} 
The following relations MSB(P) will be called Maass-Selberg relations
for the $B$-matrices associated with $P$:
\medno
MSB(P): for all $Q\in \cP$ with $\faQ = \faP$ and all $\gs \in \dMPds$
$$
B(\bar Q, \bar P, \gs, \nu)^* B(\bar Q , \bar P, \gs, \nu) = \eta(\bar Q, \bar P, \gs, \nu)
$$ 
for generic $\nu \in i\faPd.$
\end{defi}
\medbreak

From Proposition \ref{p: comp MS rels B and C} we see, for each $P \in \cP,$ that the validity of the relations MSB(P) implies the validity of the relations MSC(P)'.
The converse is not clear a priori, except in the basic setting, where
$G$ has compact center and $P$ is maximal.  This will be addressed in the next
section.
 
\section{Maass--Selberg relations in the basic setting}
\label{s: MS basic setting}
 
We consider the basic setting in which $G$ has compact center, and $P \in \cP$ is a maximal parabolic subgroup.  In this case there is precisely one $Q \in \cP$ which is adjacent to $P,$ namely $\bar P.$ From Proposition \ref{p: comp MS rels B and C} it follows that
the Maass-Selberg relations MSB(P) for the $B$-matrix imply the relations 
MSC(P)'  for the $C$-functions,
but the converse is not obvious. In the present section we will show that
the converse is obvious for the basic setting.

 \begin{lemma}
 \label{l: classes of max parab}
 The following assertions are equivalent.
 \begin{enumerate}
 \itema $|W(\faP)| =1;$
 \itemb $\bar P$ is not $W(\fa)$-conjugate to $P;$
 \itemc $[P]_\ist$ has two elements.
 \end{enumerate}
 \end{lemma}
 
 \proof 
 First of all, by using the action of $W(\fa)$ we see that we may as well assume 
 that $P$ is  standard.
 
 Suppose (a). Then $\fa_P^+$ and $-\fa_P^+$ are not $W(\faP)$-conjugate 
 hence not $W(\fa)$-conjugate and (b) follows.
 
 Suppose (b), then $\bar P$ is $W(\fa)$-conjugate to precisely one $Q \in \cPst$  
 which we know cannot be $P.$ It follows that $[P]_\ist$ has at least $2$ elements.
 If $R$ were a third element of $[P]_\ist$ then there would be an element $s\in W(\fa)$ 
 such that $s(\faP) = \fa_R.$ Then either $s(\faPp)= \fa_R^+$ or $s(-\faPp) = \fa_R^+$. In the first case it would follow that $\faQp$ and 
$\fa_R^+$ are $W(\fa)$-conjugate. But then, since $P,R\in \cPst,$ it would follow that  $R = P,$  
contradiction. In the second case, it would follow that $Q,$ $\bar P$ and $R$ are
conjugate under $W(\fa)$ hence $R = Q,$ contradiction.
 
 Finally, suppose (c). Then there is a parabolic subgroup $Q\in \cPst\setminus\{P\}.$ 
 such that $\faQ$ is $W(\fa)$-conjugate to $\faP.$ Hence, 
  $\faQp$ is conjugate to either $\faPp$ or $-\faPp.$ The first cannot be true
 since then $P = Q.$ Therefore, $\faPp$ is not conjugate to $-\faPp.$ 
 From this it follows that $- \faPp$ cannot be conjugate to $\faPp.$ 
 It follows that $\faPp$ and $-\faPp$ are not conjugate under $W(\faP).$ 
 Hence, (a) follows. 
 \qed
 
 \begin{rem}
 It follows from the proof that in any case  $[P]_\ist$  has at most two elements.
 \end{rem}
 
\begin{prop}
\label{p: MS for adjacent maximal}
Let $P \in \cP$ be a maximal parabolic subgroup of $G.$
\begin{enumerate}
\itema
If $|W(\faP)| =1 $ then $[P]_\ist$ consists of two distinct elements,  $Q_1, Q_2 \in \cPst.$ The constant
terms of  $\Wh(P, \psi, \nu)$ along  $Q_j,$ for $j=1,2,$ are of the form
$$
\Wh(P, \psi, \nu)_{Q_j} (ma) = a^{s_j\nu}C_{Q_j|P}(s_j , \nu)\psi(m),
$$ 
with $W(\fa_{Q_j}\mid \fa_P) = \{s_j\}.$ Furthermore, 
$s_1|_{\faP} = -s_2|_{\faP}.$ In this case the Maass-Selberg relations 
MSC(P) are equivalent to
\begin{equation}
\label{e: MS BC a}
\|C_{Q_1|P}(s_1 ,\nu)\psi\|^2 
=\|C_{Q_2|P}(s_2 ,\nu)\psi\|^2.
\end{equation}
for all $\psi\in \cA_{2,P}$ and a dense set of $\nu \in i\faPd.$ 
\itemb
If $|W(\faP)| = 2,$ then $[P]_\ist$ consists of  a single element $Q$ in $\cPst$ 
and $|W(\faQ\mid\faP)| = 2.$  The constant
term of  $\Wh(P, \psi, \nu)$ along  $Q$ is of the form
$$
\Wh(P, \psi, \nu)_{Q} (ma) 
 = a^{s\nu}C_{Q|P}(s , \nu)\psi(m) + a^{-s\nu}C_{Q|P}(-s , \nu)\psi(m).
$$ 
where $W(\faQ|\faP)= \{s, -s\}.$ In this case the Maass-Selberg  relations 
MSC(P) are equivalent to 
\begin{equation}
\label{e: MS BC b} 
\|C_{Q|P}(s : \nu)\psi\|^2 = \|C_{Q|P}(-s : \nu)\psi\|^2.
\end{equation}
for all $\psi \in \cA_{2, R}$ and all regular values of $\nu  \in  i\faPd.$ 
\end{enumerate}
\end{prop}

\proof
(a) $P$ is conjugate to a  standard  parabolic subgroup $Q_1.$ Clearly $Q_1 \in [P]_\ist.$ The latter set 
has two distinct elements, hence equals $\{Q_1, Q_2\}$ with $Q_2$ a second maximal parabolic subgroup
associated with $P.$ Since $W(\fa_P)$ has a single element, there exist for each $j=1,2$ a single element
$s_j \in W(\fa_{Q_j}, \fa_P).$  Since $Q_1,Q_2$ are standard and not equal, 
they cannot be $W(\fa)$-conjugate. The element $s = s_2 s_1^{-1}$ of 
$W(\fa)$ maps $\fa_{Q_1}$ to $\fa_{Q_2}$ but not $\fa_{Q_1}^+$ to $\fa_{Q_2}^+$. Therefore, there exist a point $X \in \faPp$ such that
$s_2^{-1} s_1 (X) \notin \faPp.$ Since $\faPp$ is one dimensonal
and since $s_2^{-1} s_1 $ is length preserving, it follows that $s_2^{-1} s_1(X) = -X.$ Hence, $s_1|_\faP = - s_2|_\faP.$ 

By application of Proposition \ref{p: comp MS rels B and C} we see that in this case the Maass-Selberg relations associated with $P$ are completely described by (\ref{e: MS BC a}).

We turn to case (b). By Lemma \ref{l: classes of max parab} there is a unique $Q\in \cPst$ such that $[P]_\ist =\{Q\}.$  Hence 
$\faP$ and $\faQ$ are $W(\fa)$-conjugate and $|W(\faQ|\faP)| = |W(\faP)| =2.$
The constant term of $Wh(P, \psi, \nu)$ along $Q$ is described by 
$$
Wh(P, \psi, \nu)_Q (ma) = \sum_{s \in W(\fa_Q, \fa_P)} a^{s\nu} C_{Q|P}(s:\nu)\psi(m).
$$
Take $s\in W(\faQ| \faP),$ then $W(\faQ| \faP) = \{\pm s\}$ and it follows that the description of the Maass-Selberg 
relations is complete.
\qed

\begin{lemma}
\label{l: basic setting MSC gives MSB}
Let $G$ have compact center and let
$P\in \cP$ be a maximal parabolic subgroup of $G.$ 
Assume the relations {\rm MSC(P)} are valid. Then for every $\gs \in \dMPds$ the Maass-Selberg relations {\rm MSB(P)} are valid, i.e.
\begin{equation}
\label{e: MS for B basic setting}
B(P, \bar P, \gs, -\nu)^*B(P, \bar P, \gs, -\nu) = \eta(\bar P, P, \gs, -\nu)
\end{equation}
for generic $\nu \in i\faPd.$ 
\end{lemma}

\proof
First assume that we are in case (a): $ |W(\faP)| = 1.$ Then $[P]_\ist = \{Q_1, Q_2\}$ and $s_1, s_2$ 
are as in Proposition \ref{p: MS for adjacent maximal}. Fix $\gs \in \dMPds$. Then assertion (a) of Prop. \ref{p: comp MS rels B and C}  is valid for each choice 
$(Q,s) \in \{ (Q_1, s_1), (Q_2, s_2)\}.$ It follows that assertion (b) is valid for each choice.
For one of the choices one has $s^{-1}Qs = \bar P$ hence $s^{-1} \bar Q s = P$ 
and the validity of assertion (b) now implies that (\ref{e: MS for B basic setting}).

Next assume that we are in case (b): $|W(\faP)| = 2.$ 
Then $[P]_\ist = [Q]_\ist$ and there exists a $s\in W(\faQ, \faP)$ which 
maps $\faPp$ to $- \faQp$. Then $s P s^{-1} = \bar Q.$ The condition (a) is fulfilled 
hence also (b). We find (\ref{e: MS for B basic setting}).
\qed

Thus, in order to complete the proof of the MS relations for the $B$-matrix, it suffices to give a proof of the assertions of MSC(P) for the basic
setting, as listed in Proposition \ref{p: MS for adjacent maximal}. We will do this, following a method of Harish-Chandra  \cite{HCwhit} in Sections 
\ref{s: radial part Casimir} - \ref{s: result HC}.
%%%%%%
 \section{The radial part of the Casimir operator}
 \label{s: radial part Casimir}
 By the Iwasawa decomposition $G = K A N_0,$ the multiplication
 map $m: K \times A \times N_0 \to G$ is a diffeomorphism.
 Accordingly, we may define a topological linear isomorphism 
 $\Tup: C^\infty(A, V_\tau ) \to C^\infty(\tau:G/N_0:\chi)$ 
 by 
 $$
 \Tup f (kan) = \chi(n)^{-1} \tau(k) f(a),
 $$
 for $f \in C^\infty(A)$ and $(k,a,n) \in K\times A \times N_0.$
 The inverse  of this isomorphism  is given by the restriction map $\Tdown: f \mapsto f|_A.$ 
 
 For an element $u \in U(\fg)^{N_0} $ we consider the differential operator
 $R_u$ on $C^\infty(G, \Vtau)$ given by $R_u(f)(x): = f(x\rdiff u).$
 Then $R_u$  restricts to a differential operator $r_u$ on $C^\infty(\tau:G/N_0:\chi).$ The radial part of the latter, denoted
 $\Pi(u),$ is defined by 
 $$ 
 \Pi(u) =\Tdown \after r_u \after \Tup.
 $$
 We will determine the radial part of the Casimir element $\Omega \in \fZ$ associated
 with the invariant symmetric bilinear form $B$ on $\fg,$ see (\ref{e: defi B}).
For each $\ga \in \gS^+$ we fix a basis $X_{\ga, i}, 1 \leq i \leq m_\ga,$ 
which is orthogonal with respect to the positive definite 
inner product $X,Y\mapsto -B(X, \Cartan Y).$ Furthermore, we put
$X_{-\ga, i} = -\Cartan X_{\ga, i}.$ Let $H_\ga \in \fa$ be defined
by $H_\ga \perp \ker \ga$ and $\ga(H_\ga) = 1.$ Then $\ga = B(H_\ga,\dotvar).$ 
It is readily  seen that 
$$
[X_{\ga, i} , X_{-\ga, i}] = H_\ga,\; (1 \leq i \leq m_\ga).
$$ 
The Casimir operators of $\fm$ and $\fa$, defined relative
to the restrictions of $B$ to these Lie algebra's, are denoted
by $\Omega_\fm$ and $\Omega_\fa.$ It is now well known that 
$$
\Omega = \Omega_\fm + \Omega_\fa + \sum_{\ga, i} (X_{\ga, i} X_{-\ga, i} + X_{-\ga, i} X_{\ga, i});
$$ 
here the summation ranges over $\ga \in \gS^+$ and $1 \leq i \leq m_\ga.$
The radial part of an operator $X \in U(\fg)^{N_0}$ may be calculated by
from a decomposition of the form
$$
X = \sum_j  f_j(a)  z_j ^{a^{-1}}u_j v_j,\qquad(a \in A),
$$ 
with $f_j \in C^\infty(A),$ $Z_j \in U(\fk) ,  u_j \in U(\fa), v_j \in U(\fn_0).$ 
Here the superscript $a^{-1}$ indicates that the image under $\Ad(a)^{-1}$
is taken.
Given a decomposition as above the radial component may be expressed by 
$$
[\Pi(X) \gf](a)= \sum_j f_j(a) \;\tau_*(z_j) \; \gf(a\rdiff u_j) \; \chi_*(v_j^\vee).
$$
Put $Z_{\ga, i} := X_{\ga, i} - X_{-\ga, i}.$ Then $Z_{\ga.i} \in \fk.$ 
Furthermore, for each $\ga \in \gS^+ $ and $1 \leq i \leq m_\ga,$ 
$$ 
 X_{-\ga, i} = a^{-\ga} Z_{\ga, i} ^{a^{-1}} +  a^{-2\ga} X_{\ga, i}, 
 $$ 
 It follows from this that 
 $$
 X_{-\ga, i} X_{\ga,i } = a^{-\ga} Z_{\ga, i}^{a^{-1}} X_{\ga, i} +  a^{-2\ga} X_{\ga, i}^2.
 $$
 On the other  hand, 
 $$ 
 X_{\ga, i} X_{-\ga, i} = X_{-\ga, i} X_{\ga, i} + H_\ga.
 $$ 
Hence,
$$ 
X_{\ga, i} X_{-\ga, i} + X_{-\ga, i} X_{\ga, i} =  2 X_{-\ga, i} X_{\ga, i} + H_\ga
=
2 a^{-\ga} Z_{\ga, i}^{a^{-1}} X_{\ga, i} + H_\ga + 2 a^{-2\ga} X_{\ga, i}^2
$$

\begin{lemma}
$$
\Pi(\Omega) = 
 \tau_*(\Omega_\fm ) + \Omega_\fa + \sum_{\ga} m_\ga H_\ga + 
 \sum_{\ga, i} - [ 2 a^{-\ga} \tau(Z_{\ga, i})\chi_*(X_{\ga, i}) + 2 a^{-2\ga}\chi_*( X_{\ga, i})^2]
 $$
\end{lemma}

For two functions $f, g \in C^\infty(\tau: G/N_0: \chi)$ we define the function
$[f,g]: G \to \C$ by 
\begin{equation} \label{e: brackets with Omega}
[f,g](x)  = \inp{f(x\rdiff \Omega)}{g(x)} - \inp{f(x)}{g(x\rdiff \Omega)}.
\end{equation} 
We define $\omega \in U(\fa)$ by 
$$ 
\omega = \Omega_\fa + \sum_{\ga \in \gS^+} m_\ga H_\ga.
$$ 
\begin{lemma}
\label{l: brack f g and omega}
$$
[f, g](a) = \inp{f(a;\omega)}{g(a)} - \inp{f(a)}{g(a\rdiff \omega)}.
$$
\end{lemma}

\proof
Since $\tau$ and $\chi$ are unitary, the operators $\tau_*(Y)$ for $Y \in \fk$ 
are anti-Hermitian, while $\chi_*(X) \in i \R$ for $X \in \fn_0.$ It follows from 
this that the operators
$$ 
\tau_*(\Omega), \;\;\tau(Z_{\ga, i}), \;\;\tau(Z_{\ga, i})\chi_*(X_{\ga, i})
$$ 
from $\End (V_\tau)$ are Hermitian, while $\chi_*(X_{\ga,i})^2 \in \R.$ 
It follows that 
$$ 
S(a): = \tau_*(\Omega_\fm) +  \sum_{\ga, i} - [ 2 a^{-\ga} \tau(Z_{\ga, i})\chi_*(X_{\ga, i}) + 2 a^{-2\ga}\chi_*( X_{\ga, i})^2]
$$ 
is Hermitian for all $a.$ Now $\Omega = \omega + S(a).$ 
Hence 
\begin{eqnarray*}
[f, g](a) &= &\inp{ f(a\rdiff \omega)) + S(a) f(a)}{g} - \inp{f}{g(a\rdiff \omega)+ S(a) g(a)} \\
&=&\inp{f(a\rdiff \omega)}{g(a)} - \inp{f}{g(a\rdiff \omega)} + \inp{  S(a) f(a)}{g} - \inp{f}{S(a) g(a)}\\
&=&\inp{f(a\rdiff \omega)}{g(a)} - \inp{f}{g(a\rdiff \omega)}.
\end{eqnarray*}
\qed
\begin{lemma}
\label{l: omega Omega}
$
a^{\rho} \after \omega \after a^{-\rho} =  \Omega_\fa - \inp{\rho}{\rho}.
$
\end{lemma}

\proof
Let $H_j$ be an orthonormal basis for $\fa.$ The dual inner product on $\fa^*$ makes
$B: \fa \to \fa^*$ orthogonal; in particular, $B(H_j),$ for $1 \leq j \leq  \ell$, is an orthonormal 
basis for $\fad.$ Accordingly, if $\gl, \mu\in \fad$ then $\inp{\gl}{B(X_j)} = \gl(H_j)$ 
and it follows that 
$$ 
\inp{\gl}{\mu} = \sum_{1 \leq \ell} \gl(H_j)\mu(H_j).
$$
The Casimir operator of $\fa$ is given by $\Omega_\fa = \sum_{j=1}^\ell H_j^2.$ Moreover $u\mapsto a^{\rho} \after  u \after a^{-\rho}$ 
equals the algebra automorphism  $T = T_{-\rho}$ of $U(\fa),$ determined by 
$T(H) =H - \rho(H).$  From this it follows that
$$ 
T \Omega_\fa = \sum_j [H_j - \rho(H_j)]^2 = \Omega_\fa - \sum_j 2 \rho(H_j) H_j + \inp{\rho}{\rho}.
$$
On the other hand,
$$
\sum_j  2 \rho(H_j) H_j = \sum_{\ga > 0} \sum_j m_\ga \ga(H_j)H_j =
\sum_{\ga > 0, j } m_\ga B(H_\ga, H_j)H_j =
\sum_{\ga > 0} m_\ga H_\ga
$$
from which
$$
T\Omega_\ga = \Omega_\fa - \sum_{\ga} m_\ga H_\ga + \inp{\rho}{\rho}.
$$ 
Hence, 
\begin{eqnarray*}
T(\omega) & = & T(\Omega_\fa + \sum_\ga m_\ga H_\ga)
\\
&=&  \Omega_\fa - \sum_\ga  m_\ga \rho(H_\ga) + \inp{\rho}{\rho}\\
& = & \Omega_\fa - \inp{\rho}{\rho}.
\end{eqnarray*}
\qed

Let $\gD$ denote the collection of simple roots for the positive system $\gS^+.$ 
We define 
$$
{}^\circ \fa:= {}^\circ \fg \cap \fa.
$$
This is the orthocomplement of the intersection
of root hyperplanes $\ker \ga,$ for $\ga \in \gD.$ 
Let $\{H_\ga^0\mid \ga \in \gD\}$ be the $B$-dual of the basis $\gD$ in ${}^\circ \fg.$ 
This subset of ${}^\circ \fa$ is determined  by 
$$
B (H_\ga^0, H_\gb) = \gd_{\ga\gb}, \qquad (\ga, \gb \in \gD).
$$

\begin{lemma}
Suppose $G$ has compact center, then  
\begin{eqnarray*}
\lefteqn{\!\!\!\!\!\!\!\!\!\!\!\!
\inp{f(a\rdiff \Omega_\fa)}{g(a)} - \inp{f(a)}{g(a\rdiff \Omega_\fa)}}\\
&=& 
\sum_{\ga \in \gD} R(H_\ga) ( \inp{R(H_\ga^0) f}{g} - \inp{f}{R(H_\ga^0)g})(a)
\end{eqnarray*}
\end{lemma}
\proof
Since $\{H_\ga^0\}$ is $B$-dual to $\{H_\ga\}$ 
we have
$$
\Omega_\fa = \sum_{\ga\in \gD} H_\ga H_\ga^0.
$$ 
In the following, we will abbreviate $R(H) f$ by $H f.$ 
By substituting this in the left hand side of the above equation, and by application of the Leibniz rule for differentiation, we find that the above equation holds
provided we add to the right hand side the expression 
$$
\cR(f,g, a) = \sum_{\ga\in \gD}   \inp{H_\ga^0 f(a) }{H_\ga g(a)} -\inp{H_\ga f(a)}{H_\ga^0g(a)}.
$$

We will finish the proof by showing that $\cR(f,g, a)= 0.$ Substituting
$H_\ga = \sum_i B(H_\ga, H_i)H_i$ and $H_\ga^0 = \sum_j B(H_\ga^0, H_j) H_j,$
we find that
$$
\cR(f,g) = \sum_{\ga, i, j}[ B(H_\ga^0, H_i)B(H_\ga, H_j) - B(H_\ga, H_i)B(H_\ga^0, H_j)]
 \inp{H_i f}{H_j g}.
 $$
 Fix $i, j.$ 
By duality of $\{H_\ga\}$ and $\{H_\ga^0\},$ $H_j= \sum_\ga B(H_\ga, H_j) H_\ga^0.$
In turn this implies 
$$
B(H_i, H_j) = \sum_{\ga} B(H_\ga^0, H_i)B(H_\ga, H_j).
$$
By a similar reasoning this identity holds with $i$ and $j$ interchanged. Therefore,
$$
\cR(f,g) = \sum_{i,j} [B(H_i, H_j) - B(H_j,H_i)] \inp{H_i f}{ H_j g} = 0.
$$
\vspace{-24pt}

\qed
\medbreak
Given $f, g \in C^\infty(A, \Vtau)$ and $H \in \fa,$ we define
the function $(f, g)_H: A \to \C$ by 
\begin{equation}
\label{e: defi f g sub H}
\!\!\!\!\!\!\!\!\!(f,g)_H(a) = d_0(a)^2 \left[ \inp{f(a;H)}{g(a)} - \inp{f(a)}{g(a;H)}\right], 
\end{equation}
for $a \in A.$
The following lemma is given without proof in  \cite[page 208]{HCwhit}.

\begin{lemma}
\label{l: d zero times fg} 
Let $f,g \in C^\infty(\tau: G/N_0: \chi).$ Then
$$
d_0(a)^2 [f,g](a) = \sum_{\ga \in \gD} H_\ga (f,g)_{H_\ga^0}(a), \qquad (a \in A).
$$
\end{lemma}

\proof
It follows from Lemma \ref{l: brack f g and omega} that
$$
d_0(a)^2 [f, g](a) =d_0(a)^2  \inp{\omega f (a) }{g(a)} - \inp{f(a)} {\omega g (a)}.
 $$
 Using Lemma \ref{l: omega Omega} we now find that
 \begin{eqnarray*}
d_0(a)^2 [f, g](a) 
&= &\inp{\Omega_\fa (d_0 f)}{d_0 g}(a) -  \inp{d_0f}{\Omega_\fa(D_0 g)}(a)
\\
&=& 
 \sum_{\ga \in \gD} H_\ga ( \inp{H_\ga^0 d_0 f}{d_0 g} - \inp{d_0 f}{H_\ga^0 d_0 g)}(a)\\
&=& 
\sum_{\ga \in \gD} H_\ga d_0^2( \inp{H_\ga^0 f + \rho(H_\ga)f}{g}(a) - \inp{f}{H_\ga^0 g+ \rho(H_\ga)g}(a)\\
&=& 
\sum_{\ga \in \gD} H_\ga d_0^2( \inp{H_\ga^0 f}{g}(a) - \inp{f}{H_\ga^0 g}(a))\\
&=& 
\sum_{\ga \in \gD} H_\ga (f,g)_{H_\ga^0} (a).
\end{eqnarray*}
\qed
\section{A result of Harish-Chandra}
\label{s: result HC}
We retain the assumption that $G$ has compact center.

Let $\mu \in \fad$ be defined by $\inp{\mu}{\ga} =1$ for all $\ga \in \gD.$
Equivalently,  $B^{-1} \mu = \sum_{\ga \in \gD} H_\ga^0.$ 
For $t >0$ we define $\fa[t]$ to be the subset of $\fa$ consisting of the points $H \in \fa$ 
such that for all $\ga \in \gD,$ 
$$ 
\inp{H_\ga^0}{ H} \leq t \;\mbox{{\rm and}}\;\;  \mu(H) \geq - t .
$$
Clearly, $\fa[t] = t \fa[1].$ We agree to write $A[t]=\exp \fa[t]$ and $G[t] = K A[t] N_0.$  

\begin{lemma}
The set $\fa[1]$ is a compact neighborhood of $0$ in $\fa.$
\end{lemma}

\proof
Put $\mu_\ga = B H_\ga^0 = \inp{H_\ga^0} {\dotvar};$ then $\mu = \sum_{\ga \in \gD} \mu_\ga.$ 
The set $\fa[1]$ is given by the inequalities $\mu_\ga \leq 1$ and $\mu \geq - 1,$ hence closed,
and a neighborhood of $0.$ It  remains to prove its boundedness. If $H \in \fa[1]$ then 
$$
 \mu_\ga(H) =  \mu(H) - \sum_{\gb \neq \ga} \gb(H) \geq -1 - (|\gD| -1)= -|\gD|.
$$ 
Since $\{\mu_\ga\mid \ga \in \gD\}$ form a set of linear coordinates for $\fa,$ the boundedness follows.
 \qed
Note that the argument in fact demonstrates that $\fa[1]$ is an $\ell$-dimensional simplex,
with $\ell = \dim \fa.$ 
 
In the discussion that follows we will make full use of the Euclidean structure
on $A$ obtained by transfer of structure under the exponential map $\exp: \fa \to A.$
Our notation will be in terms of the multiplicative group in order to emphasize the connection with the structure of the group $G.$  

For  $f,g \in C^\infty(\tau: G/N_0: \chi)$ it is readily checked that the function
$[f,g]$ is left $K$-invariant, and right $N_0$-invariant. 
Hence, for $t>0,$
$$
\int_{G[t]/N_0} [f, g](x)  \; d\bar x	= \int_{A[t]} d_0(a)^2	 [f,g](a)\; da.
$$ 																		
Using Lemma \ref{l: d zero times fg} we find 
\begin{equation}
\label{e: int after Gauss}
\int_{G[t]/N_0}  [f, g](x) d\bar x = 
 \int_{A[t]}  \sum_{\ga \in \gD} H_\ga (f, g)_{H_\ga^0}(a) \, da.
\end{equation}																																	
The integration over $A[t]$ coincides with the Lebesgue integration over $\fa[t]$ 
and the differentiation on $A$ induced by the right regular action coincides
with the usual directional derivative on $\fa.$ This makes Gauss' divergence theorem for the simplex $\fa[t]$ in $\fa$ available, and we obtain:
\begin{equation}
\label{e: first formulation of Gauss}
\int_{G[t]/N_0 }[f, g](x) dx = \sum_{\ga \in \gD} \int_{\partial A[T]}  \inp{\nu}{H_\ga} (f,g)_{H_\ga^0} ds(a).
\end{equation}

Here $\nu$ corresponds to the outward normal vector to the boundary $\partial \fa[t]$
and $ds$ is the $\ell -1$ dimensional Euclidean Lebesgue measure on the boundary.

We will now introduce some structure that is necessary for a proper understanding 
of the integral on the right. The boundary $\partial \fa[1]$ of $\fa[1]$ is the union of
$\ell$ simplices $\fs_\gg,$ for $\gg\in \gD \cup \{\nu\},$ of dimension $\ell -1,$ namely $\fs_\gb$ for $\gb \in \gD$ and a remaining simplex $\fs_\mu.$ 
More precisely, $\fs_\gb$ $(\gb \in \gD)$ is the intersection of $\fa[1]$ with
the hyperplane $\gs_\gb := \{H \in \fa\mid \inp{H}{H_\gb^0} = |\gb|\}$
and $\fs_\mu$ is the intersection of $\fa[1]$ with the hyperplane 
$\gs_\mu:= \{H \in \fa\mid\mu(H) = - 1\}.$
The outward normals are $\nu_{\fs_\gb} = H^0_\gb$ for $\gb \in \gD$ 
and $\nu_{\fs_\mu} = -|\mu|^{-1} B^{-1}\mu.$ We note that 
$$ 
\inp{\nu_{\fs_\gb}}{H_\ga} = \gd_{\ga \gb} \;\;\mbox{\rm and}\;\;
\inp{\nu_{\fs_\mu}}{H_\ga} = - |\mu|^{-1}.
$$
For $t >0$ we define the multiplication operator $M_t: A \to A$ by $M_t (\exp H) = \exp tH$
for $H \in \fa.$ Then $M_t$ maps $A[1]$ onto $A[t]$ and $\partial A[1]$ onto $\partial 
A[t].$  Since $\nu$ is the outward unit normal, we find that $\nu(M_t a) = \nu(a)$ 
for $a \in \partial A[1].$ The pull-back of the surface measure $ds$ by $M_t$ is given 
by $M_t^* ds = t^{\ell -1} ds.$ 
 
Write $\widehat \gD =\gD \cup \{\mu\}$ and put $c_\gg = 1$ for $\gg \in \gD$ and $c_\gg = - |\mu|^{-1}$ for $\gg = \mu.$ For $\gg \in \widehat \gD$ we put $S_\gg = \exp (\fs_\gg).$
Then $\cup_{\gg \in \widehat\gD} S_\gg = \partial A[1]$ and (\ref{e: first formulation of Gauss}) 
takes the following form.
 
\begin{lemma}
\label{l: second form Gauss}
For $t> 0$ and $f,g \in C^\infty(\tau: G/N_0: \chi)$,
$$ 
\int_{G[t]/N_0}  [f, g](x) dx = 
\sum_{\gg \in \widehat\gD} c_\gg  \int_{M_t S_\gg} (f,g)_{H_\ga^0}(a) \, ds(a).
$$ 
\end{lemma}
We will investigate the asymptotic behavior of the given integrals over the hypersurfaces $M_t S_\gg$ as $t \to \infty.$ 
The dominant asymptotic behavior will come from $\gg \in \gD$ and 
neighborhoods of the point $\exp H_\gg^0 \in S_\gg.$ The integral for $\gg = \mu$ will turn out to 
have exponential decay for $t \to \infty$.
The following lemma suggests the relevance of our discussion
for the behavior of the constant terms along maximal parabolic  subgroups.

If $f_1, f_2 \in C^\infty(A, \Vtau)$ and $H \in \fa,$ we define the function 
$\inp{f_1}{f_2}_H: A \to \C$ by 
\begin{equation}
\label{e: defi f g H}
\inp{f_1}{f_2}_H(a) = \inp{f_1(a\rdiff H) }{f_2(a)} - \inp{f_1(a) }{f_2(a\rdiff H)}, \qquad (a\in A).
\end{equation} 
The following useful lemma is easy to prove.
\begin{lemma}
Let $\xi:  A \to ]0,\infty$ be a character. Then for $f_1, f_2 \in C^\infty(A, \Vtau),$ 
$$ 
\xi^2\inp{f_1}{f_2}_H = \inp{\xi f_1}{\xi f_2}_H.
$$
\end{lemma}

If $f_1, f_2 \in C^\infty(\tau:G/N_0: \chi)$ and $H\in \fa$ then by using the isomorphism
$C^\infty(\tau:G/N_0: \chi) \simeq C^\infty(A, \Vtau)$ we define the function
$\inp{f_1}{f_2}_H: A \to \C$ as above. 
Note that by (\ref{e: defi f g sub H}) we have
\begin{equation}
\label{e: function round brackets}
(f_1, f_2)_{H} = d_0^2 \inp{f_1}{f_2}_H  = \inp{d_0f_1}{d_0f_2}_H.
\end{equation}
For $f_1, f_2 \in C^\infty(\tau: M_{1F}/(M_{1F}\cap N_0): \chi)$ we identify $f_1,f_2$ with 
functions in $C^\infty(A, \Vtau)$ and then,  
\begin{equation}
\label{e: new bracket function}
{}^*\!d^2\inp{f_1}{f_2}_H = \inp{{}^*\!df_1}{{}^*\!d f_2}_H.
\end{equation}
Here ${}^*\!d(a) = d_0(a)/d_F(a).$ 
Harish-Chandra  \cite[p. 211]{HCwhit} uses the notation $(f_1, f_2)_{H}$ for the function
in (\ref{e: new bracket function}), which he also used for the different 
function (\ref{e: defi f g H}). We tried to avoid the confusion
that may arise from this. 

Let $f_1, f_2 \in \cA(\tau: G/N_0: \chi)$, $\ga \in \gD,$ $F = F_\ga = \gD\setminus \{\ga\}.$ 
We write $f_{jF_\ga}$ for $f_{j P_F},$ for $j =1,2.$ 

\begin{lemma}
\label{l: nbhd contribution by constant term}
Let $U_\ga$ be a sufficiently small open neighborhood of $\exp H_\ga^0$ in 
$S_\ga.$ Then there exist $C, \gd > 0$ such that, for $t \geq 0,$ 
$$ 
\int_{M_t U_\ga} \left| (f_{1}, f_2)_{H_\ga^0}(a)  -  \inp{{}^*\!d\,f_{1F_\ga}}{{}^*\!d\, f_{2F_\ga}}_{H_\ga^0}(a) \right|  ds(a) 
\leq Ce^{-\gd t}.
$$ 
 \end{lemma}

\proof 
We write $F$ for $F_\ga = \gD\setminus\{\ga\}$ and define $R_j\in C^\infty(A,\Vtau),$ for $j=1,2$, by 
$R_j(a) = d_F(a) f_j(a) - f_{jF}(a).$ Then 
$$ 
(f_1, f_2)_{H_\ga^0} = \inp{d_0 \,f_1}{d_0 \,f_2}_{H_\ga^0} = \inp{{}^*d\,(f_{1F} + R_1)}{{}^*d\; (f_{2F} + R_2)}_{H_\ga^0},
$$
hence
\begin{eqnarray}
\lefteqn{\!\!\!\!\!\!
(f_1, f_2)_{H_\ga^0}- \inp{{}^*\!d\, f_{1F}}{{}^*d \, f_{2F}}_{H_\ga^0} =}\nonumber \\
&=&
\inp{{}^*\!d\, f_{1F}}{{}^*\!d\,R_2}_{H_\ga^0}+
\inp{{}^*\!d\,  R_1}{{}^*\!d\, f_{2F}}_{H_\ga^0}+
\inp{{}^*\!d\, R_1}{{}^*\!d\, R_2}_{H_\ga^0}.
\label{e: f12 min star f12}
\end{eqnarray}
By the theory of the constant term, there exists an open neighborhood 
$V$ of $H_\ga^0$ in $\fa,$ a constant $ \gd_1 > 0$ and for every  $X\in U(\fa)_1$ 
a constant $C_1 > 0$ such that for all $H \in V$ and $t \geq 0,$ one has, for $j =1,2,$
$$ 
|{}^*\!d(\exp tH) R_j(\exp tH) |\leq C_1 e^{-\gd_1 t}.
$$
Replacing $V$ by a smaller open neighborhood if necessary, we may arrange to have in addition
an estimate of the form 
$$
{}^*d(\exp tH) |f_{jF}(\exp tH)| \leq C_2 (1  + t)^N, \qquad (H \in V, t \geq 0),
$$ 
for $j =1,2.$ 
Combining these estimates with (\ref{e: f12 min star f12}) we find for a fixed $0< \gd_3 <  \gd_1$  that there exists a constant $C_3 > 0$ such that for all $H \in V$ and $t \geq 0$ we have the estimate 
$$
\left|(f_{1}, f_2)_{H_\ga^0}(\exp tH)  -  \inp{{}^*\!d\, f_{1F_\ga}}{{}^*\!d\,f_{2F_\ga}}_{H_\ga^0}(\exp tH) \right|  
\leq C_3 e^{-\gd_3 t}.
$$ 
Let now $U_\ga = \exp(V) \cap S_\ga,$ then by pulling back the integration over $M_t U$ 
 by $M_t,$ we find
$$ 
\int_{M_t U_\ga} \left| (f_{1}, f_2)_{H_\ga^0}(a)  -  \inp{{}^*\!d\,f_{1F_\ga}}{{}^*\!d\, f_{2F_\ga}}|_{H_\ga^0}(a) \right|  ds(a) 
\leq
 C_3 e^{-\gd_3 t} (1 + t)^{\ell -1}  \int_{U_\ga}  ds.
$$
The proof is now easily completed. 
\qed
Following Harish-Chandra, let $f,g \in \cA(\tau:G/N_0: \chi),$ and  suppose that
$f_F =0$ for $F \subset \gD$ with $|\gD\setminus F| \geq 2.$ For $\ga \in \gD$ 
we define the function $\{f,g\}_\ga: M_{1F_\ga} \to \C$ by
$$
\{f,g\}_\ga(m): = \inp{f_{F_\ga}(m\rdiff H_\ga^0)}{g_{F_\ga}(m)}
                                        -  \inp{f_{F_\ga}(m)}{g_{F_\ga}(m\rdiff H_\ga^0)}.
$$ 
For  $t \in \R$ we define the function 
$\{f,g\}_{\ga, t}: M_{F_\ga} \to \C$ by 
\begin{equation}
\label{e: defi f g ga t m}
\{f, g\}_{\ga, t} (m) = \{f,g\}_{\ga}(m \exp tH_\ga^0).
\end{equation}
Furthermore, we consider integral
\begin{equation}
\label{e: defi of formal J ga}
J_\ga(f,g,t):= \int_{M_{F_\ga}/(M_{F_\ga} \cap N_0)} \{f,g\}_{\ga, t}\; d\dot{m}.
\end{equation}

\begin{lemma}
\label{l: comparison with J ga}
Assume that $f_F = 0$ for $F \subset \gD$ such that $|\gD\setminus F|> 1.$ 
Let $\ga \in \gD.$ Then the integral in (\ref{e: defi of formal J ga}) converges 
absolutely. If $U_\ga$ is a sufficiently small neighborhood  of $\exp H_\ga^0$ in $S_\ga$ then
there exist constants $C>0, \gd > 0$ such that for all $t \geq 0,$ 
\begin{equation}
\label{e: comparison int with J ga}
\left| \int_{M_t U_\ga} (d^*)^2 \inp{f_{F_\ga}}{g_{F_\ga}}_{H_\ga^0 } \; ds(a)
- J_\ga(f,g,t)\right|  \leq C e^{-\gd t}.
\end{equation}
\end{lemma}

\proof
For $t \in \R$ we define the function $f_{\ga, t}: M_F \to \Vtau, m \mapsto f_{F_\ga}(m \exp t H_\ga^0).$ The function $g_{F_\ga, t}$ is defined in a similar way.  
Both of these functions behave finitely under $d/dt,$ hence can be expressed as
\begin{equation}
\label{e: tempered deco f}
f_{F_\ga, t} = \sum_{\eta\in \cE, \, 0\leq k\leq n } f_{F_\ga, \eta , k}\; t^k e^{t \eta}\qquad g_{F_{\ga, t}} = \sum_{\eta\in \cE,\, 0\leq k\leq n } g_{F_\ga, \eta, k}\; t^k e^{t \eta},
\end{equation}
where $\cE\subset i \R$ a finite subset and $n \in \N.$ Since the exponential polynomial functions
$t\mapsto t^k e^{\eta t}$ are linearly independent over $\C,$ the expressions
in (\ref{e: tempered deco f}) are unique. 
By temperedness  of $f$ and $g,$ the functions $f_{F_\ga\eta, k}$ and 
$g_{F_\ga, \eta, k}$ belong to $\cA(\tau: M_F/M_F\cap N_0: {}^*\chi).$ 
By transitivity of the constant term, each $f_{F_\ga, \eta, k}$ has constant term zero along every proper  standard  parabolic
subgroup of $M_F.$ This implies the existence of constants $\ge > 0$ and $C > 0$ 
such that for all $\eta \in \cE$ and $0\leq k \leq n,$ 
$$
| f_{F_\ga, \eta, k}(\exp {}^*H) \leq C e^{-{}^*\rho({}^* H) - \ge|{}^*H|}, \qquad ({}^*H \in {}^*\fa).
 $$ 
For $g$ there exist constants $C', N> 0$ such that for all $\eta \in \cE$ and $0\leq k \leq n$ 
we have the tempered estimates
$$ 
|g_{F_\ga, \eta, k}(\exp {}^*H)| \leq C' (1 + |{}^*H|)^N e^{-{}^*\rho({}^*H)}, \qquad ({}^*H \in {}^*\fa).
$$
It follows from the definitions that, for $m \in M_F$ and $t \in \R,$
$$
\{f, g\}_{\ga, t}(m)  = \inp{\frac{d}{dt} f_{F_\ga, t}(m)}{ g_{F_\ga, t}(m)}
-
\inp{f _{F_\ga, t}(m)}{\frac{d}{dt} g_{F_\ga, t}(m)}.
$$
From the estimates given above, we infer the existence of $C'' > 0$ and $\ge' >0$ such that
for all $t,$ 
$$
\{f, g\}_{\ga,t}(\exp {}^*H) \leq C'' (1+ t)^{2n} e^{-2 {}^* \rho ({}^*H) - \ge' |{}^*H|},\qquad ({}^*H \in {}^*\fa).
$$
This implies the estimates
$$
|\{f, g\}_{\ga,t}(m) | \leq C'' (1+ t)^{2n} e^{-2 {}^* \rho(H(m)) - \ge' |H(m)|},\qquad (m \in M_{F_\ga})
$$ 
so that the integral defining $J_\ga (f,g,t),$ see (\ref{e: defi of formal J ga}), converges absolutely, and  
$$
|J_\ga(f,g,t)| = \cO((1 + |t|)^{2n}) \qquad (t\in \R).
$$
Let now $U_\ga$ be a neighborhood of $H_\ga^0$ in $S_\ga.$ Then 
the set $V_\ga := \log U_\ga$ is a neighborhood of $H_\ga^0 $ 
in $\fs_\ga = \partial \fa[1] \cap (H_\ga^0 + {}^*\!\fa).$ Therefore,
$V_\ga = {}^*V_\ga + H_\ga^0,$ with ${}^*V_\ga:= V_\ga \cap \fa_{F_\ga}^\perp$ 
a neighborhood of $0$ in ${}^*\fa.$ It follows that 
$U_\ga = {}^*U_\ga \exp H_\fa^0$ with ${}^*U_\ga$ a neighborhood of $e$ 
in  ${}^*A.$ Hence,  $M_t U_\ga = M_t({}^*U_\ga) \exp (t H_\ga^0).$ 
The Euclidean measure $ds$ on $M_t U_\ga$ is the translate of the Euclidean measure
$d{}^*\!a$ on $M_t ({}^*U_\ga)$ by $\exp t H_\ga^0.$ Consequently,
the integral on the left of (\ref{e: comparison int with J ga}) may be rewritten as 
$$ 
\int_{M_t {}^*U_\ga} {}^*\!d(a)^2  \{f,g\}_\ga({}^*\!a \exp tH_\ga^0)\, d\,{}^*\!a
=
\int_{M_t {}^*U_\ga} {}^*\!d(a)^2  \{f,g\}_{\ga, t}({}^*\!a ) d\,{}^*\!a.
$$
In view of (\ref{e: defi f g ga t m}) the latter integral may be rewritten as 
\begin{equation}
\label{e: partial int for J ga} 
\int_{\cO_t} \{f,g\}_{\ga, t} (m) \;d\bar m,
\end{equation}
where $\cO_t$ is the image of $K_F M_t ({}^*U_\ga)$ in  $M_F/M_F \cap N_0.$ 
The difference of (\ref{e: partial int for J ga}) with $J_\ga(f, g, t)$ is the  integral
with $\cO_t$ replaced by its complement $\cO^c_t$in $M_F/M_F \cap N_0.$ 
To finish the proof, it  suffices to show that there exists a $\gd > 0$ such that 
$$ 
\int_{\cO^c_t} | \{f,g\}_{\ga, t} (m) | \;d\bar m = \cO(e^{-\gd t}). 
$$ 
Choose $r > 0$ such that the ball $B_r \subset {}^*\!\fa$ with center $0$ and radius $r> 0$  is
contained in $\log  {}^*\!U_\ga.$ Then $\cO_t$ contains the image of $K_F \exp (t B_r) $ in $M_{1F}/M_{1F} \cap N_0,$ 
hence its complement $\cO_t^c$ is contained in $K_F\exp( t B_r^c) N_0,$ so that
\begin{eqnarray*}
 \int_{\cO_t^c} |\{f,g\}_{\ga, t} (m)| d\bar m &\leq&
 \int_{\exp t B_r^c} | {}^*d(a)^2 \{f,g\}_{\ga, t} ({}^* a)|  d{}^*\!a \\
&=&  \int_{{}^*\fa\setminus t B_r} C'' (1+ t)^{2n} e^{- \ge' |{}^*H|} d{}^*\!H.
\end{eqnarray*}
Now fix  $0< \gd < \ge'/r;$ then by using polar coordinates one readily checks that 
there exists a constant $C > 0$ such that the latter integral is bounded by $C e^{-\gd t},$ 
for $t \geq 0.$ 
\qed
 
\begin{lemma}
\label{l: generic estimation of integrals}
With assumptions as in Lemma \ref{l: comparison with J ga}, let
$H_0$ be any point of $\partial \fa[1]$ different from the 
points $H_\ga^0,$ for $\ga \in \gD.$ Then there exists an open neighborhood $U$ 
of $\exp H_0$ in $\partial A[1]$ and constants $C, \gd >0$ such that
\begin{equation}\label{e: generic estimate of int}
\int_{M_t U}  \; \left| (f,g)_{H_\ga^0}(a) \right| \;\;ds(a) \leq C e^{-\gd t}.
\end{equation} 
\end{lemma}

\proof
We will show that there exists an open neighborhood $V$ of $H_0$ in $\fa$ and constants
$C_1, \gd_1 > 0$ such that for all $H \in V$ and $t \geq 1$ we have the estimate
\begin{equation}
\label{e: estimate integrand surface int}
|(f,g)_{H_\ga^0}(\exp (tH))| \leq C_1 e^{- \gd_1 t}.
\end{equation}
Before proving this estimate we will first show that it implies the required estimate (\ref{e: generic estimate of int}). Indeed, let $U = V\cap A[1].$ By pulling back under $M_t$ and applying substitution 
of variables we obtain
\begin{eqnarray*}
\int_{M_t U}  \; \left| (f,g)_{H_\ga^0}(a) \right| \;\;ds(a)
&= & 
\int_U  \left| (f,g)_{H_\ga^0}(M_t a) \right| t^{\ell -1} ds(a)\\
&\leq &
\int_U  C_1 e^{-\gd_1 t} t^{\ell -1} ds(a).
\end{eqnarray*}
From this the result follows for any $0 < \gd < \gd_1$ and for $C$ suitably chosen.

We now turn to the proof of (\ref{e: estimate integrand surface int}).
Let $H_0 \in \fa[1]$ and assume that $H_0 \neq H_\ga^0$ for all $\ga \in \gD.$
First we assume that $H_0 \notin \cl (\fa^+).$ 
By the argument of  \cite[Cor.~2.4]{Bunitemp} it follows that for any $r> 0$ there exists an open neighborhood $V$ of $H_0$ and a constant $C'>0$ such that for $\gf$ equal to one of the functions $f, L_{H_\ga^{0\vee}}f, g$ or $L_{H_\ga^{0\vee}}g,$ we have 
$$ 
|\gf (t H)| \leq C' e^{- rt/2}, \qquad (H \in V, t \geq 1). 
$$
In view of (\ref{e: defi f g sub H}) the above estimate implies
$$
|(f,g)_{H_\ga^0}| \leq  (C')^2 d_0(\exp tH)^{-2} e^{-rt} =  (C')^2 e^{ - r' t}
$$
with $r' < r - 2 \sup_{H \in V} \rho(H).$ 

Note that if $H_0 \in S_\mu$ then $\mu(H_0) < -1$ so that 
$H_0 \notin \cl (\fa^+)$ and we are in the setting just discussed.

We now assume that $H_0 \in \cl(\fa^+)$. As just noticed, $H_0\notin S_\mu$ 
so that $H_0 \in \cup_{\gg \in \gD} S_\gg.$ 
Let $F$ be the collection of $\ga \in \gD$ vanishing on $H_0$. 
Since $H_0 \neq H_\ga^0$ for every $\ga \in \gD,$ it follows that $|\gD
\setminus F| > 1, $ so that by assumption, $f_F = 0.$ It now follows by application of 
\cite[Lemma 3.8]{Bcterm}that also 
$[L_{H_\ga^{0\vee}}f]_F =0.$ Thus for $\gf$ equal to 
$f$ or $L_{H_\ga^{0\vee}}f$ it follows from the estimation
of $d_F \gf  - \gf_F$ in  \cite[Lemma 3.5]{Bcterm} that 
for  a sufficiently small open neighborhood $V$ of $H_0$ in $\fa$ there exist constants
$C, N, \eta > 0$ such that, for $H \in V$ and $t\geq 0,$ 
$$
|d_F \gf(tH) |\leq C(1 + t)^N e^{-{}^*\rho_F(t{}^*H)}e^{-\eta t}.
$$
Hence, for $0< \eta' < \eta $ there exists $C'  >0$ such that
$$
|\gf(tH)| \leq C'  e^{-t \rho(H)-t\eta'},\qquad  (H\in V , t' \geq 0).
$$

Combining this with the tempered estimates for $g$ and $L_{H_\ga^{*\vee}} g,$ 
and using (\ref{e: defi f g H})
we obtain (\ref{e: estimate integrand surface int}) with $0< \gd_1 < \eta'$ and a suitable
$C_1 > 0.$
\qed

\begin{thm}{\rm [Harish-Chandra]}
\label{t: integral asympt HC} 
Let $f,g \in \cA(\tau: G/N_0:\chi)$ and assume that $f_F= 0$ for each subset $F \subset \gD$
with $|\gD\setminus F|> 1.$ Then there exists $\gd > 0$ such that
$$ 
\int_{G[t]/N_0} [f, g](x) \; d\bar x = \sum_{\ga \in \gD} J_{\ga}(f,g,t) + \cO(e^{-\gd t}), \qquad (t \to \infty).
$$ 
\end{thm}
\proof
For each $\ga \in \gD$ let $U_\ga$ be an open subset of $S_\ga$ with the properties
of Lemma \ref{l: nbhd contribution by constant term}. Put $U_\mu = \emptyset.$ Then it follows from combining 
Lemmas \ref{l: second form Gauss}, \ref{l: nbhd contribution by constant term} and \ref{l: comparison with J ga},
that
\begin{eqnarray*}
\lefteqn{\!\!\!\!\!\!\!\!\!\!\!\! \int_{G[t]/N_0} [f, g](x) \; d\bar x - \sum_{\ga \in \gD} J_\ga(f,g, t) = }\\
&=& \sum_{\gg \in\widehat\gD} c_\gg \int_{M_t (S_\gg \setminus U_\gg)} (f_{1}, f_2)_{H_\ga^0}(a)\,
ds(a)  + \cO(e^{-\gd t}).
\end{eqnarray*}
We now consider the compact set $\cK=  \partial A[1] \setminus \cup_{\gg \in \widehat \gD} U_\gg.$
It follows by application of Lemma \ref{l: generic estimation of integrals} and compactness that there exist  constants
$C, \gd > 0$ such that the estimate of the lemma is valid with $U$ replaced by $\cK.$ 
This implies the existence of $\gd > 0$ such that
\begin{eqnarray*}
\lefteqn{\!\!\!\!\!\!\!\!\!\!
\sum_{\gg \in \gD} c_\gg \int_{M_t (S_\gg \setminus U_\gg)} (f_{1}, f_2)_{H_\ga^0}(a)\; ds(a)=}\\
&=& \sum_{\gg \in \gD} c_\gg \int_{M_t (S_\gg \cap \cK)} (f_{1}, f_2)_{H_\ga^0}(a)\; ds(a)= 
 \cO(e^{-\gd t}).
\end{eqnarray*} 

We will now show that in the basic setting 
Theorem \ref{t: integral asympt HC} 
implies the Maass-Selberg relations, see also  \cite[p.~206]{HCwhit}.
In the basic setting, $G$ has compact center. $P\in \cPst$ is a maximal 
standard parabolic subgroup $P$ of $G.$
Suppose $\gs \in \dMPds$ and $\psi \in \cA_{2,P,\gs}.$ Let $\gL\in {}^*\fh_{P\iC}^*$ be the infinitesimal character of $\gs.$ Then $\gL$ is real and regular in the sense that the inner products with the roots of ${}^*\fh_P$ in $\fm_{P\iC}$ are real and non-zero.
In addition we fix $\nu \in i\fa_P^*$ such that $\nu \neq 0,$ $\gL +\nu$ is regular.
Put $f = f_\nu = \Wh(P,\psi,\nu).$ 
The infinitesimal character of $\Ind_{\bar P}^G(\gs \otimes - \nu)$
is given by $Z\mapsto \gg(Z, \gL -\nu).$ From Definition 
\ref{d: defi j Q} one sees that 
$$
Z f = R_Z f  = \gg(Z, \gL -\nu)f , \qquad (Z \in \fZ).
$$
It follows that the Casimir $\Omega$ acts on $f$ by the real eigenvalue
$\inp{\gL}{\gL} + \inp{\nu}{\nu} - \inp{\rho_P}{\rho_P}.$ In turn this
implies that 
$$
[f,f] = 0;
$$
see (\ref{e: brackets with Omega}). 
Since $P$ is maximal, it follows from the discussion below Lemma \ref{l: Wh is twohol} that the constant term of $f_\nu$ along a standard parabolic subgroup $Q$ is zero if $Q$ is not maximal. From Theorem 
\ref{t: integral asympt HC} it now follows that
\begin{equation}
\label{e: sum J to zero}
\sum_{\ga \in \gD} J_\ga(f,f,t) \to 0 \qquad (t \to \infty).
\end{equation}
For $\ga \in \gD$ let $F_\ga = \gD\setminus\{\ga\}.$ 
We write $f_{F_\ga}$ for the constant term of $f$ along the standard parabolic subgroup $P_{F_\ga}$ whose split component is $\R H_\ga^0.$ 

According to Poposition \ref{p: MS for adjacent maximal} there are two possibilities, 
(a): $|W(\fa_P)|=1$ and (b): $|W(\fa_P)|=2. $

In case (a) there exist precisely two distinct roots $\ga_1, \ga_2 \in\gD$ 
for which $Q_j: = P_{F_{\ga_j}}\sim P.$ 
Moreover, $W(\fa_{Q_j}| \fa_P) = \{s_j\}$ and $s_2 = - s_1,$ so, for $m_j \in M_{F_{\ga_j}}, t \in \R,$
$$
f_{F_{\ga_1}}(m_1 \exp t H_\ga^0) =  e^{i\gl t} \gf_1(m_1) ,\quad  f_{F_{\ga_2}}(m_2 \exp t H_\ga^0) = e^{-i\gl t}\gf_2(m_2)
$$
with $\gl = - i \nu(H_\ga^0) \in \R\setminus\{0\}.$ 
It now follows that 
$$
J_{\ga_1} (f,f,t) + J_{\ga_2}(f,f,t) = - \gl^2\inp {\gf_1}{\gf_1} + \gl^2\inp {\gf_2}{\gf_2} \to 0,
$$
from which we conclude that $\|\gf_1\|^2 = \|\gf_2\|^2.$ 

In case (b) we have $|W(\fa_P)| = 2$ and there is precisely one simple root $\ga \in \gD$
such that $Q:= P_{F_\ga} \sim P.$ It follows that $W(\fa_Q|\fa_P)$ consists of two elements, $s$ and $-s.$ 
Moreover, the constant term of $f$ along $Q$ is of the form
$$
f_{F_\ga} (m \exp tH_\ga^0) = e^{i \gl t} \gf_1(m) + e^{- i \gl t} \gf_2(m) 
$$
for $m \in M_{F_\ga},$ $t \geq 0.$
It follows that
\begin{eqnarray*}
J_\ga(f,f, t) & = & i\gl \inp{e^{i \gl t}\gf_1 - e^{-i \gl t} \gf_2}{e^{i \gl t}\gf_1 +e^{-i \gl t}\gf_2} - \\
&& - i\gl \inp{e^{i \gl t}\gf_1 + e^{-i \gl t} \gf_2}{e^{i \gl t}\gf_1 -e^{-i \gl t}\gf_2}\\
&=& 2i\gl (\|\gf_1\|^2- |\gf_2\|^2).
\end{eqnarray*}
From $J_\ga(f,f,t) \to 0$ it now follows that $\|\gf_1\|^2 =\|\gf_2\|^2.$\qed

{\em Completion of the proof of Theorem \ref{t: MS for B}.\ } In view of Proposition \ref{p: MS for adjacent maximal}, we have now completed the proof
of the Maass-Selberg relations MSC(P) for the basic setting. According to Lemma \ref{l: basic setting MSC gives MSB}  this implies
the validity of MSB(P) for the basic setting. 
By Lemma \ref{l: basic case MS} this implies the validity of the Maass-Selberg relations for the $B$-matrix as 
formulated in Theorem \ref{t: MS for B}. \qed
\medno
{\em Proof of the Maass-Selberg relations  MSC(P).\ } 
By Proposition \ref{p: comp MS rels B and C} we now conclude the validity of all Maass-Selberg relations for the $C$-functions
as formulated in Lemma \ref{l: MSrels C P}. 
\qed

\medbreak

For $P\in \cP$ we define the meromorphic function
$\Eta(P, \bar P) : \faPdc\to \End(\cA_{2,P})$ by
$$ 
\Eta(P, \bar P, \nu)|_{\cA_{2,P,\gs}} = \eta(P, \bar P, \gs, \nu) \, {\rm id}|_{\cA_{2,P,\gs}}.
$$
We may now formulate the validity of the entire collection of Maass-Selberg relations 
for the $C$-functions  as follows.
\begin{thm}
\label{t: MS for C}
Let $P \in \cP$, $Q \in \cPst$ and suppose that $Q\sim P.$ Then for each $s \in W(\faQ|\faP),$
$$
C_{Q|P}(s, -\bar \nu)^*C_{Q|P}(s,  \nu) = \Eta (P, \bar P, - \nu)
$$
as an identity of meromorphic functions in $\nu \in \faPdc.$ 
\end{thm}

\proof
The expressions on both sides of the equation define meromorphic functions of $\nu \in \faPdc.$ 
Hence, it suffices to prove the identity
$$
C_{Q|P}(s, \nu)^*C_{Q|P}(s,  \nu) = \Eta (P, \bar P, - \nu)
$$
 for generic $\nu \in i\faPd.$ This identity is equivalent to the MSC(P) as formulated in Lemma \ref{l: MSrels C P},
 which were proven to be valid in the text preceding the theorem.
 \qed
 
%%%%%%%%%%%%%%%%%%%%%%%%%%%%%%%%%%%%%
\section{The normalized Whittaker integral}
\label{s: normalized Whittaker}
For $P\in \cPst$ (standard   is mandatory) we consider the meromorphic function
$\nu \mapsto C_{P|P}(1,\nu), \faPdc \to \End(\cA_{2, P}).$ 
We recall from Lemma \ref{l: C Q Q}
that for each $\gs \in \dMPds$ and all 
$T \in C^\infty(\tau: K/K_P:\gs_P)\otimes \ccH^{-\infty}_{\gs, \chi_P}$ we have 
$$
C_{P|P}(1,\nu)\psi_T = \psi_{(A(P,\bar P, \gs, -\nu) \otimes I)T}
$$
as meromorphic functions of $\nu \in \faPdc$ with values in $\End(\cA_{2,P,\gs}).$ 
We recall that, for $R > 0,$ 
$$
\faPd(P,R):= \{\nu \in \faPdc\mid \inp{\Re \nu}{\ga} > R\;\;(\forall \ga \in \gS(P))\}.
$$
For $Q \in \cP$ we define $\Pi_{\gS(Q), \R}(\faQd)$ to be the set of polynomial functions
$q \in P(\faQd)$ which can be written as a product of linear factors of the form
$ \inp{\ga}{\dotvar} -c $ with $\ga \in \gS(Q)$ and $c \in \R.$

The following result is due to Harish-Chandra.
\begin{lemma}
\label{l: estimate C PP}
For every $R \in \R$ there exist $C,N > 0$ and  a polynomial function 
$q \in \Pi_{\gS(P), \R}(\faPd)$  such that the meromorphic function $\nu \mapsto q(\nu) C_{P\mid P}(1, \nu)$ 
is regular on $\faPd(P,R)$ and such that  
$$ 
\|q(\nu) C_{P\mid P}(1, \nu)\|_{\rm op} \leq C (1 + \|\nu\|)^N, \qquad (\nu \in \faPd(P, R).
$$
\end{lemma}
\proof
It suffices to prove this for the restriction of $C_{P|P}(1, \nu)$ to $\cA_{2,P, \gs},$ for 
each representation of the finite set of  $\gs \in \dMPds$ for which $\cA_{2,P, \gs}\neq 0.$ 
Since $T \mapsto \psi_T$ is a linear isomorphism of finite dimensional spaces,
it suffices to prove a similar estimate for $A(P, \bar P,\gs, -\nu)$ restricted to 
the finite dimensional space  $C^\infty(\tau: K/K_P:\gs_P).$ By equivalence of
norms on the latter space, that estimate is a consequence of \cite[Cor. 1.4]{BSestc} which
in turn is a straightforward consequence of  \cite{VW}, see \cite{Wrrg2}.
\qed
 
\begin{lemma}
\label{l: estimate eta inv}
Let $\gs \in \dMPds.$
There exist constants $\ge, C, N > 0$ such that 
the meromorphic function $\nu \mapsto \eta(\bar P, P, \gs, \nu)^{-1}$ 
is regular on $\faPd(\ge)$ and  
$$
|\eta(\bar O, P, \gs, \nu)^{-1}| \leq C ( 1 + |\nu|)^N , \qquad (\nu \in \faPd(\ge)).
$$ 
\end{lemma}

\proof This result is due to Harish--Chandra for $\gs$ a representation of the discrete series
of $M_P.$ His notation for $\eta(\bar P, P, \gs, \nu)^{-1}$ is $\mu_{P,\gs}(\nu).$ 
Under the weaker assumption that $\gs$ is unitary with real infinitesimal character,
the same result is proven in \cite[p. 235]{VW}. 
\qed 

Recall the definition of $\Eta(\bar P, P, \nu)\in \End(\cA_{2,P})$ in the text preceding
Theorem \ref{t: MS for C}.

\begin{cor}
\label{c: estimate C PP inverse}
There exist a polynomial function $q\in \Pi_{\gS(P), \R}(\faPd)$
and constants $\ge, C, N > 0$ such that 
the meromorphic function $\nu \mapsto  q(\nu) C_{P|P}(1,\nu)^{-1}$ is regular on $\faPd(\ge)$ and  
$$
\|q(\nu) C_{P|P}(1, \nu)^{-1}\| _{\rm op} \leq C(1 +|\nu|)^N , \qquad (\nu \in \faPd(\ge)).
$$ 
\end{cor}

\proof
From Theorem 12.8 it follows that
$$
C_{P|P}(1,\nu)^{-1} = \Eta(P, \bar P, \nu)^{-1} C_{P|P}(1,- \bar \nu)
$$
as meromorphic functions of $\nu \in \faPdc.$ 
The result now follows from Lemmas \ref{l: estimate C PP} and \ref{l: estimate eta inv}.
\qed
 For $P \in \cPst$ we define the associated normalized Whittaker integral by 
$$
\Wh^\circ(P, \psi, \nu)(x) := \Wh(P, C_{P|P}(1,\nu)^{-1} \psi, \nu)(x),
$$ 
for $\psi \in \cA_{2,P},$ $\nu \in \faPdc, x \in G.$ 

Let $\ge, q$ as in Corollary \ref{c: estimate C PP inverse}, and let $\gs \in \dMPds.$ It is readily verified that for suitable $0< \ge' < \ge,$ $r >0$ and for $\psi \in \cA_{2,P,\gs}$ 
the function
\begin{equation}
\label{e: q times Wh circ}
\nu \mapsto q(\nu) \Wh^\circ(P, \psi,\nu) \in C^\infty(\tau:G/N_0:\chi)
\end{equation}
belongs to $\twohol(\gL_\gs, \faPd, \ge, r,\tau).$
Here $\gL_\gs$ denotes the infinitesimal character of $\gs.$ 

The constant term of (\ref{e: q times Wh circ}) along  a  standard  parabolic subgroup  $Q \in \cPst$ associated with $P$ is given by 
$$ 
q(\nu) \Wh_Q^\circ(P, \psi, \nu, ma)= q(\nu) \sum_{s \in W(\faQ| \faP)} a^{s\nu} [C^\circ_{Q|P}(s, \nu)\psi](m)(a),
$$ 
for $\nu \in i \faPd,$ $m \in M_Q,\;$ $a \in A_Q.$ Here 
$$
C^\circ_{Q|P}(s, \nu) := C_{Q|P}(s, \nu)C_{P|P}(1 , \nu)^{-1}
$$
are meromorphic $\Hom(\cA_{2,P}, \cA_{2,Q})$-valued functions of $\nu \in \faPdc.$ The following result is an
important 
manifestation of the Maass-Selberg relations.
\begin{lemma}
\label{l: MS for normalized C}
Let $P,Q \in \cPst.$ Then for all $s \in W(\faQ|\faP),$ 
$$
C_{Q|P}^\circ(s ,-\bar\nu)^* C_{Q|P}^\circ(s,\nu) = {\rm id}_{\cA_{2, P}}
$$
as meromorphic functions of $\nu\in \faPdc.$ 
\end{lemma}

\proof
This follows from Theorem \ref{t: MS for C}.
\qed

\begin{lemma}
\label{l: est norm Ccirc}
\begin{enumerate}
\itema
There exists a constant $\ge>0$ such that $\nu \mapsto C_{Q|P}^\circ(s,\nu)$ is a holomorphic
$\Hom(\cA_{2,P}, \cA_{2,Q})$-valued function on $\faPd(\ge).$ 
\itemb
The constant $\ge>0$ can be chosen such that there exist $C,N > 0$ such that
$$
\|C_{Q|P}^\circ(s,\nu)\| \leq C (1+|\nu|)^N,\qquad (\nu \in \faPd(\ge)).
$$
\itemc 
The constant $\ge > 0$ can be chosen such that for all $u\in U(\faPd)$
there exist constants $C_u, N_u >0$ such that
$$
\|C_{Q|P}^\circ(s,\nu\rdiff u)\|\leq C_u (1 +|\nu|)^{N_u}.
$$
\end{enumerate}
\end{lemma}

\proof
From Cor.~\ref{c: estimate C PP inverse} and \cite[Lemmas 10.1, 10.2]{Bcterm} it follows that there exists a $q \in \Pi_{\gS(P)}(\faPd)$ and a constant $\ge > 0$ such 
that $\Gamma: \nu \mapsto q(\nu) C^\circ_{Q|P}(s, \nu)$ 
is holomorphic on $\faPd(\ge)$  and satisfies the estimate 
\begin{equation}
\label{e: estimate func Gamma}
\|\Gamma(\nu)\|_{\rm op} \leq C (1 +|\nu|)^N\qquad (\nu \in \faPd(\ge)).
\end{equation}

Let $\ell: \nu \mapsto \inp{\ga}{\nu} - c$ be a linear factor of $q.$ Then $\ell^{-1}(0)$
consists of $\nu \in \faQdc$ such that $\Im \inp{\ga}{\nu} = \Im c$ and 
$\Re \inp{\nu}{\ga} = \Re c.$ If $\Re c \neq 0$ then for $0< \ge < |\Re c|(1 +|\ga|)^{-1}$
we have $\ell^{-1}(0) \cap \faPd(\ge) = \emptyset.$ Furthermore,
$$
|\ell(\nu)|^{-1} \leq (|\Re c| -\ge)^{-1} \qquad (\nu \in \faPd(\ge)).
$$ 
We may write 
$q=q_0 q_1$ with $q_0$ equal to the product of the linear  factors $\ell$ with
$\Re c = 0$ and with $q_1$ equal to the product of the remaining factors. 
By choosing $\ge>0 $ sufficiently small we may arrange  that $|q_1 |^{-1}$  is bounded from above on $\faPd(\ge).$ Then $q_0(\nu) C^\circ_{Q|P}(s, \nu)  = q_1\nu)^{-1} \Gamma (\nu) $
is holomorphic in $\nu \in \faPd(\ge)$ and we have an estimate like (13.2) with  $q_1^{-1}\Gamma$ in place of $\Gamma.$ Thus, we may as well assume that 
$q = q_0$ from the start.

From Lemma \ref{l: MS for normalized C} it follows that the Hilbert-Schmid norm  $\|C^\circ_{Q|P}(s,\nu)\|$
is bounded for $\nu \in i\faPd\setminus q^{-1}(0).$ The latter set is open and dense
in $i\faPd.$ 
Let $\ell $ be a linear factor of $q$, and let $\nu_0 \in \ell^{-1}(0) \cap i\faPd.$
There exists a sequence $\mu_j$ in $i\faPd\setminus q^{-1}(0)$ with limit $\nu_0.$
The sequence $\|C^\circ_{Q|P}(s,\mu_j)\|$ is bounded, hence $\Gamma(\mu_j) =
q(\mu_j) C^\circ_{Q|P}(s,\mu_j)$ tends to zero for $j \to \infty.$ It follows
that $\Gamma(\nu_0) = 0.$ Hence $\Gamma = 0$ on $\ell^{-1}(0) \cap i\faPd.$ The latter set is
a hyperplane in the real linear space $i \faPd.$ Furthermore. $\ell^{-1}(0) \cap \faPd(\ge)$ 
is a connected open part of the complex hyperplane $\ell^{-1}(0).$ By analytic
continuation it follows that $\Gamma = 0$ on $\ell^{-1}(0) \cap \faPd(\ge).$
We claim that this implies that $\ell^{-1}\Gamma$ extends to a holomorphic function
$\tilde \Gamma$ on $\faPd(\ge).$ Indeed, by choosing suitable (affine linear) coordinates $z_j$ 
on $\faPdc$ we may  arrange that $\ell = z_1.$ By using local power series expansions
we find that $z_1$ divides $\Gamma.$ 

By a straightforward application of Cauchy's integral formula we infer that
$\tilde \Gamma$ satisfies an estimate of type (\ref{e: estimate func Gamma}) with $\tilde \ge = \ge/2$ in place
of $\ge.$ Repeating this process we reduce $q$ to a non-zero constant, so that (a) and (b) are valid.
 
Finally,
(c) follows from (b) by an easy application of Cauchy's integral formula.
\qed

We observe that on account of Lemma \ref{l: est norm Ccirc} we have

\begin{lemma}
\label{l: multipl C on cS}
Let $P,Q \in \cPst$ be associated and suppose that $s \in W(\faQ|\faP).$ 
Then the map 
$$
\gf \mapsto s^{-1*}[ C_{Q|P}^\circ(s, \dotvar)\gf]
$$
 is continuous linear from $\cS(i\faPd, \cA_{2,P})$  to $\cS(i\faQd, \cA_{2, Q}).$
Here $s^{-1*}$ denotes pull-back under $s^{-1}: i \faPd \to i\faQd.$ 
\end{lemma}

\begin{cor}
\label{c: wh o two hol pr}
There exist $\ge > 0, r>0$  
such that for each $\gs \in \dMPds$ the function $\nu \mapsto 
Wh^\circ(P, \psi, \nu)$ belongs to $\twoholpr(\faP, \gL_\gs, \ge, r , \tau).$ 
\end{cor}
\proof This follows from \cite[Cor.\ 11.5]{Bcterm}, in view of 
Lemma \ref{l: est norm Ccirc}.
\qed

%%%%%%%%%%%%%%%%%%%%%%%%%%%%%%%%%%%%%%%%%%%% 
\section{Fourier transform and Wave packets}
\label{s: nFou and Wave}
Let $P \in \cPst$.
Since $\cA_{2,P}$ is the finite orthogonal direct sum of the finite dimensional non-zero
subspaces $\cA_{2, P, \gs}$ with $\gs \in \dMPds$ it follows from Cor. 13.6,
that the normalized Whittaker integral $\Wh^\circ(P, \psi)$ satisfies the uniformly tempered 
estimates of \cite[Thm. 16.2]{Bunitemp}. 

In analogy with the definition of the Fourier transform $\cF_P$ in \cite[\S 16]{Bunitemp}
we define the normalized Fourier transform $\cF^\circ_P: \cC(\tau:G/N_0:\chi)
\to C^0(i\faPd, \cA_{2, P})$ by 
$$ 
\inp{\cF^\circ_P(f)(\nu)}{\psi}  = \inp{f}{ \Wh^\circ(P, \psi, \nu)}_2
:=
\int_{G/N_0} \inp{f(x)}{\Wh^\circ(P,\psi, \nu)(x)}_\tau \; d\dot{x},
$$ 
for $\psi \in \cA_{2,P},$ $\nu \in i\faPd.$ 

\begin{thm}
The normalized Fourier transform $\cF^\circ$ defines a continuous linear map.
\begin{equation}\label{e: normalized cF P}
\cF^\circ_P: \cC(\tau:G/N_0:\chi)
\to \cS(i\faPd, \cA_{2, P}).
\end{equation}
\end{thm}

\proof 
This result is the analogue of  \cite[Thm. 16.6]{Bunitemp}. The proof is identical, provided
one uses the uniformly tempered estimates for the normalized Whittaker integral.
\qed

Later on it will be convenient to employ a characterization of the normalized Fourier transform
in terms of an integral kernel. 
For this point of view it is convenient to view the normalized Whittaker integral $\Wh^\circ(P, \dotvar, \nu),$ for $P \in \cPst$ and $\nu \in \faPdc$, 
as a function $G \to \Hom(\cAtwoP, \Vtau).$ Accordingly we write
$$
\Wh^\circ(P, \nu)(x) \psi := \Wh^\circ(P, \psi, \nu, x),
$$ 
for $x \in G$ and $\psi \in \cAtwoP.$ Note that $\nu \mapsto \Wh^\circ(P, \nu)$
may thus be viewed as a meromorphic function with values in the Fr\'echet space
$C^\infty(G, \Hom(\cAtwoP, \tau)).$ 

We adopt the similar point of view for the unnormalized Whittaker integral
$\Wh(P, \nu, x)$ and note that the two are related by
$$
\Wh^\circ(P, \nu, x) = \Wh(P, \nu, x) C_{P|P}(1, \nu)^{-1}
$$

We proceed to the promised characterization  of $\nFou_P$ with an integral kernel.
For $A \in \Hom(\cAtwoP, \Vtau)$ we denote by $A^*$ the Hermitian
adjoint in $\Hom(\Vtau, \cAtwoP)$ 
with respect to the given Hilbert structures on $\cAtwoP$ and $\Vtau.$
Next, we define the dual Whittaker integral
$\Wh^*(P, \nu)$ by 
$$ 
\Wh^*(P, \nu)(x) : = \Wh^\circ(P, - \bar \nu, x)^*, \qquad (\nu \in \faPdc, x \in G).
$$ 
We note that $\nu \mapsto \Wh^*(P, \nu)$ 
is a meromorphic function $\faPdc \to C^\infty(G , \Hom(\Vtau, \cAtwoP))$ 
satisfying the transformation laws
$$ 
\Wh^*(P, \nu, k x n) = \chi(n)  \Wh^*(P, \nu, x) \tau(k)^{-1},
$$
for $x\in G, k\in K$ and $n \in N_0.$ 
If $f \in \cC(\tau : G/N_0:\chi) $ we use the notation
$
\Wh^*(P, \nu)f
$
for the function $G/N_0 \to \cAtwoP$ defined by 
$$
\Wh^*(P, \nu)f: x \mapsto \Wh^*(P, \nu)(x) f(x).
$$ 
It is now readily checked that the normalized Fourier transform 
(\ref{e: normalized cF P})
is given by 
\begin{equation}
\label{e: charn normalized FT}
\nFou_P f (\nu )                                                                                                                                                                                                                                                                                                                                                                                                                                                                                                                                                                                                                                                                                                                                                                                                                                                                                                                                                                                                                                                                                                                                                                                                                                                                                                                                                                                                                                                                                                                                                                                                                                                                                                                                                                                                                                                                                                                                                                                                                                                                                                                                                                                                                                                                                                                                                                                                                                                                                                                                                                                                                                                                                                                                                                                                                                                                                                                                                                                                                                                                                                                                                                                                                                                                                                                                                                                                                                                                                                                                                                                                                                                                                                                                                                                                                                                                                                                                                                                                                                                                                                                                                                                                                                                                                                                                                                                                                                                                                                                                                                                                                                                                                                                                                                                                                                                                                                                                                                                                                                                                                                                                                                                                                                                                                                                                                                                                                                                                                                                                                                                                                                                                                                                                                                                                                                                                                                                                                                                                                                                                                                                                                                                                                                                                                                                                                                                                                                                                                                                                                                                                                                                                                                                                                                                                                                                                                                                                                                                                                                                                                                                                                                                                                                                                                                                                                                                                                                                                                                                                                                                                                                                                                                                                                                                                                                                                                                                                                 = \int_{G/N_0}  \Wh^*(P, \nu, x) f(x)\; dx,\qquad (\nu \in i\faPd).
\end{equation}

We retain the assumption that $P\in \cPst$ and denote
by $d\nu$ the Lebesgue measure on the real linear space $i\faPd$,
normalized in such a way that the usual Euclidean Fourier transform 
$\cF_e:  \cS(A_P) \to \cS(i\faPd)$ given by $\cF_e(f)(\nu)= \int_A f(a) a^{- \nu}\, da$
is an isometry for the obvious $L^2$-inner products on $\cS(A_P)$ and $\cS(i\faPd).$ 

\begin{defi}
The inverse  transform 
$\Wave_P : \cS(i\faPd,\cAtwoP) \to C^\infty(\tau: G/N_0:\chi),$ 
also called Wave packet transform, is 
defined by the formula
\begin{equation}
\label{e: wave packet formula}
\Wave_P \gf(x) = \int_{i \faPd} \Wh^\circ(P, \nu , x)\gf(\gl) \; d\gl,
\end{equation}
for $\gf \in \cA_{2,P}, x \in G.$ 
\end{defi}

We note that by the integral (\ref{e: wave packet formula}) is absolutely 
convergent and defines a smooth function of
$x$ in view of the uniformly tempered estimates for the normalized Whittaker integral.

\begin{thm}
The normalized Wave packet transform $\Wave_P,$ for $P \in \cP$, 
defines a continuous linear map 
$$ 
\Wave_P :  \cS(i \faPd, \cA_{2,P}) \to \cC(\tau: G/N_0 : \chi).
$$ 
\end{thm}
\proof
We fix $\gs \in \dMPds$ such that $\cA_{2,P,\gs} \neq 0.$ Let $\psi$ be an element
of the latter space. Then by linearity and finite dimensionality of $\cA_{2, P}$ 
it suffices to show that the  map 
$$
\gf \mapsto \int_{i\faPd} \gf(\nu) \Wh^\circ(P,\psi,\nu, x)\; dx
$$
is continuous linear  $\cS(i\faPd) \to \cC(\tau: G/N_0:\chi).$ Since $\Wh^\circ(P,\psi)$ is 
a family in $\twoholpr(\faPd, \gL_\gs, \ge, r, \tau),$ see Cor. \ref{c: wh o two hol pr}, this follows from \cite[Thm.\ 12.1]{Bunitemp}.
\qed

\begin{lemma} Let $P \in \cPst.$ The transforms $\nFou_P$ and $\Wave_P$ are 
conjugate in the sense that
$$ 
\inp{\nFou_P f}{\gf}_2 = \inp{f}{\Wave_P \gf}_2,
$$
for $f\in \cC(\tau: G/N_0:\chi)$ and $\gf \in \cS(i\faPd)\otimes \cAtwoP.$
\end{lemma}

\proof
The brackets on the left indicate the $L^2$-type inner product on $L^2(i\faPd, \cA_{2, P})$
and the brackets on the right indicated the inner product on $L^2(\tau: G/N_0: \chi).$
The inclusions $\cS(i\faPd, \cA_{2,P}) \to L^2(i\faPd, \cA_{2, P})$ and 
$\cC(\tau: G/N_0: \chi) \to L^2(\tau: G/N_0: \chi)$ are continuous linear. It follows that the pairings
$(f,\gf) \mapsto \inp{\nFou_P f}{\gf}$ and $(f,\gf)\mapsto  \inp{f}{\Wave_P \gf}$ are continuous sesquilinear
$\cC(\tau: G/N_0:\chi) \times \cS(i\faPd, \cA_{2,P}) \to \C.$ By density of $C_c^\infty(\tau: G/N_0:\chi)$
in $\cC(\tau: G/N_0:\chi)$ and of $ C^\infty_c(i\faPd, \cAtwoP) $ in $\cS(i\faPd, \cA_{2,P})$
it suffices to prove the identity for $f \in C_c^\infty(\tau:G/N_0:\chi)$ and $\gf \in C_c^\infty(i\faPd, \cA_{2, P}).$
For such $f$ and $\gf$ we have
\begin{eqnarray*}
\inp{\nFou_P f}{\gf}_2 &=& \int_{i\faPd} \inp{\int_{G/N_0} \Wh^*(P,\nu,x)f(x) dx}{ \gf(\nu)}\; d\nu\\
&=& \int_{i\faPd} \int_{G/N_0} \inp{\Wh^*(P,\nu,x)f(x)}{ \gf(\nu)}\; dx \;d\nu\\
&=& \int_{G/N_0}  \int_{i\faPd}\inp{f(x)}{ \Wh^\circ (P,\nu,x)\gf(\nu)} \;d\nu \; dx \\
&=&\inp{f}{\Wave_P\gf}_2
\end{eqnarray*}
\vspace{-24pt}

\qed
\medbreak

\section{The functional equations}
\label{s: functional eqn}
Based on the results of the previous section, we shall
now derive Harish-Chandra's functional equation,
see  \cite[\S 1.7]{HCwhit}. 
 
\begin{lemma}
\label{l: conclude identity}
Let $P,Q \in \cPst$ be associated and let $s \in W(\faQ\mid\faP).$ Then
\begin{equation}
\label{e: functional equation}				
\Wh^\circ(P, \nu) = \Wh^\circ(Q, s\nu)\,C^\circ_{Q\mid P}(s, \nu)
\end{equation}
as an identity of meromorphic functions of the variable $\nu \in \faPdc.$
\end{lemma}

\proof 
We will first establish the existence of a meromorphic function
$F: \faPdc \to \Hom(\cA_{2,P}, \cA_{2,Q})$ such that
\begin{equation}
\label{e: Wh circ P Q F}
\Wh^\circ(P, \nu) = \Wh^\circ(Q, s\nu) F(\nu), \qquad (\nu \in \faPdc).
\end{equation} 
Indeed, by Lemma \ref{l: action of s on Whittaker int} there exists an element $F_s \in \Hom(\cA_{2, P}, \cA_{2,Q})$ such that
$$
\Wh(P , \nu) = \Wh(s P s^{-1}, s \nu) F_s,\qquad (\nu \in \faPdc).
$$  
Since $s P s^{-1}$ and $Q$ have the same split component, it follows from
Lemma \ref{l: standard int ops on Whit int} that there exists a meromorphic
$G: \faQdc \to \End (\cA_{2,Q})$ such that 
$$ 
\Wh(s^{-1} P s, s\nu) = \Wh(Q, s\nu) G(\nu).
$$
Combining these assertions it follows that (\ref{e: Wh circ P Q F})   
is valid with everywhere $\Wh$ in place of $\Wh^\circ$ and with $G (\nu)  F_s$ 
in place of  $F(\nu).$ 
If we combine this observation with the definition of the normalized Whittaker
integrals, we find that (\ref{e: Wh circ P Q F}) is valid with 
$$ 
F(\nu) := C^\circ_{P\mid P}(1, \nu)^{-1}\after G(\nu) \after F_s \after C^\circ_{Q \mid Q}(1, s \nu),\qquad
 (\nu \in \faPdc).
$$ 
Let $\Omega\subset i\faQd$ be the set of all 
$\nu \in i \faQd$ which are a regular point of  all meromorphic functions
in the above expression for $F,$ and for which the elements $v\nu,$ for $v\in W(\fa_Q),$
are mutually distinct. Then $\Omega$ is open dense in $i\faQd.$ 
For $\nu \in \Omega,$ the functions on the left and right of
(\ref{e: Wh circ P Q F}) are tempered and belong to $\cA(\tau: G/N_0:\chi).$ By taking the constant terms of these functions along $Q$ and comparing 
the appearing exponential functions with exponent
$s\nu,$ we find that $F(\nu) = C^\circ_{Q\mid P}(s, \nu)$ for all $\nu \in \Omega.$
By analytic continuation this identity is valid as an identity of meromorphic
functions of $\nu \in \faQdc,$ and (\ref{e: functional equation}) follows.
\qed

\begin{cor}
Let $P,Q,R\in \cPst$  all be associated to each other.
For all  $s \in W(\faQ\mid \faP)$ and $t \in W(\faR\mid \faQ),$  
\begin{equation}
\label{e: func eq C func}
C_{R \mid P}^\circ(ts ,  \nu)  = C_{R \mid Q }^\circ(t, s\nu) \, C_{Q \mid P}^\circ(s , \nu)
\end{equation}
as an identity of meromorphic functions of $\nu \in \faPdc.$
\end{cor}
\proof
Let $\nu \in i\faPd$ be a regular point for each of the three meromorphic functions
appearing in equation (\ref{e: functional equation}). Then (\ref{e: func eq C func}) follows by taking the constant terms along $R$ of the tempered functions on both sides of  (\ref{e: functional equation}). The proof is finished by application of analytic continuation.
\qed

The functional equations have important consequences for the 
Fourier and Wave packet transforms.

\begin{cor}
\label{c: FE nFou}
Let $P ,Q \in \cPst$ be associated and suppose that $s \in W(\faQ|\faP).$ 
Then for all $f \in \cC(\tau:G/N_0 : \chi),$
\begin{equation}
\label{e: functional transformation FT}
C^\circ_{Q \mid P}(s, \nu) \nFou_P f(\nu) =  \nFou_Q f(s\nu),\qquad (\nu \in i\faPd)
\end{equation}
\end{cor}

\proof
From Lemma \ref{l: conclude identity} we deduce that, for $\nu \in i\faPd,$ 
$$
\Wh^*(P, \nu) = C^\circ_{Q\mid P}(s, -\bar \nu)^*\, \Wh^*(Q, s\nu) =
C^\circ_{Q\mid P}(s, \nu)^{-1} \Wh^*(Q, s\nu).                                                                                                                         
$$
This implies that
$$ 
C^\circ_{Q \mid P}(s, \nu)\Wh^*(P, \nu) =\, \Wh^*(Q, s\nu).
$$
The identity (\ref{e: functional transformation FT}) now follows in view of the characterization  of the normalized Fourier transform in (\ref{e: charn normalized FT}).
\qed

\begin{cor}
\label{c: FE Wave}
Let $P ,Q \in \cPst$ be associated and suppose that $s \in W(\faQ|\faP).$ 
Then for all $\gf \in \cS(i\faPd, \cA_{2, P}),$
$$
\Wave_P \gf = \Wave_Q s^{-1*}[C^\circ_{Q|P}(s: \dotvar) \gf].
$$
Here $s^{-1*}$ denotes pull-back by $s^{-1} : i\faQd \to i\faPd.$
\end{cor}

\proof
By using the definition of $\Wave_P$ and Lemma \ref{l: conclude identity}, taking into account Lemma \ref{l: multipl C on cS} we find
\begin{eqnarray*}
\Wave_P \gf & = & \int_{i\faPd} \Wh^\circ(Q, s \nu) C^\circ_{Q|P}(s, \nu)\gf(\nu) \; d\nu\\
&= &
\int_{i\faQd} \Wh^\circ(Q,  \nu) s^{-1*}[ C^\circ_{Q|P}(s, \dotvar)\gf](\nu)\; d\nu.
\end{eqnarray*}
Now use the definition of $\Wave_Q.$ 

\qed
\begin{cor}
Let $P,Q\in \cPst.$ Then 
$$
\Wave_P \nFou_P f = \Wave_Q \nFou_Q f, \qquad (f \in \cC(\tau: G/N_0: \chi)).
$$
\end{cor}

\proof
By Corollaries  \ref{c: FE Wave} and \ref{c: FE nFou} we have 
$$
\Wave_P \nFou_P f =
\Wave_Q s^{-1*}[ C_{Q|P}^\circ (s, \dotvar)
\nFou_P f] = \Wave_Q s^{-1*}[s^*(\nFou_Q f)] = \Wave_Q \nFou_Q f.
$$

\qed

\section{Appendix: criterion of smoothness for distributions}
\label{a: smoothness of distributions}
In this appendix, we will prove a result needed in Section \ref{s: smoothness J}.
Let $U$ be
a smooth manifold of dimension $n.$ For $\cK \subset U$ a compact
subset and $u \in \cD'(U) = C_c^\infty(U)'$ we denote by $u_\cK$  the restriction
of $u$ to $C^\infty_\cK(\R^n).$ We will say that $u$ has order $r$ on $\cK$ if 
$u_\cK$ extends to a continuous linear functional on $C_\cK^r(U).$

We will write $C_\cK^\infty(U)$ for the Fr\'echet space of smooth functions 
$\gf \in C^\infty(U)$ with support contained in $\cK.$ For $p \in \N$  we have
$C^\infty_\cK (U) \subset C^p_\cK(U).$ The latter space has a natural Banach topology,
for which the inclusion is continuous with dense image. Transposition induces a
natural inclusion of the dual spaces equipped with the strong dual topologies,
$$ 
C^p_\cK (U)' \subset C^\infty_\cK(U)'.
$$ 
The topology on the first of these spaces is the dual Banach topology.

Let $\Omega$ be an open subset of $\R^q$ and suppose a map 
 $
T: \Omega \to C_c^{\infty}(U)', \;\;  y \mapsto T_y =T(y)
$
is given. 
For a compact
subset $\cK \subset U$ we define 
\begin{equation}
\label{T sub K}
T_\cK: \Omega \to C_\cK^\infty(U)'
\end{equation}
by
$$
T_\cK(y):= T_y|_{C_\cK^\infty (U)}, \qquad (x \in U).
$$ 
We say that $T_\cK$ maps to the Banach space $C_\cK^p(U)'$
if $T_\cK(y) \in C^p_\cK(U)'$ for all $y \in \Omega.$ 

Let $v_1, \ldots , v_q$ be a collection of smooth vector fields on $U,$ 
such that for each $x \in U$ the vectors $v_j(x)$ span the tangent space
$T_x U.$ If $\ga \in \N^q,$ we denote by $v^\ga$ the differential operator
$ v_1^{\ga_1} \cdots v_q^{\ga_q}$  on $U.$

\begin{thm}
Let $\Omega \subset \R^n$ be open and let $T: \Omega \to C^\infty_c(U)'$  be a map. 

If for each $\ga \in \N^q$ and every compact set $\cK \subset U$ the map
\begin{equation}
\label{e: condition T cK}
y \mapsto [v^\ga(T_y)]_\cK
\end{equation}
is continuous from $\Omega $ to the Banach space $C_\cK^p(U)'$, 
then
\begin{enumerate}
\itema
 for every $y \in \Omega$ the density $T_y$ is of the form $\tau_y dx$
 with $\tau_y \in C^\infty(U)$ and $dx$ a smooth positive density on $U;$
\itemb
for each $\ga \in \N^n$  
the function  $y \mapsto \partial^\ga (\tau_y),$ $ \Omega \to C^\infty(U),$ 
is continuous.
\end{enumerate} 
If $\Omega \subset \R^k \times \C^\ell$ and for every compact $\cK \subset U$ 
 the map (\ref{e: condition T cK}) 
is continuous from $\Omega $ to the Banach space $C_\cK^p(U)'$
and in addition holomorphic in the variable from $\C^\ell,$ then for each $\ga \in \N^q$  the function 
$y \mapsto v^\ga(\tau_y)$, $\Omega \to C^\infty(U),$  is continuous and in addition holomorphic in the variable from $\C^\ell.$ 
\end{thm}

\proof
If $\chi \in C_c^{\infty}(U)$ then by a straightforward application of the Leibniz rule
for differentiation, it follows that the map $\widetilde T: y \mapsto \chi T_y$ 
fulfills the hypotheses with $\widetilde T$ in place of $T.$ Clearly, it suffices to prove the
conclusions (a), (b) and the final assertion with $\widetilde T$ in place of $T,$ 
for any choice of $\chi.$
 
Therefore, we may as well assume that there exists a compact set $\cK_0$ contained in 
the interior of a compact subset $\cK \subset U$ such that $\supp \,T_y \subset \cK_0$ 
for all $y \in \Omega.$ Then the hypothesis that for every $\ga \in \N^n$ the 
function (8.2) is continuous is fulfilled for this particular $\cK.$ 
We will keep the sets $\cK_0$ and $\cK$ fixed from now on 
and fix a cut off function $\chi \in C_\cK^\infty(\R^n)$ which is $1$ on an open neighborhood
of $K_0.$ Then $\chi T_\circ = T_\circ$ for every distribution $T_\circ \in C^\infty(U)'$  
with support contained in $\cK_0.$ In particular this is true for $ T_\circ = T_y$  $(y \in \Omega).$ 

We denote by $\cF$ the Euclidean Fourier transform which on a function 
$f$ from the Schwartz space $\cS(\R^n)$ is given by 
$$ 
\cF f(\xi) = \int_{\R^n} f(x) e^{-i \xi \cdot x}\; dx, \quad (\xi \in \R^n).
$$ 
This operator is a topological linear isomorphism from $\cS(\R)$ onto itself.
The inverse transform is given by $\cS(\R^n) \ni g \mapsto S \cF (g),$ where $Sg(x) = g(-x).$ 
(No constant is appearing provided Lebesgue measure is replaced by a suitable
positive multiple.)

On the space $\cS'(\R^n)$ of tempered distributions, the Fourier transform is given by transposition of the transform on Schwartz functions, hence a topological linear isomorphism
from the space $\cS'(\R^n)$ onto itself; it is also denoted $\cF.$ This is compatible with
the notation for Schwartz functions if we embed $\cS(\R^n)$ into $\cS'(\R^n)$ 
by $f\mapsto f dx,$ where $dx$ denotes the standard smooth density on $\R^n.$ 
On the space of tempered distributions, the inverse Fourier transform is given by the transpose of $S \cF$. 

For every $y\in \Omega$ the distribution $T_y$ is compactly supported, with support in $\cK_0.$ 
Let $T_\circ\in C_c^{\infty}(U)'$ be any distribution with compact support contained
in $\cK_0$ and such that $(\partial^\ga T_\circ)_\cK$ belongs to $C_\cK^p(U)'$,
for every $\ga \in \N^n.$ 
It particular, such a $T_\circ$ is a tempered distribution and the associated Euclidean Fourier
transform is the tempered distribution given by the analytic function $\R^n \to \C$ defined
by  
$$ 
\cF(T_\circ): \xi \mapsto T_\circ(e^{-i\xi}),  \qquad(\xi \in \R^n). 
$$ 
Here $e^{-i\xi}$ denotes the exponentional function $x \mapsto e^{-i\xi\cdot x}, \R^n \to \C.$ 
By the inversion formula, 
\begin{equation}
\label{e: inverse Fourier}
T_\circ = \cF\after S ( \cF T_\circ). 
\end{equation}

We fix a norm $\|\dotvar \|_{\cK, r}$
which gives the Banach topology on $C_\cK^p (U),$ 
$$  
\|\gf\|_{\cK, p} := \max_{|\ga|\leq p} \sup_\cK |\partial^\ga \gf|,
$$ 
and we denote the dual norm on  $C_\cK^p(U)'$ by $\|\dotvar\|^*.$ 
One readily verifies that there exists a constant $c  > 0$ such that for every $\xi \in \R^n,$ 
$$
\|\chi e^{-i\xi}\|_{\cK, p} 
\leq c (1+\|\xi\|)^p,\qquad (\xi \in \R^n).
$$

Suppose now that for every $\ga \in \N^n$ the distribution $\partial^\ga(T_\circ)$ 
satisfies $(\partial^\ga T_\circ)_\cK \in C_\cK^p(U)'$.
Put $\gD := \sum_{j=1}^n \partial_j^2.$
Then for every $N \in \N$ and all $\xi \in \R^n,$ 
\begin{eqnarray*}
\label{e: estimate Fourier of partial ga}
|(1 + \|\xi\|^2)^N \cF(T_\circ)(\xi)| & = & |\cF((1- \gD)^N  T_\circ )(\xi)| \\
&=&  |\inp{(1- \gD)^N T_\circ }{\chi e^{- i\xi}}|\\ 
&\leq & \| (1- \gD)^N T_\circ\|^* \; \|\chi e^{-i\xi}\|_{\cK, p}
 \leq  C_N(T_\circ)  (1 + \|\xi\|)^p,
\end{eqnarray*}
where 
$$ 
C_N (T_\circ)  := c \|(1- \gD)^N T_\circ\|^*.
$$ 
This leads to the estimate
\begin{equation}
\label{e: estimate for cF T circle}
|\cF(T_\circ)(\xi)| \leq C_N(T_\circ)  (1 + \|\xi\|^2)^{p/2 - N}, \qquad (\xi \in \R^n)
\end{equation}
for every $N \in \N.$  
It follows from this estimate that the inverse  Fourier transform of $\cF (T)$ is the 
continuous density given by $T_\circ =  \tau_\circ dx$ where $dx$ is the  standard  smooth
density on $\R^n,$ and where 
\begin{equation}
\label{e: integral for T circ}
\tau_\circ(x) = S \cF \cF (T_\circ)(x) = \int_{\R^n} e^{i\xi \cdot x} \cF(T_\circ)(\xi) \; d\xi.
\end{equation} 

Using the estimate (\ref{e: estimate for cF T circle}) --  with $N$ sufficiently large -- for domination under the integral sign, we infer that $\tau_\circ$ is a smooth function
and for each $\gb \in \N^n, $ the derivative  
$\partial^\gb \tau_\circ$  is  given by differentiation under the integral sign. Substituting $T_y$ for $T_\circ$ in the resulting expression, and writing $T_y = \tau_y dx,$ 
we obtain
\begin{equation}
 \label{e: integral for partial T y}
\partial^\gb \tau_y(x) = \int_{\R^n} e^{i\xi \cdot x} (i\xi)^\gb \; \cF(T_y)(\xi) \; d\xi ,\quad (y \in \Omega).
\end{equation}

We note that $y\mapsto (T_y)_\cK$ is continuous $\Omega \to C_\cK^p(U)'$ 
by assumption. On the other hand, it is straightforward that
$\xi \mapsto \chi e^{- i\xi}$ is continuous  $\R^n \to C^p_\cK(U).$ Since the natural pairing
$C_\cK^p(U)' \times  C_\cK^p(U) \to \C$ is continuous bilinear it follows that
the map 
$$
\Omega \times \R^n \to \C,  (y, \xi) \mapsto  \inp{T_y}{\chi e^{-i\xi}} = \cF(T_y)(\xi)
$$ 
is continuous. It now easily follows that the integrand $I(y, x, \xi)$  of 
(\ref{e: integral for partial T y}) is a continuous function of $(y, x, \xi) \in \Omega \times U \times \R^n.$ 
On the other hand, the integrand is dominated by 
\begin{equation}
\label{e: dominating function}
C_N(T_y)  (1 + \|\xi\|)^{|\gb|}  (1 + \|\xi\|)^2)^{p/2 - N}
\end{equation}
while $C_N(T_y)$ is locally bounded in $y.$ Since this is valid for 
$N$ arbitrarily high, 
we conclude that the function $(y, x) \mapsto \partial^\gb \tau_y(x)$ belongs 
to $C(\Omega \times U).$ This in turn is equivalent to the assertion 
that $y \mapsto  \partial^\gb \tau_y$ is continuous from $\Omega$ to $C(U).$  As
$\beta\in \N^n$ was arbitrary we conclude that $y \mapsto \tau_y$ 
is continuous $\Omega \to C^\infty(U).$ 

We now turn to the statement about holomorphy. Write $y = (z , \gl)$ according to the decomposition 
$\R^q \simeq \R^k \times \C^{\ell}.$ 
It remains to prove the final assertion about holomorphy in the variable $\gl.$ For this
we first investigate the holomorphy of $\cF(T_\gl).$ 

By assumption the map $y\mapsto T_y$ is continuous $\Omega \to C_\cK^p(U)'$
and holomorphic in $\gl.$ This implies that $(y, \xi) \mapsto \inp{T_{y}}{e^{-i\xi}}$ 
is continuous $\Omega \times \R^n \to \C,$ with holomorphy in the variable $\gl.$ 
The integrand $I(y, x, \xi)$ introduced above is continuous, and holomorphic 
in $\gl,$ while it is dominated by (\ref{e: dominating function}). It is a well-known result that 
this implies that $(x,y) \mapsto \partial^\gb (\tau_y)(x)$ 
is holomorphic in the variable $\gl.$
It follows from this that the map $y \to \tau_y$ is continuous from $\Omega$ to $C^\infty(U)$ 
and holomorphic in the variable $\gl.$ 
\qed

\section{Appendix: divergence for a convex polyhedron}
\label{a: divergence}
This appendix gives a rigorous proof of Gauss' divergence theorem
for a compact convex polyhedral set in Euclidean space. For a more systematic 
treatment of Stokes' theorem on manifolds with singularities, we refer the reader to
\cite{Kstokes}.

By an affine hyperplane $\gs$ in $\R^n$ we mean a translate of a linear subspace
of codimension $1.$ Its complement $\R^n \setminus \gs$ is a disjoint
union of two open half-spaces. The closures of these are called the
closed half-spaces associated with $\gs.$ The latter can be retrieved
as the intersection of its closed half-spaces.

Any affine hyperplane can be described by a formula of the form
$\gs = \{ x \in \R^n\mid \xi(x) = c\},$ where $\xi \in \R^{n*}\setminus \{0\}$ 
and $c \in \R.$ The associated closed half-spaces are then described
by $\xi \leq c$ and by $\xi \geq c.$ 

In this section we assume that $C$ is a compact convex polyhedral subset
of $\R^n,$ i.e., a compact finite intersection of closed half-spaces.

Assume that $C$ has non-empty interior. Then by a hyperplane facet of
$C$ we mean an affine hyperplane $\gs$ such that $C$ is contained
in precisely one of the two closed half-spaces determined by $\gs$
and such that $C_\gs:= C \cap \gs$ has non-empty interior (denoted by
$C_\gs^\circ$) as a subset of $\gs.$ The sets $C_\gs$ and $C^\circ_\gs$ 
will be called the closed and open facets associated with $\gs.$ Note
that $\gs$ is the affine span of the open facet $C_\gs^\circ.$ The collection 
$\gS = \gS(C)$ of affine facets of $C$ is finite. Furthermore, it is readily verified
that 
$$ 
\partial C = \cup_{\gs \in \gS} \;\; C_\gs
$$ 
If $\gs_1$ and $\gs_2$ are distinct affine facets, then $C^\circ_{\gs_1} \cap C^\circ_{\gs_2} = \emptyset.$ 
If $\gs \in \gS$ then by $\nu_\gs$ we denote the unit vector in $\gs^\perp$ 
which points away from the closed half-space associated with $\gs$ that contains
$C$. We write $\nu: \partial C\to \R^n$  for the partially defined function 
determined by $\nu(s) =\nu_\gs$, for $s \in C_\gs^\circ.$ This function
is called the outward unit normal to $\partial C.$

If $f: U \to \R$ is a $C^1$-function on an open 
subset $U\subset \R^n$ and $H \in \R^n$ then 
we define the directional derivative $\partial_H f$ by 
$\partial_H f(x) = d/dt f(x + tH)|_{t=0},$  $(x \in U).$ 

\begin{lemma}
\label{l: Gauss for convex}
Let $f: C \to \R$ be a continuous function which is partially differentiable
on ${\rm int}(C)$ with partial deriatives that extend continuously to $C.$ Then for every $H \in \R^n,$
$$ 
\int_C \partial_H f(x) dx = \int_{\partial C} \inp{\nu(s)}{H} \; f(s)\; ds.
$$ 
\end{lemma}

\begin{rem} Let $\gS$ denote the set of affine facets of $\partial C.$ Then
$$
 \int_{\partial C} \inp{\nu(s)}{H} \; f(s)\; ds = \sum_{\gs \in \gS} 
\inp{ \nu_\gs }{H} \int_{C\cap \gs} f(s)\; ds.
$$ 
\end{rem}

The lemma will be proven in the rest of this section.
We start with investigating partial differentiation for integration
over $\R_+^n,$ where $\R_+ =] \,0,\infty\, [.$ 

\begin{lemma}
Let $f: \R^n_+ \to \R$ be a $C^1$-function, with partial derivatives up to order
$1$ extending continuously to $[0,\infty[^n$ and with support bounded in $\R^n.$
Then for every $H \in \R^n,$ 
$$ 
\int_{\R_+^n} \partial_H f(x)\; dx = \sum_{j=1}^n 
\int_{\R_+^{j-1} \times \R_+^{n-j}}   - H_j\;  f(x,0,y)\; dx dy.
$$
\end{lemma}

\proof
This is elementary from Fubini, and the formula $-f(x,0,y) = 
\int_0^\infty \partial_t f(x,t,y) dt.$
\qed

As a next step, we consider the cone $\Gamma$ in 
$\R^n$ generated by an $n$-tuple of distinct vectors $\gg_1, \ldots, \gg_n.$ 
We denote by $\partial_j \Gamma$ the cone generated by 
the points $\gg_k,$ $k \neq j.$ Then $\partial \Gamma$ is the union
of the cones $\partial_j \Gamma.$ The outward unit normal to 
$\partial_j \Gamma $ is denoted by $\nu_j.$ 

\begin{lemma}
Let $f: {\rm int} (\Gamma) \to \R$ be a $C^1$ function with partial
derivatives up to order $1$ that extend continuously to $\Gamma.$ In 
addition it is assumed that $f$ has bounded support.
Then 
$$ 
\int_\Gamma \partial_H f(x)\; dx = \sum_{j=1}^n \inp{\nu_j}{H}
\int_{\partial_j \Gamma} f(s)\; ds.
$$ 
\end{lemma}

\proof
The previous lemma is a special case. To obtain the more generala result 
from the present lemma,
let $T: \R^n \to \R^n $ be the linear map determined by 
$T e_j = \gamma_j,$ for $1 \leq j \leq n.$ Then $T(\R_+^n) = \Gamma.$ 
The boundary part $\partial_j \Gamma$ is the image under $T$ of the boundary 
part $\partial_j:= \partial_j \R_+^n.$ The outward normal vector $n_j$ at points of 
the interior of $\partial_j$ relative to the affine span of $\partial_j$ is related to 
$\nu_j$ as follows. 

A half-line $x + \R_+ H$ emanates from $\partial_j \Gamma$ in the outward 
direction if and only if $T^{-1} x + \R_+ T^{-1}H$ emanates from $\partial_j$ in the outward direction. From this it follows
 that for all $H \in \R^n\setminus \{0\},$ we have that 
$\inp{\nu_j}{H} > 0$ if and only if $\inp{n_j}{T^{-1}H} >0$. This implies that
there exists a constant  $c_j  > 0$ such that 
\begin{equation}
\label{e: nuj vs nj}
\nu_j = c_j \;T^{-1*}(n_j). 
\end{equation}
Since $\nu_j$ has unit length it follows that $c_j = \|T^{-1*}n_j\|^{-1}.$ 
 
 By linear substitution of variables and application of the previous lemma 
we obtain
\begin{eqnarray}
\lefteqn{\int_\Gamma \partial_H f(x) \; \; dx = 
 \int_{\R_+^n}  \partial_H f(T(z)) \; dz} \nonumber\\
&=& 
|\det T|\;  \int_{\R_+^n}  \partial_{T^{-1}H } \; [T^*f](z)\; dz\nonumber \\
&=&
\label{e: parametrized cone int}
| \det T|  \; \sum_{j=1}^n  \int_{\R_+^{j-1} \times \R_+^{n-j}} - (T^{-1}H)_j\; (T^*f)(x,0,y)  dx dy,
\end{eqnarray}
where $\hat x^j = (x_1\ldots x_{j-1} x_{j+1} \ldots x_n).$
We consider the $j-th$ term, and compare with the surface integral over $\partial_j \Gamma.$ Let $U_j:= \R_+^{j-1} \times \{0\} \times R_+^{n-j}.$ Then 
a regular parametrisation of $\partial_j \Gamma $ is given by the map
$
T_j: U_j \to \partial_j \Gamma, (x,y) \mapsto T(x, 0, y).
$
Now 
\begin{equation}
\label{e: definition of cone int}
\int_{\partial_j \Gamma }\inp{\nu_j}{H)} f(s) \; ds = 
\int_{U_j} \inp{\nu_j}{H}\; f(T(x,0,y))\; \; |J_j(T)| \;dx\, dy ,
\end{equation}
where $J_j(T) = Te_1 \times \cdots \times \widehat {T e_j} \times \cdots \times T e_n.$
By the definition of the exterior product, we have, for all $v \in \R^n,$ that
\begin{eqnarray*}
\inp{J_j(T)}{v} &= &\det (v, Te_1 ,\ldots, \widehat {T e_j} ,\ldots T e_n)\\
&=&
\det T \cdot  \det (T^{-1}v, e_1,\ldots, \widehat {e_j}, \ldots , e_n)\\
&= & (-1)^{j-1} \det T \cdot (T^{-1}v)_j = (-1)^{j} \cdot \det T\cdot  \inp{n_j}{T^{-1}v},
\end{eqnarray*}
from which we infer that
$
J_j(T) = (-1)^{j} \det T \cdot T^{-1*}n_j,
$
so that $c_j |J_j(T)|= |\det T|,$ see the line below (\ref{e: nuj vs nj}).
Using  (\ref{e: nuj vs nj}) it now follows that
$$ 
\inp{\nu_j}{H} \, |J_j(T)| =  c_j \inp{n_j }{T^{-1}H} |J_j(T)| \,  =  \inp{n_j }{T^{-1}H}|\det T|.
$$ 
Hence the integrals on the right of (\ref{e: parametrized cone int})
and of (\ref{e: definition of cone int}) are identical.
\qed

{\em Completion of the proof of Lemma \ref{l: Gauss for convex}.}
By an easy translation argument, Lemma \ref{l: Gauss for convex} is seen to be valid with domain 
$a +\Gamma$ $(a \in \R^n)$ in place of $\Gamma.$ 

We will now give the proof for $C$ an $n$-dimensonal simplex in $\R^n,$
i.e., the convex hull $\co(a_0, \ldots, a_n) $ of $n+1$ points in $\R^n$ 
whose affine span is $\R^n.$ For $0\leq j \leq 0$ we denote by $\Gamma_j$ 
the cone in $\R^n$ spanned by $a_i- a_j,$ for $i\neq j.$ 
Then $C \subset a_j + \Gamma_j.$ 

By a simple argument there exists a cover of $C$ by bounded open subsets 
$\cO_j \subset \R^n,$ for $0 \leq j \leq n,$ such that $a_j \in \cO_j$ and  
$\overline \cO_j \cap \co(a_0, \ldots, \widehat a_j , \ldots, a_n) = \emptyset.$
We fix a $C^1$-partition of unity $\{\psi_j\mid 0\leq j \leq 0\}$ 
over $C$ subordinate to the given cover. Thus, $\psi_j \in C_c^1(\cO_j),$
$0 \leq \psi_j \leq 1,$ and $\sum_{j=0}^n \psi_j = 1$ on an open neighborhood of 
$C.$ Let $f$ satisfy the hypothesis of Lemma \ref{l: Gauss for convex} and define, for
each $j,$ $f_j = \psi_j f.$ Then $f = \sum_{j=0}^n f_j $ and by linearity it suffices to prove 
the assertion of the lemma for each $f_j;$ fix $j.$

The function $f_j: {\rm int}(C) \to \R$ is continuous and its
partial derivatives up to order 1 extend continuously to $C.$ Moreover,
since $f_j = 0$ on an open neighborhood of $\co(a_0, \ldots, \widehat a_j, \ldots,
a_n)$ it follows that the extension of $f_j$ to $a_j + \Gamma_j,$
by requiring it to be zero outside $C$ has bounded support and its partial 
derivatives of order at most $1$ extend continuously to $a_j + \Gamma_j.$
By the first part of this proof, it follows that
$$
\int_{a_j + \Gamma_j} \partial_H f_j(x) \; dx =
\int_{a_j + \partial \Gamma_j}  \inp{\nu_{\Gamma_j}}{H} f(s) ds.
$$ 
The integrands of both integrals are zero on an open neighborhood
of $(a_j +\Gamma \setminus C)$. Therefore,
$$
 \int_{C} \partial_H f_j(x) \; dx =
\int_{C \cap (a_j + \partial \Gamma_j)}
  \inp{\nu_{\Gamma_j}}{H} f(s) ds.
$$ 
The domain of the latter integration equals
$\partial C \setminus \co(a_0, \ldots, \widehat a_j, \ldots, a_n).$
Since $f_j$ vanishes on $\co(a_0, \ldots, \widehat a_j, \ldots, a_n)$ 
it follows that the value of the latter integral remains unchanged if the domain is replaced by $\partial C.$ This completes the proof for $C$ a simplex.

Let now $C\subset \R^n$ be a compact convex polyhedral set
with non-empty interior and  fix a simplicial decompostion $\gS$ of $C$. Let
let $\gS_n$ denote the (finite) set of $n$-dimensional simplices $S\in \gS.$ 
Then $C = \cup \gS_n$ while all points of overlap are contained in 
the union of the simplices from $\gS_{n-1}.$ 
The integral of $\partial_H f$ over $C$ equals
\begin{eqnarray}
\sum_{S \in \gS_n}\int_{S} \partial_H f(x)\; dx &=& \sum_{S \in \gS_n}\int_{\partial S}
\inp{\nu_S}{H} f(s)ds \nonumber\\
& = & \sum_{S\in \gS_n} \sum_{\gs \in S_{n-1}(\partial S)} \int_{\gs}\inp{\nu_\gs}{H} f(s)ds.
\label{e: final sum of ints}
\end{eqnarray}
where $S_{n-1}[\partial S]$ denotes the set of $\gs \in S_{n-1}$ which
are contained in $\partial S$ (the appearing unit normal $\nu_\gs$ points out of $S$).
The double sum presenting the last integral can be rewritten as a sum 
of integrals over $\gs \in S_{n-1}.$ The elements of $S_{n-1}[\partial C]$ cover $\partial C$ with overlap contained in the negligable set $\cup S_{n-2}$. The remaining
elements, from $S_{n-1}\setminus S_{n-1}[\partial C],$ can be grouped 
in pairs of simplices in $\gS_{n-1}[S \cap S'] $, $(S,S'\in S_n)$ equipped with
opposite unit normals. 
As the contributions of these pairs cancel each other, the final sum in (\ref{e: final sum of ints})
can be rewritten as
$$
\sum_{\gs \in \gS_{n-1}[\partial C]}  \int_{\gs}\inp{\nu_\gs}{H} f(s)ds
=
 \int_{\partial C} \inp{\nu_\gs}{H} f(s)ds.
 $$
 \qed
 %------------------------------------------------------------
%\bibliographystyle{plain}
%\bibliography{MSrels.bib}

%Names -------------------------------------------------------------------
\def\adritem#1{\hbox{\small #1}}
\def\distance{\hbox{\hspace{3.5cm}}}
\def\apetail{@}
\def\addVdBan{\vbox{
\adritem{E.~P.~van den Ban}
\adritem{Mathematical Institute}
\adritem{Utrecht University}
%\adritem{PO Box 80 010}
%\adritem{3508 TA Utrecht}
\adritem{The Netherlands}
\adritem{E-mail: E.P.vandenBan{\apetail}uu.nl}
}

}
\mbox{}
\vfill
\hbox{\vbox{\addVdBan}\vbox{\distance}}
\end{document}